\documentclass{pspum-l}

\usepackage{amssymb}

\usepackage{graphicx}

\usepackage[cmtip,all]{xy}

\usepackage{amsmath, amssymb, fancyhdr, verbatim, stmaryrd}
\usepackage{enumerate}
\usepackage[usenames,dvipsnames]{xcolor}
\usepackage{mathrsfs}
\usepackage{tikz-cd}
\usepackage{framed, hyperref}
\usepackage[titletoc]{appendix}
\usepackage{bbm}
\usepackage{lipsum}
\usepackage{adjustbox}
\usepackage{devanagari}

\def\acts{\curvearrowright}
\def\qacts{\stackrel{?}{\acts}}

\DeclareFontFamily{U}{matha}{\hyphenchar\font45}
\DeclareFontShape{U}{matha}{m}{n}{
      <5> <6> <7> <8> <9> <10> gen * matha
      <10.95> matha10 <12> <14.4> <17.28> <20.74> <24.88> matha12
      }{}
\DeclareSymbolFont{matha}{U}{matha}{m}{n}
\DeclareFontFamily{U}{mathx}{\hyphenchar\font45}
\DeclareFontShape{U}{mathx}{m}{n}{
      <5> <6> <7> <8> <9> <10>
      <10.95> <12> <14.4> <17.28> <20.74> <24.88>
      mathx10
      }{}
\DeclareSymbolFont{mathx}{U}{mathx}{m}{n}

\DeclareMathSymbol{\obot}         {2}{matha}{"6B}

\RequirePackage{xspace}

\newcommand{\aA}{\ensuremath{\mathbf{A}}\xspace}

\newcommand{\BG}{\ensuremath{\mathbb{G}}\xspace}
\newcommand{\fg}{\ensuremath{\mathfrak{g}}\xspace}

\newcommand{\sH}{\ensuremath{\mathscr{H}}\xspace}

\newcommand{\tw}[1]{\langle #1 \rangle}
\newcommand{\cc}{\mathfrak{c}}

\newcommand{\F}{\mathbf{F}}
\newcommand{\RR}{\mathbf{R}}
\newcommand{\CC}{\mathbf{C}}

\newcommand{\wt}[1]{\widetilde{#1}}

\newcommand{\Q}{\mathbf{Q}}
\newcommand{\QQ}{\mathbf{Q}}
\newcommand{\Z}{\mathbf{Z}}
\newcommand{\ZZ}{\mathbf{Z}}

\newcommand{\mf}[1]{\mathfrak{#1}}

\newcommand{\Gal}{\operatorname{Gal}}

\newcommand{\R}{\mathbb{R}}

\newcommand{\ul}[1]{\underline{#1}}
\newcommand{\ol}[1]{\overline{#1}}

\newcommand{\Cal}[1]{\mathcal{#1}}
\newcommand{\A}{\mathbf{A}}
\DeclareMathOperator{\et}{\text{\'{e}t}}
\newcommand{\mbf}[1]{\mathbf{#1}}

\newcommand{\ft}{{}^{\tau}} 
\newcommand{\co}{\colon}
\newcommand{\mrm}[1]{\mathrm{#1}}
\newcommand{\msf}[1]{\mathsf{#1}}
\newcommand{\bs}{\backslash}

\newcommand{\PP}{\mathbf{P}}

\newcommand{\EE}{\mathbb{E}}
\newcommand{\DD}{\mathbb{D}}
\newcommand{\TT}{\mathbb{T}}

\newcommand{\inj}{\hookrightarrow}
\newcommand{\surj}{\twoheadrightarrow}
\newcommand{\dotimes}{\stackrel{\mbf{L}}\otimes}
\newcommand{\Zn}{{\mathbf{Z}/p^n\mathbf{Z}}}

\newcommand{\Zm}{{\mathbf{Z}/p^m\mathbf{Z}}}
\newcommand{\Zp}{{\mathbf{Z}_p}}
\newcommand{\cRHom}{\Cal{R}\Cal{H}om}

\DeclareMathOperator{\GL}{GL}
\DeclareMathOperator{\SL}{SL}
\DeclareMathOperator{\Frob}{Frob}

\DeclareMathOperator{\coker}{coker}

\DeclareMathOperator{\Tr}{Tr}

\DeclareMathOperator{\PGL}{PGL}

\DeclareMathOperator{\Hom}{Hom}

\DeclareMathOperator{\rank}{rank}

\DeclareMathOperator{\Aut}{Aut}
\DeclareMathOperator{\Rep}{Rep}
\DeclareMathOperator{\Mod}{Mod}

\DeclareMathOperator{\Spec}{Spec\,}
\DeclareMathOperator{\Spf}{Spf\,}

\DeclareMathOperator{\Lie}{Lie}
\DeclareMathOperator{\End}{End}

\DeclareMathOperator{\ad}{ad}

\DeclareMathOperator{\Frac}{Frac}

\DeclareMathOperator{\Stab}{Stab}
\DeclareMathOperator{\Bun}{Bun}
\DeclareMathOperator{\Ext}{Ext}
\DeclareMathOperator{\Pic}{Pic}

\DeclareMathOperator{\Id}{Id}
\DeclareMathOperator{\Ad}{Ad}

\DeclareMathOperator{\Sym}{Sym}

\DeclareMathOperator{\pt}{pt}

\DeclareMathOperator{\Sht}{Sht}

\DeclareMathOperator{\cla}{cl}

\DeclareMathOperator{\pr}{pr}

\DeclareMathOperator{\Hk}{Hk}

\DeclareMathOperator{\Ch}{CH}
\DeclareMathOperator{\CH}{CH}

\DeclareMathOperator{\Herm}{Herm}

\DeclareMathOperator{\Map}{Map}

\DeclareMathOperator{\RHom}{RHom}
\DeclareMathOperator{\Tor}{Tor}
\DeclareMathOperator{\res}{res}
\DeclareMathOperator{\obs}{obs}
\DeclareMathOperator{\Et}{\acute{E}t}

\DeclareMathOperator{\vir}{vir}

\DeclareMathOperator{\naive}{naive}

\DeclareMathOperator{\crys}{crys}
\DeclareMathOperator{\Betti}{Betti}
\DeclareMathOperator{\Hitch}{Hitch}

\DeclareMathOperator{\Tot}{Tot}

\newcommand{\bL}{\mathbf{L}}
\newcommand{\bT}{\mathbf{T}}

\newcommand{\limit}{\varprojlim}
\newcommand{\colim}{\varinjlim}

\newcommand\rightthreearrow{%
        \mathrel{\vcenter{\mathsurround0pt
                \ialign{##\crcr
                        \noalign{\nointerlineskip}$\rightarrow$\crcr
                        \noalign{\nointerlineskip}$\rightarrow$\crcr
                        \noalign{\nointerlineskip}$\rightarrow$\crcr
                }%
        }}%
}
\newcommand\leftthreearrow{%
        \mathrel{\vcenter{\mathsurround0pt
                \ialign{##\crcr
                        \noalign{\nointerlineskip}$\leftarrow$\crcr
                        \noalign{\nointerlineskip}$\leftarrow$\crcr
                        \noalign{\nointerlineskip}$\leftarrow$\crcr
                }%
        }}%
}

\DeclareMathOperator{\red}{red}
\DeclareMathOperator{\Sing}{Sing}
\DeclareMathOperator{\spec}{spec}

\newcommand{\DGHA}{\sH}

\DeclareMathOperator{\cInd}{c-Ind}

\newcommand{\isom}{\stackrel{\sim}{\to}}

\newcommand{\Ql}{\Q_{\ell}}

\newcommand{\htimes}{\stackrel{h}\times}
\renewcommand{\j}[1]{\langle{#1}\rangle}

\newcommand{\bu}{\bullet}

\newcommand\xr{\xrightarrow}
\DeclareMathOperator{\op}{op}

\renewcommand\a\alpha
\renewcommand\b\beta
\newcommand\g\gamma
\renewcommand\d\delta
\newcommand\D\Delta

\newcommand{\s}{\sigma}

\newcommand{\om}{\omega}

\newcommand{\da}{{\langle\alpha\rangle}}

\newcommand\cA{\mathcal{A}}
\newcommand\cB{\mathcal{B}}
\newcommand\cC{\mathcal{C}}
\newcommand\cD{\mathcal{D}}
\newcommand\cE{\mathcal{E}}
\newcommand\cF{\mathcal{F}}
\newcommand\cG{\mathcal{G}}
\newcommand\cH{\mathcal{H}}
\newcommand\cI{\mathcal{I}}
\newcommand\cJ{\mathcal{J}}
\newcommand\cK{\mathcal{K}}

\newcommand\cM{\mathcal{M}}

\newcommand\cO{\mathcal{O}}
\newcommand\cP{\mathcal{P}}

\newcommand\cR{\mathcal{R}}

\newcommand\cW{\mathcal{W}}
\newcommand\cX{\mathcal{X}}
\newcommand\cY{\mathcal{Y}}
\newcommand\cZ{\mathcal{Z}}

\newcommand\frP{\mathfrak{P}}

\newcommand\frS{\mathfrak{S}}

\newcommand{\fk}{\mathfrak{k}}
\newcommand{\fp}{\mathfrak{p}}

\newcommand\fra{\mathfrak{a}}

\newcommand\fm{\mathfrak{m}}

\newcommand\frp{\mathfrak{p}}

\newcommand{\isoarrow}{{~\overset\sim\longrightarrow~}}

\newcommand\dnv{\text{\dn{v}}}

\newcommand\ra{\rightarrow}

\newcommand{\chG}{\check{G}}

\newcommand\rA{\mathrm{A}}
\newcommand\rB{\mathrm{B}}
\newcommand\rH{\mathrm{H}}
\newcommand\rD{\mathrm{D}}
\newcommand\rM{\mathrm{M}}
\newcommand\rR{\mathrm{R}}
\newcommand\rS{\mathrm{S}}
\newcommand\rZ{\mathrm{Z}}

\newcommand{\sSet}{\msf{sSet}}

\newtheorem{thm}{Theorem}[subsection]
\newtheorem{lemma}[thm]{Lemma}
\newtheorem{prop}[thm]{Proposition}
\newtheorem{cor}[thm]{Corollary}

\newtheorem{conj}[thm]{Conjecture}
\newtheorem{hyp}[thm]{Hypothesis}
\newtheorem{defn-prop}[thm]{Definition-Proposition}

\theoremstyle{remark}
\newtheorem{fact}[thm]{Fact}
\newtheorem{remark}[thm]{Remark} 
\newtheorem{defn}[thm]{Definition}

\newtheorem{example}[thm]{Example}
\newtheorem{exercise}[thm]{Exercise}
\newtheorem{question}[thm]{Question}

\makeatletter
\def\th@remark{%
  \thm@headfont{\bfseries}%
  \normalfont 
  \thm@preskip \thm@preskip 
  \thm@postskip\thm@preskip
}

\numberwithin{equation}{subsection}

\begin{document}
\title[Derived structures in the Langlands correspondence]{Derived structures in the Langlands correspondence}


\author{Tony Feng}
\address{University of California Berkeley, Department of Mathematics, Berkeley, CA 94720, USA}
\curraddr{}
\email{fengt@berkeley.edu}
\thanks{TF was supported by NSF grants DMS-1902927, DMS-2302520, DMS-2441922, and a grant from the Simons Foundation.}

\author{Michael Harris}
\address{Columbia University, Department of Mathematics, New York, NY 10027}
\curraddr{}
\email{harris@math.columbia.edu}
\thanks{MH was partially supported by NSF grants
DMS-2302208 and DMS-2001369}


\date{}

\begin{abstract}
We survey several recent examples of derived structures emerging in connection with the Langlands correspondence. Cases studies include derived Galois deformation rings, derived Hecke algebras, derived Hitchin stacks, and derived special cycles. We also highlight some open problems that we expect to be important for future progress.
\end{abstract}

\maketitle

\tableofcontents

\section{Motivation and overview}

\subsection{What do we mean by ``derived structures in the Langlands correspondence''?}

The adjective ``derived'' has come to be applied to various constructions in mathematics, often with quite different meanings. For example, ``derived functors'', ``derived categories'', and ``derived algebraic geometry'' are some of the instances that we will encounter. Generally speaking, the word ``derived'' refers to an enhancement of mathematical constructions that incorporates homotopy theory. The past ten years have seen the rapid development of applications of homotopy theory to the Langlands correspondence, on several different fronts. Let us give an overview of those aspects which will be touched upon in this survey.  

\subsubsection{Derived Hecke algebras} Perhaps the first examples of ``derived functors'' that one encounters are the $\Tor$ and $\Ext$ functors, which are constructed by deriving $\otimes$ and $\Hom$. In the classical Langlands correspondence, a central role is played by the notion of Hecke operators, which form the Hecke algebra acting on automorphic forms. The Hecke algebra can be viewed as a certain space of endomorphisms; by replacing endomorphisms with ``derived endomorphisms'' (i.e., Ext groups) one obtains a notion of \emph{derived Hecke algebra}. This concept has arisen in two rather different contexts: the work of Ollivier-Schneider on the $p$-adic Langlands correspondence \cite{OS}, and the work of Venkatesh et al. on cohomology of arithmetic groups \cite{V19, GV, PV}.

\subsubsection{Derived moduli spaces} Derived functors in the above sense fall under the umbrella of classical homological algebra, where one derives functors on \emph{abelian} categories such as the category of modules over a ring. There is also a theory of derived functors on \emph{non-abelian} categories, which goes back to Quillen's \emph{homotopical algebra} \cite{Qu67}. When applied to the category of commutative rings, it enters the realm of ``derived algebraic geometry''. 

In derived algebraic geometry, one constructs derived enhancements of the usual constructs of algebraic geometry. For example, we shall discuss ``derived schemes'', ``derived stacks'', etc. One can imagine derived schemes as classical schemes enhanced with some additional homotopical information. 

It was recently realized that such derived enhancements provide a natural explanation of structures in the Langlands correspondence. Two different types of examples are the derived Galois deformation spaces of Galatius-Venkatesh \cite{GV}, which were used to explain the structure of cohomology of arithmetic groups, and the derived Hitchin stacks of Feng-Yun-Zhang \cite{FYZ2}, which were used to construct virtual fundamental cycles for special cycles related to the Kudla program. 

\subsubsection{Other directions, that will not be covered here}
Let us mention some further aspects of the Langlands correspondence where derived notions play a critical role, although they will \emph{not} be elaborated upon in these notes. 
\begin{itemize}
\item The Geometric Langlands correspondence features the moduli stack of local systems on a curve, and in particular the category of coherent sheaves on it. To obtain the ``correct'' version of this category, meaning the one which has a Langlands dual interpretation, one needs to understand this moduli stack in a derived way in general. Such matters are treated in Sam Raskin's article \cite{Ras} in these proceedings.
\item Derived commutative algebra plays an important technical role in recent work of Bhatt-Morrow-Scholze \cite{BMS2} and Bhatt-Scholze \cite{BS} on integral $p$-adic Hodge theory. 
\end{itemize}

\subsection{Why derive?}

Given the significant technical background required for derived algebraic geometry, it is natural to ask why and how it is useful for researchers interested in the more classical aspects of the Langlands correspondence, such as the reciprocity between automorphic forms and Galois representations.\footnote{We exclude branches like Geometric Langlands theory where derived algebraic geometry is baked into the very formulations of the core problems, hence relevant in an obvious way.} 

There are several independent reasons, some of which are touched upon in this survey. One is that higher cohomology (of locally symmetric spaces) has become central to the Langlands correspondence, because cohomology provides natural integral structures on the spaces of automorphic forms, and integral structures are needed to bring in $p$-adic methods. In many classical situations, the relevant cohomology groups are concentrated in a single degree, which allows one to ignore derived aspects to a large extent, but in general they are spread out over many degrees and it becomes necessary to work with complexes, derived functors, derived categories, etc. 

So far our discussion only requires deriving \emph{abelian structures}, which is a theory that has been developed and widely used since the 1970s. But in order to make the connection to Galois representations, we want to invoke moduli spaces of Galois representations called \emph{Galois deformation rings}. In order to make versions of these objects which are correspondingly ``spread out over many degrees'', it becomes necessary to derive the notion of a commutative ring. This is the starting point of derived algebraic geometry, and requires more sophisticated homotopical methods.

A third, independent thread treated in these notes concerns enumerative geometry related to automorphic forms. A celebrated example is the \emph{Kudla program}, which seeks to develop an incarnation of theta functions in arithmetic geometry. More generally, one would hope to develop an incarnation of \emph{relative Langlands duality}, in the sense of \cite{BZSV}, in arithmetic geometry. This means that we seek some incarnation of automorphic forms as cycle classes in the Chow groups or higher cohomology groups of Shimura varieties and moduli of shtukas. We explain here that such cycle classes should arise from derived geometry: the cycles in classical algebraic geometry are poorly behaved in general, while their derived versions have the desired properties. 

There are also other motivations which are \emph{not} treated in these notes, including the conjectural description of cohomology of moduli of shtukas explained in the articles of Xinwen Zhu \cite{Zhu21} and Emerton-Gee-Hellman \cite{EGH}. In certain situations discussed in those articles, one needs to incorporate derived structure in order to get the ``correct'' answers. 

\subsection{The style of these notes}

As the above summary is intended to convey, the past few years have witnessed a rapid and remarkably broad influx of homotopy theory into number theory, in the form of $\infty$-categories, derived algebraic geometry, stable homotopy theory, etc. 

The technical foundations for this theory are formidable, and complicated further by the multitude of different approaches that divide the literature. For example, the literature is split between the language of model categories and $\infty$-categories (the latter has acquired more popularity recently), and even within the framework of $\infty$-categories one finds different approaches (quasi-categories, DG-categories, simplicial categories). Derived algebraic geometry also comes in different flavors, being built out of commutative differential graded algebras, simplicial commutative rings, and $\EE_{\infty}$-algebras. These issues are significant, but often completely orthogonal to those encountered in applying the tools to number theory. 

Therefore, in view of the vast literature that already exists for homotopy theoretic background, we opted to write these notes as a kind of ``user's manual'' for how it is applied to number theory. We do not even attempt to give a thorough technical treatment of the foundational material. Occasionally, we will point out when different approaches exist and comment on their ``intuitive'' differences. Mostly, we will leave things somewhat vague and informal, relying on analogies and intuitions rather than formal definitions. 

\subsection{Outline of this article} In \S \ref{sec: DAG} we provide a moral introduction to derived algebraic geometry. We do not provide any technical details -- in fact, we make almost no mathematically precise statements; rather, we try to give some guide on ``how to think about'' certain words and phrases which come up in derived algebraic geometry. This is counterbalanced by Appendix \ref{appendix}, which focuses purely on the basic definitions of derived algebra. 

In \S \ref{sec: derived moduli} we discuss the notion of derived moduli spaces, which are our raison d'\^{e}tre for using derived algebraic geometry in the first place. We try to give a general pattern of heuristics which one uses in practice to recognize when derived moduli spaces should be relevant, and how to construct them. We illustrate these heuristics in several examples of relevance to the Langlands correspondence. 

In \S \ref{sec: derived special cycles} we survey in more detail the emergence and application of derived moduli spaces related to the Kudla program, focusing on the work of Feng-Yun-Zhang \cite{FYZ, FYZ2, FYZ3} on higher theta series for unitary groups over function fields. 

In \S \ref{DGdr} we survey the work of Galatius-Venkatesh \cite{GV} on derived Galois deformation rings in extremely informal fashion. 

In \S \ref{sec: DHA}, we introduce local derived Hecke algebras, touching on the $\ell \neq p$ results of \cite{V19} and also briefly the $\ell = p$ results of Schneider \cite{Sc15}, Ollivier, and Ronchetti. 

In \S \ref{sec: LSS} we discuss the cohomology of locally symmetric spaces and the global derived Hecke algebra. Then in \S \ref{venk} we explain Venkatesh's motivic conjectures.  

Finally, in \S \ref{sec: open problem} we formulate several open problems at the interface of derived algebraic geometry and the Langlands program, which should have important consequences. 

\subsection{Acknowledgments} We thank Henri Darmon, Dennis Gaitsgory, Soren Galatius, Chandrashekhar Khare, Barry Mazur, Arpon Raksit, Victor Rotger, Akshay Venkatesh, Jonathan Wang, Zhiwei Yun, and Wei Zhang for conversations about this material and for collaborations which are discussed here. We thank Gurbir Dhillon, Justin Wu, and the anonymous referees for comments and corrections on a draft. 

 A preliminary version of Appendix \ref{appendix} was tested on MIT graduate students in lectures by the first-named author in Spring 2022. We are grateful to the participants for their feedback.

\section{A guide to derived algebraic geometry}\label{sec: DAG}

\subsection{What is derived algebraic geometry?}

Algebraic geometry is the study of geometry built out of commutative rings.  \emph{Derived} algebraic geometry is the study of geometry built, in an analogous way, out of objects that we will informally call ``derived commutative rings'' for now. Derived commutative rings are some sort of enhancement of commutative rings; they govern what one might call ``derived affine schemes'', which form the local building blocks for \emph{derived schemes}. Therefore, the passage from algebraic geometry to derived algebraic geometry mostly amounts to the passage from commutative rings to derived commutative rings. 

\subsection{What are derived commutative rings?}
Let us caution at the start that there are multiple different approaches to the notion of ``derived commutative ring''. All of them are quite technical, and so we will not go into the formal definitions of any of them here. Instead, we will content ourselves with explaining intuitions, analogies, and how we think about things in practice.

\subsubsection{Simplicial commutative rings}For the bulk of \emph{these} notes, our notion of a ``derived commutative ring'' will be that of \emph{simplicial commutative ring}. The adjective ``simplicial'' can be thought of as synonymous with ``topological'', so it is a reasonable first intuition to think of a simplicial commutative ring as a topological commutative ring.\footnote{One reason we avoid working literally with topological commutative rings, which of course are well-defined mathematical concepts, is that we want to rule out ``pathological'' topological spaces. So a more refined slogan for ``simplicial'' might be ``topological and nice''.}

In particular, a simplicial commutative ring $\cR$ has associated (abelian) homotopy groups $\pi_0(\cR), \pi_1(\cR), \pi_2(\cR)$, etc. These should be thought of as the homotopy groups of the underlying ``topological space'' of $\cR$, but moreover the commutative ring structure on $\cR$ equips $\pi_*(\cR)$ with the structure of a graded-commutative ring. Understanding this graded-commutative ring is an approximation to ``understanding'' $\cR$, in the same sense that  understanding the usual topological homotopy groups of a topological space is an approximation to ``understanding'' the space.

The technical definition of a simplicial commutative ring is quite involved, and for this reason will be quarantined to Appendix \ref{appendix}. The main body of this text is written to be comprehensible with the notion of ``derived commutative rings'' treated as a black box, and it might be advisable to read it as such on a first pass. We also remark that we will shortly (in \S \ref{sssec: category CR SCR}) change our terminology for derived commutative rings from ``simplicial commutative rings'' to \emph{animated commutative rings}.

\subsubsection{Other models for derived commutative rings} We briefly mention other possible models for ``derived commutative rings'' that are commonly encountered in the literature. 
\begin{itemize}
\item \emph{Commutative differential graded algebras} (CDGAs). These are the most concrete to define and write down ``explicitly''. For example, it is relatively easy to write down CDGAs by generators and relations, while any such presentation would be too huge to write down for almost any simplicial commutative ring. The crucial flaw of CDGAs from our perspective is that they only work well in characteristic zero, while we definitely want to work integrally or in characteristic $p$. However, \emph{in characteristic $0$}, CDGAs could be used just as well as simplicial commutative rings. 
\item \emph{$\EE_{\infty}$-algebras}. Roughly speaking, an $\EE_{\infty}$-algebra is a ring whose multiplication is commutative up to homotopy coherence. The precise definition of $\EE_{\infty}$-algebras requires even more homotopy theory than that of simplicial commutative rings. In practice, they are typically more relevant for homotopy theory than algebraic geometry. It is true that $\EE_{\infty}$ lead to a functional theory in all characteristics and could be used just as well as simplicial commutative rings in characteristic zero, but away from characteristic zero they have different behavior which seems to be less relevant for us. 
\end{itemize}

\subsubsection{The category of derived commutative rings}\label{sssec: category CR SCR}

For the purposes of moduli theory, it is not the notion of a simplicial commutative ring itself which is essential, but the correct notion of its ambient category. In terms of the intuition for simplicial commutative rings as being like topological commutative rings, what we would like to do is to take the category of topological commutative rings \emph{up to homotopy equivalence}; informally speaking, we want to regard ``homotopic'' things as being equivalent.

 For this discussion, we would like to make an analogy to the following objects: abelian groups, complexes of abelian groups, and the derived category of abelian groups. In our analogy a commutative ring is parallel to an abelian group. To ``derive'' functors on the category of abelian groups, one considers ``simplicial abelian groups'' and the category of such is equivalent to the category of connective chain complexes. However, for many purposes (such as in defining derived functors), one is not really interested in the category of abelian groups but rather the \emph{derived category of abelian groups}. Informally, one can construct this category so that the objects are chain complexes, and then the morphisms are obtained from the morphisms among chain complexes by inverting quasi-isomorphisms. 

Simplicial commutative rings are parallel to chain complexes of abelian groups. The category we are really after, called \emph{the $\infty$-category of simplicial commutative rings}, is not the literal category of simplicial commutative rings but rather a category obtained from this by inverting ``weak equivalences'', which should be thought of as the analogue of quasi-isomorphisms between simplicial commutative rings.

Following terminology proposed by Clausen, and explained formally in \cite[\S 5.1]{CS}, we will refer to the desired category from the previous paragraph as the \emph{category of animated commutative rings}, and we will use ``animated commutative ring'' when thinking of objects of this category. Thus, an animated commutative ring is the same object as a simplicial commutative ring, but the terminology connotes a difference in the ambient \emph{category} implied, hence also in the notion of (iso)morphism between such objects. For example, isomorphisms of animated commutative rings correspond to what might be called ``quasi-isomorphisms'' of simplicial commutative rings, compared to a stricter notion of isomorphisms between simplicial commutative rings that is defined in the Appendix. 

We should remark that the concept of ``animated commutative ring'' is not new, but has traditionally just been treated under the name ``simplicial commutative rings'' in the literature. More recently, it has become more common (at least, in the areas that this survey touches upon) to use the terminology of animation to distinguish the $\infty$-categorical version. Animation works more generally for other algebraic structures as well, such as abelian groups, where it reproduces an equivalent notion to that of (connective) chain complexes.

\subsection{Visualizing derived schemes}\label{ssec: visualizing derived schemes} 

In algebraic geometry, one learns to think about commutative rings geometrically. For example, in introductory textbooks on scheme theory (such as \cite{Vakil}, whose \S 4.2 is the inspiration for the title of this subsection), one learns that non-reducedness of a ring -- which at first seems ``ungeometric'' because it is invisible at the level of field-valued points -- can be pictured geometrically as ``infinitesimal fuzz''. In other words, one visualizes $\Spec R$ as an infinitesimal thickening of $\Spec R_{\red}$. In this subsection, we will explain how one can similarly visualize a derived scheme as a type of infinitesimal thickening of a classical scheme. The guiding slogan is: 

\begin{quote}
The relationship between derived schemes and classical schemes is analogous to the relationship between (classical) schemes and reduced schemes. 
\end{quote}

Let us first explain the formal similarities. A scheme has an ``underlying reduced'' scheme, whose formation is functorial and defines a right adjoint functor to the inclusion of the category of reduced schemes into the category of all schemes. 

\begin{equation}\label{eq: RS S}
\begin{tikzcd}
\{\text{schemes}\} \ar[r, bend right, "Y \mapsto Y_{\red}"']  & \{\text{reduced schemes}\} \ar[l, bend right, "X \mapsfrom X"'] 
\end{tikzcd}
\end{equation}
At the level of rings, the formation of underlying reduced scheme corresponds locally to quotienting a ring by its nilradical, and this is left adjoint to the forgetful functor including the category of reduced commutative rings into  the category of all commutative rings.

\begin{equation}\label{eq: RCR CR}
\begin{tikzcd}
\{\text{commutative rings}\}  \ar[r, bend left, "R \mapsto R/\mrm{Nil}"]  & \{\text{reduced commutative rings}\} \ar[l, bend left, "S \mapsfrom S"]   
\end{tikzcd}
\end{equation}

Now let us explain the analogous picture for animated commutative rings. Any commutative ring $R$ can be viewed as an animated commutative ring $\ul{R}$, intuitively by ``equipping it with the discrete topology''. Then $\pi_0(\ul{R}) = R$, while $\pi_i(\ul{R}) = 0$ for $i>0$. This induces a fully faithful embedding of the category of commutative rings into the category of animated commutative rings, and in particular justifies the perspective that animated commutative rings form an ``enlargement'' of commutative rings. 

On the other hand, for any animated commutative ring $\cR$, its 0th homotopy group $\pi_0(\cR)$ has a natural commutative ring structure. These functors fit into an adjunction analogous to \eqref{eq: RCR CR}:
\begin{equation}\label{eq: CR SCR}
\begin{tikzcd}
\{\text{animated commutative rings}\} \ar[r, bend left, "R \mapsto \pi_0(\cR)"]  &  \{\text{commutative rings}\}   \ar[l, bend left, " \ul{\rR} \mapsfrom \rR"]  
\end{tikzcd}
\end{equation}

This construction glues, so that any derived scheme $Y$ has a \emph{classical truncation} $Y_{\cla}$, whose formation is functorial and defines a right adjoint to the functor from schemes to derived schemes:
\begin{equation}
\begin{tikzcd}
\{\text{derived schemes}\} \ar[r, bend right, "Y \mapsto Y_{\cla}"']  & \{\text{schemes}\} \ar[l, bend right, "X \mapsfrom X"'] 
\end{tikzcd}
\end{equation}
Analogously to \eqref{eq: RS S}, one can visualize a derived scheme as a classical scheme plus some ``derived infinitesimal fuzz''.\footnote{One difference in practice is that ``derived fuzz'' can carry \emph{negative (virtual) dimension}. This is analogous to how the natural extension of dimension to complexes is the \emph{Euler characteristic}, which can be negative.}

\begin{example}
We have not defined the \'{e}tale site of a derived scheme, but just as the \'{e}tale site of a scheme is isomorphic to that of its underlying reduced scheme, it turns out that the \'{e}tale site of a derived scheme is isomorphic to that of its classical truncation. This further supports the analogy between derived and non-reduced structure. 
\end{example}

Going forward we will use calligraphic letters such as $\cA,  \cB, \cR$ for simplicial commutative rings, and Roman letters such as $\rA, \rB, \rR$ for classical commutative rings. By the above remarks, we may regard classical commutative rings as animated commutative rings, but we still prefer to use this convention to emphasize when an animated commutative ring comes from a classical commutative ring. Similarly, we use roman letters like $\rM$ or $M$ for classical moduli spaces, and calligraphic letters like $\cM$ for derived moduli spaces. 

\subsection{Characteristics of simplicial commutative rings} We would like to make some remarks about how to express an animated commutative ring. In terms of the intuition for an animated commutative ring as a ``topological commutative ring up to homotopy'', it is natural to compare this question to that of expressing a topological space up to homotopy. A topological space is a kind of amorphous object, but algebraic topology furnishes two approaches to measuring it by algebraic invariants: homotopy groups, and (singular) homology groups. We shall explain two parallel invariants for ``measuring'' animated commutative rings: homotopy groups, and the cotangent complex. 

\subsubsection{Homotopy groups}\label{sssec: scr homotopy groups} An animated commutative ring $\cR$ has a sequence of homotopy groups $\pi_0(\cR), \pi_1(\cR), \ldots$, which are all abelian groups. Furthermore, the multiplication on $\cR$ equips 
\[
\pi_*(\cR) := \bigoplus_{i=0}^{\infty} \pi_i(\cR)
\]
with the additional structure of a graded-commutative ring. The construction is detailed in \ref{sssec: homotopy groups}. Under our intuitive analogy of an animated commutative ring being like a topological commutative ring, you can think of $\pi_*(\cR)$ as being the homotopy groups of the underlying topological space, with the ring structure induced by the multiplicative structure of $\cR$. 

\begin{example}
If $\cR$ comes from a topological commutative ring, then $\pi_*(\cR)$ agrees with the graded-commutative ring of homotopy groups (defined in the usual topological sense) of the topological commutative ring (with respect to the basepoint at the zero element).
\end{example}

For a topological space $X$, the homotopy groups of $X$ are a good measure of ``understanding'' $X$ itself. For example, if $X$ and $X'$ are ``nice'' topological spaces then a morphism $f \co X \rightarrow X'$ is a homotopy equivalence if and only if it induces a bijection on $\pi_i$ for all $i$, suitably interpreted to account for connected components. 

Similarly, we regard $\pi_i(\cR)$ as a good approximation to ``understanding'' an animated commutative ring $\cR$. In particular, a morphism of animated commutative rings $f \co \cR \rightarrow \cR'$ is an isomorphism if and only if it induces a bijection on $\pi_i$ for all $i$. 

\subsubsection{Cotangent complex}

Let $f \co A \rightarrow B$ be a morphism of classical commutative rings. Then there is a cotangent complex $\bL_f$, defined in \cite{Ill71, Ill72}, which was originally introduced for applications to deformation theory.  For our purposes, we regard $\bL_f$ as an animated $B$-module. We do not give the definition here because it requires a significant amount of development; Appendix \ref{appendix} gives a quick introduction to the cotangent complex and its applications. More generally, if $f$ is a morphism of schemes, or derived schemes, or even ``nice'' stacks, then it has a cotangent complex $\bL_f$.  

If $A = \Z$, we write $\bL_B$ for the cotangent complex of $B$ with respect to the unique map from $\Z$. More generally, for a scheme $X$ (or derived scheme, or ``nice'' stack) we write $\bL_X$ for the cotangent complex of the unique map $X \rightarrow \Spec \Z$.  

\begin{example}\label{ex: infinitesimal deformations}
Let $f \co A \rightarrow B$ be a morphism of commutative rings and $\wt{A} \rightarrow A$ a square-zero thickening with ideal $I$ (i.e., $I^2 = 0$). Let $J = I \otimes_A B$. We consider $\wt{A}$-algebra deformations of $B$, i.e. flat $\wt{A}$-algebras $\wt{B}$ such that $\wt{B} \otimes_{\wt{A}} \wt{A}/I \xrightarrow{\sim} B$. It is explained in \S \ref{ssec: deformation theory 1} that this deformation theory problem is ``controlled'' by $\bL_f$, in the sense that: 
\begin{itemize}
\item There is a class $\obs(I) \in \Ext^2_B(\bL_f, J)$ which vanishes if and only if such a deformation $\wt{B}$ exists. 
\item If $\obs(I)=0$, then the set of deformations $\wt{B}$ has the natural structure of a torsor for $\Ext^1_B(\bL_f, J)$.
\item The automorphism group of any deformation $\wt{B}$ is naturally isomorphic to $\Ext^0_B(\bL_f, J) \cong \Hom_B(\Omega_f, J)$. 
\end{itemize}
\end{example}

\begin{example}If $f$ is smooth, then $\bL_f$ has cohomology concentrated in degree $0$, i.e., is represented by a $B$-module, which is none other than the K\"ahler differentials $\Omega_{B/A}$. 
\end{example}

In fact, the construction of the cotangent complex uses simplicial commutative rings, and was one of the earliest applications for this theory. Roughly speaking, the idea is that when $B$ is not smooth over $A$ then one ``resolves'' $B$ by a smooth simplicial $A$-algebra, and forms the K\"{a}hler differentials of this resolution. This is reminiscent of how one defines derived functors in homological algebra, by replacing inputs with ``resolutions'' by certain types of chain complexes with good properties. While chain complexes only make sense for objects of abelian categories, simplicial objects make sense in much more generality. 

\begin{example}
Given a sequence of maps 
\[
\begin{tikzcd}
A \ar[r, "f"] \ar[rr, bend left, "h"]  &  B \ar[r, "g"] & C
\end{tikzcd}
\]
with $h = g \circ f$, we get an exact sequence 
\[
\Omega_f \otimes_B C  \rightarrow \Omega_h \rightarrow  \Omega_g \rightarrow 0.
\]
This will extend to an exact \emph{triangle} (in the derived category of $C$-modules) of cotangent complexes 
\[
\bL_f \dotimes_B C \rightarrow \bL_h \rightarrow \bL_g \xrightarrow{+1} 
\]
which in particular gives a long exact sequence in cohomology, 
\[
\ldots \rightarrow H^{-1}(\bL_h) \rightarrow H^{-1}(\bL_g) \rightarrow \Omega_f \otimes_B C  \rightarrow \Omega_h \rightarrow  \Omega_g \rightarrow 0.
\]
\end{example}

By definition, a morphism of animated commutative rings is an isomorphism if and only if it induces an isomorphism on homotopy groups. By contrast, it is a priori unclear ``how much'' information the cotangent complex carries. It turns out that it carries ``all'' the information of the derived structure. This is made precise by the following lemma of Lurie: 

\begin{lemma}[{\cite[Corollary 3.2.17]{Lur04}}]\label{lemm: Lurie Hurewicz} Suppose a morphism of simplicial commutative rings $f \co \cA \rightarrow \cB$ induces an isomorphism on $\pi_0$ and an isomorphism $\bL_{\cA} \otimes_{\cA} \cB \xrightarrow{\sim}  \bL_{\cB}$. Then $f$ is an isomorphism. 
\end{lemma}

We will not give the formal argument for Lemma \ref{lemm: Lurie Hurewicz}, but we will sketch the intuition. From a functor of points perspective, what we have to show is that $f$ induces an equivalence between morphisms to an arbitrary simplicial commutative ring from $\cA$ and from $\cB$. As discussed in \S \ref{ssec: visualizing derived schemes}, an arbitrary simplicial commutative ring can be thought of as a ``derived nilpotent thickening'' of its classical truncation. Since morphisms to a classical ring are controlled by $\pi_0$, the morphisms to any classical ring from $\cA$ and $\cB$ are identified by $f$. Then, by (a generalization of) Example \ref{ex: infinitesimal deformations}, morphisms to a derived nilpotent thickening of a classical ring are controlled by the respective cotangent complexes, which are again identified by $f$. 

\begin{remark}
We regard the cotangent complex as a ``homology theory'' for simplicial commutative rings. From this perspective, Lemma \ref{lemm: Lurie Hurewicz} is analogous to the Hurewicz theorem, which says that for ``nice'' spaces, it is equivalent for a map to induce a bijection on homotopy sets and homology groups. 
\end{remark}

\subsubsection{Tangent complex}\label{sssec: tangent complex}
Suppose $f \co \cA \rightarrow \cB$ is a morphism of animated commutative rings whose cotangent complex $\bL_f$ is \emph{perfect} as an animated $\cB$-module. This means that $\bL_f$ can be locally represented by a finite complex of free $\cB$-modules. Then we define the \emph{tangent complex} of $f$ to be $\bT_f = \cRHom_{\cB}(\bL_f, \cB)$, the (derived) dual of the cotangent complex. 

This construction can be globalized: if $f \co \cY \rightarrow \cX$ is a morphism of (derived) schemes (or nice stacks, so that $\bL_f$ exists) such that $\bL_f$ is perfect over $\cO_{\cY}$, then we define $\bT_f := \cRHom_{\cY}(\bL_f, \cO_{\cY})$.  

\subsection{Examples of derived schemes} Affine derived schemes are anti-equivalent to animated commutative rings. For an animated commutative ring $\cA$, we write $\Spec \cA$ for the corresponding affine derived scheme. General derived schemes are glued from affine derived schemes, analogously to how schemes are glued from affine schemes. We will now give some examples of natural constructions which lead to derived schemes in practice. 

\subsubsection{Derived fibered products} Even when starting out with purely classical schemes, derived schemes arise naturally through the \emph{derived fibered product} operation. 

The classical avatar of this operation is the usual fibered product of two morphisms of schemes $X \rightarrow Z$ and $Y \rightarrow Z$. The resulting fibered product $X \times_Z Y$ is characterized by a universal property among schemes with maps to $X$ and $Y$, such that the composite maps to $Z$ agree. It is built affine-locally by the tensor product operation on commutative rings. 

Similarly, the derived fibered product, which we also sometimes call the \emph{homotopy fibered product} $X \htimes_Z Y$, is characterized by an analogous universal property. It can be built affine-locally by the \emph{derived tensor product} operation on animated commutative rings (see below). More generally, the same considerations allow to construct the derived fibered product of derived schemes. 

Given morphisms of animated commutative rings $C \rightarrow A$ and $C \rightarrow B$, the derived tensor product $A \dotimes_C B$ can be characterized by its universal property as the coproduct of $A$ and $B$ in the category of animated $C$-algebras. An explicit model can be built in terms of simplicial commutative rings, by choosing ``resolutions'' of $A$ and $B$ as $C$-algebras; this is explained in \S \ref{ssec: derived tensor}. In particular, this explicit model shows that the homotopy groups of $A \dotimes_C B$ are the $\Tor$-groups, 
\begin{equation}\label{eq: derived tensor product}
\pi_i(A \dotimes_C B) \cong \Tor_i^C(A,B);
\end{equation}
although \eqref{eq: derived tensor product} does not illuminate the multiplicative structure on $\pi_*(A \dotimes_C B)$, so it really only encodes information about the underlying animated abelian group of $A \dotimes_C B$. From \eqref{eq: derived tensor product}, we see that the derived fibered product $X \htimes_Z Y$ of classical schemes is classical if and only if $\cO_X$ and $\cO_Y$ are Tor-independent over $\cO_Z$. In particular, this holds if either $X \rightarrow Z$ or $Y \rightarrow Z$ are flat. 

\begin{example}\label{ex: tangent complex of derived fibered product}
The derived fibered product interacts well with formation of (co)tangent complexes, which is actually a technical advantage of working in derived algebraic geometry. More generally, a slogan is that ``formation of tangent complexes commutes with homotopy limits''. In particular, for derived schemes $\cX \rightarrow \cZ \leftarrow \cY$ we have that 
\[
\bT_{\cX \times_{\cZ} \cY}  \cong \mrm{Fib} \left( \bT_{\cX}|_{\cX \times_{\cZ} \cY} \oplus \bT_{\cY}|_{\cX \times_{\cZ} \cY} \rightarrow \bT_{\cZ}|_{\cX \times_{\cZ} \cY} \right), 
\]
with the fiber formed in the derived category of quasicoherent sheaves on $\cX \times_{\cZ} \cY$ (this notation means the derived fibered product, as there is no ``naive fibered product'' for derived schemes).  
\end{example}

\subsubsection{Derived vector bundles}\label{sssec: derived vector bundles}
Let $X$ be a proper variety over a field $k$ and $\cE$ a locally free coherent sheaf on $X$. Consider the functor sending $S \in \mrm{Sch}_{/k}$ to the space of global sections of $\cE$ on $X_S$; it is represented by $\ul{H^0(X, \cE)}  = H^0(X, \cE) \otimes_k \A^1_k$, the $k$-vector space $H^0(X, \cE)$ regarded as an affine space over $k$. It can be described as $\Spec (\Sym H^0(X, \cE)^{*})$ where $H^0(X, \cE)^{*}$ is the dual of $H^0(X, \cE)$ over $k$.  

A perspective one learns in homological algebra is to view $H^0(X, \cE)$ as the zeroth cohomology group of a complex $\rR\Gamma(X, \cE)$, whose higher cohomology groups are $H^i(X, \cE)$ for $i \geq 1$. Analogously, there is a derived affine scheme 
\[
\ul{\rR\Gamma(X, \cE)} = ``\Spec (\Sym \rR\Gamma(X, \cE)^{*})",
\]
whose classical truncation is $\ul{H^0(X, \cE)} $. This represents the functor which, informally speaking, sends a derived scheme $S$ to $\rR\Gamma(X_S, \cE)$. What we have seen here is that derived algebraic geometry allows us to construct derived schemes which extend classical moduli of global sections in the same way that derived functors extend classical global sections. Now, we have put quotations in the formula above because we have not explained what it means to take the symmetric power of an animated $k$-module. It can be characterized as the left adjoint to the forgetful functor from animated commutative $k$-algebras to animated $k$-modules.

\begin{remark}[Derived vector bundles]\label{rem: derived vector complex}More generally, one can perform an analogous construction ``relative to a base'': given a perfect complex $\cK$ (of quasicoherent sheaves) on a (derived) stack $S$, one can form a ``total space'' $\Tot_S(\cK)$ as a derived stack over $S$. Such constructions are called \emph{derived vector bundles} in \cite{FYZ3}, since they generalize vector bundles, which are the total spaces arising in the special case where $\cK$ is a locally free coherent \emph{sheaf} (viewed as a perfect complex concentrated in degree $0$). 
\end{remark}

\begin{example}
Let $S = \Bun_n$, the moduli stack of rank $n$ vector bundles on smooth projective curve $C$. Consider the space $\rM \rightarrow S$, with $R$-points the groupoid of rank $n$ vector bundles $\cF$ on $C$ plus a global section of $\cF$. Although the fibers of $\rM \rightarrow S$ are vector spaces, the morphism $\rM \rightarrow S$ is not even flat, since the fibers have different dimensions (even when restricted to connected components of the base). This is because $\cF \mapsto H^0(\cF)$ behaves ``discontinuously'', for example because the dimension of $H^0(\cF)$ varies discontinuously with $\cF$. 

However, the perfect complex $\rR\Gamma(\cF)$ behaves ``continuously'', e.g., the Euler characteristic of $\rR\Gamma(\cF)$ varies continuously with $\cF$. This foreshadows the fact that we can assemble $R \Gamma(\cF)$ into a perfect complex on $S$, whose total space $\cM \rightarrow S$ has classical truncation $\rM \rightarrow S$ and whose derived fiber over $\cF$ is the derived scheme $R \Gamma(\cF)$ in the above sense. This $\cM$ is a variant of the ``derived Hitchin stacks'' to be discussed in \S \ref{sec: derived moduli} and \S \ref{sec: derived special cycles}. 
\end{example}

\subsubsection{Derived moduli spaces} Perhaps the most interesting examples of derived schemes (or stacks) are \emph{derived moduli spaces}. The entirety of \S \ref{sec: derived moduli} is devoted to such examples, so let us just give a brief teaser here. The idea here is that one begins with a classical moduli space, which by definition is a functor defined on commutative rings, satisfying appropriate descent conditions. Then one finds a way to possibly reformulate the moduli problem so that it makes sense on animated commutative rings as well, and in general outputs anima (i.e., animated sets) instead of sets. The representing object is called a derived moduli space. We shall see later some of the benefits of promoting moduli spaces to derived moduli spaces. 

\subsection{Quasismoothness}
In this subsection, we define a property of morphisms of derived schemes which is called \emph{quasismoothness}, which turns out to play a very important role in derived algebraic geometry. It is the derived generalization of a \emph{local complete intersection} morphism. For motivation, we recall what this means: 

\begin{defn}
Let $A$ be a noetherian ring and $f \co A \rightarrow B$ be a finite type morphism. We say that $f$ is a \emph{local complete intersection} (LCI) if it is Zariski-locally on $A$ the composition of a map of the form $A \rightarrow A[x_1, \ldots, x_n]$ followed by a quotient by a regular sequence. 
\end{defn}

Local complete intersections admit a clean characterization in terms of the cotangent complex. 

\begin{thm}[Avramov, \cite{Avr99}]\label{thm: LCI}
Let $f \co A \rightarrow B$ be a homomorphism of Noetherian rings. Then $f$ is LCI if and only if $\bL_f$ is represented by a complex with cohomology concentrated in degrees $[-1,0]$. 
\end{thm}

\begin{defn}\label{def: quasismooth}
A morphism $f \co \cA \rightarrow \cB$ of animated commutative rings is \emph{quasismooth} if $\bL_f$ has Tor-amplitude in $[-1,0]$, i.e., is Zariski-locally on $\cB$ represented by a complex of vector bundles in degrees $-1, 0$. (Here we are using \emph{cohomological} indexing.)  

More generally, a morphism $f \co \cX \rightarrow \cY$ of derived stacks is \emph{quasismooth} if $\bL_f$ exists and has tor-amplitude in $[-1, \infty)$, i.e., is locally represented by a perfect complex of vector bundles in degrees $[-1, n]$ for some $n<\infty$. 
\end{defn}

If $f$ is quasismooth, then we define the \emph{relative dimension of $f$} to be $\chi(\bL_f)$. 

Theorem \ref{thm: LCI} explains why we say that quasismoothness is a derived generalization of ``LCI''. This property enables certain important constructions, for example the construction of ``Gysin pullbacks'' in intersection theory. Moreover, there is a sense in which it is more common to encounter quasismooth morphisms than LCI morphisms in nature; this is actually one of the main motivations for the entrance of derived algebraic geometry into enumerative algebraic geometry. For example, once one knows the definitions it is easy to show: 

\begin{prop}\label{prop: quasismooth preservation}
Quasismooth morphisms are preserved by compositions and derived base changes. 
\end{prop}

By contrast, note that LCI morphisms are \emph{not} preserved by base change. 

\begin{example}\label{ex: quasismooth schemes}
Any finite type affine scheme can be realized as the classical truncation of a quasismooth derived affine scheme, for trivial reasons. Indeed, given any choice of presentation $R = \Z[x_1, \ldots, x_n]/(f_1, \ldots, f_m)$, the derived fibered product of the diagram 
\[
\begin{tikzcd}
& \Spec \Z[x_1, \ldots, x_n] \ar[d, "{(f_1, \ldots, f_m)}"] \\
\Spec \Z \ar[r, "0"] & \Spec \Z[y_1, \ldots, y_m]
\end{tikzcd}
\]
has classical truncation isomorphic to $\Spec R$. This shows that the property of being quasismooth imposes no restriction on the underlying classical truncation. 
\end{example}

If $f \co \cX \rightarrow \cY$ has a tangent complex $\bT_f$, then it can be used to reformulate quasismoothness: $f$ is quasismooth if and only if $\bT_f$ is represented by a perfect complex with tor-amplitude in $(-\infty, 1]$. 

\begin{example}\label{ex: quasismooth as derived fibered product} 
By Example \ref{ex: tangent complex of derived fibered product}, a derived scheme which is Zariski-locally isomorphic to a derived fibered product of smooth schemes is quasismooth. The converse is also true.
\end{example}

\section{Derived moduli spaces related to the Langlands program}\label{sec: derived moduli}

The relevance of derived algebraic geometry in the Langlands program is through \emph{derived moduli spaces}. In this section we will give an overview of some examples, discuss how to think about them, and hint at why they are useful.

\subsection{The hidden smoothness philosophy}
Historically, the notion of derived moduli spaces was anticipated in the 1980s, even before the advent of derived algebraic geometry, by (according to our understanding) Beilinson, Deligne, Drinfeld, and Kontsevich. They emphasized a principle called the \emph{hidden smoothness philosophy}, which predicts that natural moduli spaces should have derived versions which are \emph{quasismooth} (and of the ``expected dimension''). 

The backdrop for this principle is that in certain areas of classical algebraic geometry, one frequently encounters classical moduli spaces -- that is to say, classical schemes or stacks -- which are not LCI or have the ``wrong'' dimension. In such situations, the hidden smoothness philosophy predicts that these are the classical truncations of derived moduli spaces which do have the ``correct'' geometric properties. It can then be useful to find the definition of these derived enhancements. In this section we will discuss derived moduli spaces from this perspective: what to look for when trying to ``upgrade'' a classical moduli space to a derived one. Subsequent sections will (hopefully) illustrate why this might be a useful thing to do.

As far as we know, the motivations for the hidden smoothness philosophy are empirical, and we will try to illustrate it through examples. 

\begin{example}[Moduli of Betti local systems]\label{ex: betti local systems 0} Let us begin by brainstorming informally about a somewhat older example in derived algebraic geometry, which is relevant for (Betti) Geometric Langlands \cite{BZN15}: the moduli space of (Betti) local systems on a smooth projective (connected) curve $C/\CC$ of genus $g$. After choosing a basepoint, one could view a Betti local system as a representation of the fundamental group of $C$. There is a more canonical way to phrase things without choosing a basepoint, but for concreteness we make such a choice $c_0 \in C$, and then $\pi_1(C, c_0)$ has a standard presentation
\begin{equation}\label{eq: pi_1 presentation}
\pi_1(C,c_0) \approx \langle a_1, \ldots, a_g, b_1, \ldots, b_g \co \prod_{i=1}^g [a_i, b_i] = 1  \rangle.
\end{equation}
This presentation suggests a reasonable guess for how to write down the moduli space of rank $n$ (Betti) local systems on $C$, as a fibered product 
\begin{equation}\label{eq: betti locsys diagram}
\begin{tikzcd}
& \GL_n^{2g} \ar[d, "{\prod_{i=1}^g [x_i, y_i]}"]  \\
e \ar[r] & \GL_n
\end{tikzcd}
\end{equation}
quotiented by the conjugation action of $\GL_n$. Since all the spaces in the diagram \eqref{eq: betti locsys diagram} are smooth (and in particular the bottom horizontal arrow is a regular embedding), it suggests that the moduli space of local systems might be a local complete intersection, of dimension equal to $2g \dim \GL_n - 2\dim \GL_n$. This turns out to not always be the case, however, because a regular sequence cutting out $\{e\} \inj \GL_n$ does not always pull back to a regular sequence in $\GL_n^{2g}$. 

However, if instead of taking the fibered product of diagram \eqref{eq: betti locsys diagram} in classical schemes, one takes the \emph{derived} fibered product, then by Proposition \ref{prop: quasismooth preservation} the resulting derived scheme is quasismooth and has the ``expected'' dimension $2(g-1) \dim \GL_n$.
\end{example}

The Hidden Smoothness philosophy is justified in many examples by the fact that one can find natural (local) presentations for naturally occurring moduli spaces, similar to \eqref{eq: betti locsys diagram}. More precisely, the quasismoothness and dimensionality conditions are suggested by deformation theory, as shall be explained below. 

This section explains heuristics for how to solve the following problem: 
\begin{quote}
Given a classical moduli space $\rM$, find the ``correct'' derived moduli space $\cM$ of which it is a classical truncation. 
\end{quote}
By definition, a classical moduli space is a type of functor defined on the category of commutative rings, while a derived moduli space is a type of functor on the category of simplicial commutative rings, so the content of this problem is one of extending the domain of definition for the functor in a ``good'' way. (The target should also be extended, from sets or groupoids to simplicial sets.) In analogy to the notion of derived functors, we will informally call this problem ``deriving the classical moduli space $\rM$''. 

Of course, this problem is extremely ill-defined: a given classical moduli space can be realized as the classical truncation of many different (quasismooth) derived moduli spaces (see Example \ref{ex: quasismooth schemes}), just as a given reduced scheme can be realized as the underlying reduced scheme of many different possible schemes. Therefore, what we mean by ``correct'' derived moduli space is heuristic, and is dictated by the needs of the particular problem at hand. In practice, the problem may be rigidified by the desire to obtain a quasismooth $\cM$, or even to obtain an $\cM$ with a particular tangent complex.

In this section we will describe certain patterns and heuristics for how this process tends to be carried out in practice. We will illustrate it through two main examples: the \emph{derived Galois deformation rings} of Galatius-Venkatesh \cite{GV}, and the \emph{derived Hitchin stacks} of Feng-Yun-Zhang \cite{FYZ2}.

\subsection{Heuristics for derived moduli spaces}\label{ssec: heuristics}

Suppose that we begin with a certain classical moduli space $\rM$, which we want to promote to a \emph{derived} moduli space $\cM$. The following outline describes the most standard pattern for finding the ``correct'' derived enhancement $\cM$. 

Suppose that we are in the following situation: 
\begin{itemize}
\item The tangent space to $\rM$ at a point $m \in \rM$ can be calculated, and has some interpretation as a natural cohomology group. Furthermore, the automorphisms lie in a ``previous'' cohomology group, and obstructions lie in the ``next'' cohomology group. 
\end{itemize}
Then we should try to construct $\cM$ so that its tangent complex (cf. \S \ref{sssec: tangent complex}) is the cohomology complex computing the cohomology groups from the bullet point. 

We will give some explicit examples to make this concrete. Before that, however, we will discuss heuristics for when this process is unnecessary. 

\subsubsection{Heuristics for classicality}\label{ssec: classicality}

In this subsection we give some guiding principles for when one expects to find a derived enhancement which is not represented by a classical scheme or stack. The basic philosophy is captured by the following slogan. 

\begin{quote}
\textbf{Slogan.} If a moduli space is LCI of the ``correct'' dimension, then it does not need to be derived. 
\end{quote}

The notion of ``correct'' dimension is itself heuristic, and might arise from an intersection-theoretic setup such as in \eqref{eq: betti locsys diagram}, or from calculating the Euler characteristic of the expected cotangent complex. 

A little more precisely, the statement that a moduli space ``does not need to be derived'' means that the ``correct'' derived moduli space should just be isomorphic to its classical truncation.

The basis for this principle is as follows. Suppose $\cM$ is a derived scheme with classical truncation $\rM \inj \cM$. Following the Hidden Smoothness philosophy, we suppose that $\cM$ is quasismooth of the ``correct'' dimension, which at $m \in \cM$ is the Euler characteristic of the cotangent complex at $m$, $\chi(\bL_{\cM,m})$. 

\begin{lemma}\label{lem: classical truncation map of L}
The canonical map $\iota \co \rM \rightarrow \cM$ induces an isomorphism $H^0(\iota^*  \bL_{\cM})\xrightarrow{\sim} H^0(\bL_{\rM}) $ and a surjection $H^{-1}(\iota^* \bL_{\cM}) \surj H^{-1}(\bL_{\rM})$. 
\end{lemma}

We will not give a proof of Lemma \ref{lem: classical truncation map of L}, which one can find in \cite[\S 7.2]{GV}, but we will at least give some intuition for it. At the level of simplicial commutative rings, classical truncation is implemented by adding simplices in degrees $\geq 2$ to ``kill off'' higher homotopy groups. In topology, the Hurewicz theorem says roughly that a map of topological spaces induces an isomorphism of low-degree homotopy groups if and only if it induces an isomorphism of low-degree homology groups. Lemma \ref{lem: classical truncation map of L} follows from an analogous estimate for simplicial commutative rings. 

\begin{remark}
There is a variant of Lemma \ref{lem: classical truncation map of L} when $\cM$ is a derived stack, which says that the map $\iota^*  \bL_{\cM} \rightarrow \bL_{\rM}$ induces an isomorphism in degrees $\geq 0$ and a surjection in degree $-1$. 
\end{remark}

Now, suppose that $\rM$ is LCI of the ``correct'' dimension, which we take to mean $\chi(\bL_{\cM})$. Then the cotangent complex of $\rM$ has tor-amplitude in $[-1,0]$ by Theorem \ref{thm: LCI}, and its Euler characteristic is the same as that of $\bL_{\cM}$. Then Lemma \ref{lem: classical truncation map of L} forces $\iota^* \bL_{\cM} \rightarrow \bL_{\rM}$ to be an isomorphism, and then Lemma \ref{lemm: Lurie Hurewicz} implies that the map $\rM \rightarrow \cM$ is itself an isomorphism. 

We have not turned the above argument into a general formal statement, since the ``correct'' dimension is determined heuristically. The dimension used in the proof is $\chi(\bL_{\cM})$, but this is circular as it presupposes the existence of $\cM$. In practice, ``correct'' is determined by success in applications. We caution that examples do come up ``in nature'' where $\rM$ is LCI (or even smooth) but \emph{not} of the ``correct'' dimension, such as the derived Galois deformation rings (with crystalline conditions) studied in \cite{GV}. 

\begin{remark}
In practice, it can be difficult to verify that a moduli space $\rM$ is LCI of the correct dimension, unless $\rM$ is furthermore smooth of the correct dimension. This is because the cotangent complex of $\rM$ is typically not easy to compute, unless the obstruction group vanishes (so that it is just a sheaf, in which case $\rM$ is smooth). By contrast, for a \emph{derived} moduli space $\cM$, the cotangent complex can usually be described concretely because it has an interpretation in terms of ``higher derived dual numbers''; see \S \ref{ssec: local conditions} for more discussion of this. 
\end{remark}

\subsubsection{Running examples of classical moduli spaces} In the outline at the beginning of \S \ref{ssec: heuristics}, we said: 
\begin{quote}
One begins with a certain classical moduli space $\rM$, which one wants to promote to a derived moduli space $\cM$.
\end{quote}
We list some examples below, which will serve as running illustrations that we will return to repeatedly. 

\begin{example}[Moduli of Betti local systems]\label{ex: betti 1} Let $C$ be a smooth projective curve of genus $g$ and $G$ a reductive group over $\CC$. Motivated by Example \ref{ex: betti local systems 0}, we define the \emph{moduli stack of Betti $G$-local systems on $C$} to be the stack $\rM^{\Betti}_{G}$ obtained by taking the fibered product
\begin{equation}\label{eq: betti locsys G^}
\begin{tikzcd}
\rM^{\Betti, \square}_{G} \ar[d] \ar[r] & G^{2g} \ar[d, "{\prod_{i=1}^g [x_i, y_i]}"]  \\
e \ar[r] & G
\end{tikzcd}
\end{equation}
and then defining $\rM^{\Betti}_{G} := \rM^{\Betti, \square}_{G}/G$ to be the quotient by the conjugation $G$-action.

This moduli space is the subject of the ``spectral side'' of the Betti Geometric Langlands conjecture, which was proposed by Ben-Zvi -- Nadler \cite{BZN15}. 
\end{example}

\begin{example}[Global unrestricted Galois deformation ring]\label{ex: gal def 1}
Let $F$ be a global field and $\cO_F \subset F$ its ring of integers. Let $S$ be a finite set of primes of $\cO_F$ and $\cO_F[1/S]$ be the localization of $\cO_F$ at $S$. Let $G$ be a split reductive group and $\ol \rho \co  \pi_1(\cO_F[1/S]) \rightarrow G(\F_p)$ be a Galois representation.\footnote{We suppress the choice of basepoint; for a better perspective see \S \ref{sssec: derived global deformation functor}.} We define the \emph{(global unrestricted) Galois deformation functor of $\ol{\rho}$} to be the stack $\rM^{\Gal}_{G}$ on the category of Artinian $\Z_p$-algebras $A$ equipped with an augmentation $A \xrightarrow{\epsilon} \F_p$, with $A$-points being the groupoid of lifts $\rho \co \pi_1(\cO_F[1/S]) \rightarrow G(A)$ which are sent to $\ol{\rho}$ by the augmentation $\epsilon$, as in the diagram below:
\begin{equation}\label{eq: deformation lift diagram}
\begin{tikzcd}
& G(A) \ar[d, "\epsilon"] \\ 
\pi_1(\cO_F[1/S])\ar[r, "\ol \rho"'] \ar[ur, "\rho"] & G(\F_p)
\end{tikzcd}
\end{equation}
Such deformation functors were introduced by Mazur in \cite{Maz89}, and are at the foundation of the automorphy lifting theorems pioneered by Taylor-Wiles; see the article of Caraiani--Shin \cite{CS23} for more about this. (Here ``global'' indicates that we are looking at representations of the fundamental group of a ring of $S$-integers in a \emph{global} field, while the ``unrestricted'' refers to that we have not imposed any local conditions, which we would eventually want to do for applications to the Langlands correspondence.) 

\end{example}

\begin{remark}
The need for $\rM^{\Gal}_{G}$ to be derived depends subtly on $F$ and on $G$. We assume that $\ol{\rho}$ is absolutely irreducible and \emph{odd} (meaning that the trace of all complex conjugations is minimal -- see \cite{Gr}). For example, if $F = \Q$, and $G = \GL_1, \GL_2$ then one does not expect any interesting derived enhancement, while if $G = \GL_3, \GL_4, \ldots$ then one does. If $F $ is a number field that is not totally real and $G$ is semisimple, then $\rM^{\Gal}_G$ should be derived, whereas if $F$ is a function field then $\rM^{\Gal}_{G}$ need not be derived for any $G$ by \cite{deJ01, Gai07}. We mention finally that if $G = \GL_1$ (and $p \neq 2$), then $\rM^{\Gal}_{G}$ should be derived unless $F = \Q$ or a quadratic imaginary field; this will be discussed more in Example \ref{ex: gl_1}.
\end{remark}

\begin{example}[Local unrestricted Galois deformations]
Below we will discuss deriving (unrestricted) Galois deformation functors of global fields. One could also ask whether (unrestricted) Galois deformation functors of local fields should also be derived. It turns out that one can prove that these are always LCI of the correct dimension, so that this is unnecessary. See \cite{DHKM} for the $\ell \neq p$ case, and \cite{BIP} for the $\ell=p$ case. By contrast, the question of whether local deformation functors \emph{with conditions} should be derived is unclear at present, and seems to be an important problem for the future of derived Galois deformation theory; this is discussed more in \S \ref{ssec: local conditions}. 
\end{example}

\begin{example}[Hitchin stacks]\label{ex: hitchin 1}
Let $C$ be a smooth projective curve over $\F_q$ of characteristic $p>2$. For integers $m,n \geq 0$, we define the \emph{Hitchin stack $\rM^{\Hitch}_{m,n}$} to be the stack over $\F_q$ with $R$-points being the groupoid of tuples $(\cE_1, \cE_2, \cF, t_1, t_2)$ where 
\begin{itemize}
\item $\cE_i$ is a rank $m$ vector bundle on $C_R = C \times_{\F_q} \Spec R$, 
\item $\cF$ is a rank $n$ vector bundle on $C_R$
\item $t_1 \in \Hom(\cE_1, \cF)$ and $t_2 \in \Hom(\cF, \cE_2^\vee)$ where $\cE_2^\vee \cong \cE_2^* \otimes \omega_C$ is the Serre dual of $\cE_2$. 
\end{itemize}
This is the specialization of the unitary Hitchin spaces $\rM$ from \cite{FYZ, FYZ2} to the ``split case'', where the unitary group is split. It was introduced in order to analyze special cycles on moduli stacks of shtukas, which feature into a function field version of arithmetic theta series.\footnote{The spaces $\rM^{\Hitch}_{m,n}$ are not the spaces of Higgs bundles which were originally studied by Hitchin, but are other instances of a similar construction that has risen to importance in automorphic representation theory -- see \cite{Y18, FW} for more discussion.} 
\end{example}

\subsubsection{Hints from deformation theory} We now suppose the classical moduli space $\rM$ to be given. Next, in \S \ref{ssec: heuristics}, we said: 

\begin{quote}
The tangent space to $\rM$ at a point $m \in \rM$ can be calculated, and has some interpretation as a natural cohomology group. Furthermore, the automorphisms lie in a ``previous'' cohomology group, and obstructions lie in the ``next'' cohomology group.
\end{quote}

We will now explain this. The point is that the tangent space to $\rM$ at $m \in \rM$ has a ``concrete'' description in terms of first-order deformations, by virtue of the definition of $\rM$ as a moduli space, and the interpretation of the tangent space in terms of dual numbers. 

\begin{example}[Moduli of Betti local systems]\label{ex: betti 2}
Continuing the notation of Example \ref{ex: betti 1}, let $L$ be a (Betti) $G$-local system on $C$. We may regard $L$ as a $\CC$-point of $\rM^{\Betti}_{G}$. Then the tangent space to $\rM^{\Betti}_{G}$ at $L$ is $H^1(\pi_1(C), \mf{g}_L )$ where $\mf{g}_L := L\times^{G}  \mf{g} $ is the local system on $C$ obtained by taking the quotient of $L \times \mf{g}$ by the diagonal $G$-action. This local system is also sometimes denoted $\Ad L$ or $\mf{g}_L$. Furthermore, 
\begin{itemize}
\item The automorphisms of $L$ may be canonically identified with $H^0(\pi_1(C), \mf{g}_L)$, and \item One can show that the obstruction to first-order deformations lie in $H^2(\pi_1(C), \mf{g}_L)$. 
\end{itemize}
\end{example}

\begin{example}[Global unrestricted Galois deformations]\label{ex: gal def 2}
Continuing the notation of Example \ref{ex: gal def 1}, view $\ol{\rho}$ as an $\F_p$-point of $\rM^{\Gal}_{G}$. Then the tangent space to $\rM^{\Gal}_{G}$ at $\ol{\rho}$ is $H^1(\pi_1(\cO_F[1/S]) , \mf{g}_{\ol{\rho}})$, where $\mf{g}_{\ol{\rho}}$ has underlying vector space $\mf{g}_{\F_p}$ with the action of $\pi_1(\cO_F[1/S])$ being through $\ol{\rho} \co \pi_1(\cO_F[1/S]) \rightarrow G(\F_p)$ and the adjoint action of $G(\F_p)$ on $\mf{g}_{\F_p}$. Furthermore, 
\begin{itemize}
\item The automorphisms of $\ol{\rho}$ may be canonically identified with $H^0(\pi_1(\cO_F[1/S]) , \mf{g}_{\ol{\rho}})$, and 
\item One can show that the obstructions to first-order deformations lie in $H^2(\pi_1(\cO_F[1/S]), \mf{g}_{\ol{\rho}})$ \cite[Proposition 2]{Maz89}. 
\end{itemize}
Note the similarity to Example \ref{ex: betti 2}. 
\end{example}

\begin{example}[Hitchin stack]\label{ex: hitchin 2}
Continuing the notation of Example \ref{ex: hitchin 1}, let $(\cE_1,\cE_2, \cF, t_1,t_2)  \in \rM^{\Hitch}_{m,n}(k)$ for a field $k/\F_q$. Let $\msf{T}(\cE_1,\cE_2, \cF, t_1,t_2) $ be the complex
\[
\End(\cE_1) \oplus \End(\cE_2) \oplus \End(\cF) \xrightarrow{d} \Hom(\cE_1, \cF) \oplus \Hom(\cF, \cE_2^\vee)
\]
where the left term lies in degree $-1$, the right term lies in degree $0$, and the differential sends $(A_1, A_2, B) \in \End(\cE_1) \oplus \End(\cE_2) \oplus \End(\cF)$ to $(Bt_1 - t_1A_1, A_2^\vee t_2- t_2B) \in \Hom(\cE_1, \cF) \oplus \Hom(\cF, \cE_2^\vee)$. Then the tangent space to $\rM^{\Hitch}_{m,n}$ at $(\cE_1,\cE_2, \cF, t_1,t_2)$ may be canonically identified with $H^1(C_k, \msf{T}(\cE_1,\cE_2, \cF, t_1,t_2))$. Furthermore, as shown in \cite[\S 4.14]{Ngo10}, 
\begin{itemize}
\item The automorphisms of $(\cE_1,\cE_2, \cF, t_1,t_2)$ may be canonically identified with $H^0(C_k, \msf{T}(\cE_1,\cE_2, \cF, t_1,t_2))$, and 
\item The obstruction to first-order deformations lies in $H^2(C_k, \msf{T}(\cE_1,\cE_2, \cF, t_1,t_2))$. 
\end{itemize}
\end{example}

\subsubsection{Construction of derived moduli stacks} At this point we have our classical moduli space $\rM$ together with a candidate guess for the (co)tangent complex of the derived moduli stack $\cM$. Finally, in \S \ref{ssec: heuristics}, we said: 

\begin{quote}One then tries to construct $\cM$ so that its tangent complex is the cohomology complex computing the cohomology groups from the previous step.
\end{quote}

We will illustrate how this plays out in our running examples. 

\begin{example}[Derived moduli of Betti local systems]\label{ex: betti 3}
We maintain the notation of Example \ref{ex: betti 2}. In fact, there is already a subtle wrinkle here: the cohomology groups in Example \ref{ex: betti 2} were group cohomology of $\pi_1(C)$ with coefficients in a local system $\mf{g}_L$. The only dependence of this on $C$ is through its fundamental group $\pi_1(C)$, which can be viewed as capturing the information of the Postnikov truncation $\tau_{\leq 1} C$, but the right thing to do is to incorporate the whole ``homotopy type'' of $C$. This means that we want the tangent complex at the local system $L$ to be the cohomology complex $\rR\Gamma_{\Betti}(C, \mf{g}_L[1])$ rather than $\rR\Gamma_{\mrm{group}}(\pi_1(C), \mf{g}_L[1])$. 

Let us briefly discuss the difference between these two objects. We have 
\[
H^0(C,  \mf{g}_L) \cong H^0(\pi_1(C), \mf{g}_L)
\]
and
\[
H^1(C, \mf{g}_L) \cong H^1(\pi_1(C), \mf{g}_L)
\]
so the difference lies in $H^2$. Furthermore, if $C$ has genus $g \geq 1$ then in fact $C$ is a $K(\pi_1(C), 1)$, so that $\rR\Gamma_{\Betti}(C, \mf{g}_L[1]) \cong \rR\Gamma_{\mrm{group}}(\pi_1(C), \mf{g}_L[1])$. Therefore this distinction only arises (among smooth proper complex curves) in the case where $C  \cong \PP^1$. In that case $\pi_1(C)$ is trivial, so the group cohomology is trivial. However, the higher Betti cohomology of $\PP^1$ is non-trivial (note that $L$ is necessarily ``the'' trivial local system, since $\pi_1(\PP^1) = 0$). 

One way to construct the desired derived moduli stack $\cM^{\Betti}_{G}$ is to take the \emph{derived} fibered product of \eqref{eq: betti locsys diagram}, and then quotient out by the (diagonal) conjugation action of $G$. To calculate the tangent complex of the resulting object, use that the formation of tangent complex respects homotopy limits (Example \ref{ex: tangent complex of derived fibered product}), and the presentation of the homotopy type of $C(\CC)$ as the homotopy pushout 
\begin{equation}\label{eq: C homotopy pushout}
\begin{tikzcd}
S^1 \ar[r] \ar[d] & \bigvee_{i=1}^{2g} (S^1) \ar[d] \\
\pt \ar[r] & C(\CC)
\end{tikzcd}
\end{equation}

Note that for $g = 0$, where $C \cong \PP^1$, our presentation of $\cM^{\Betti}_{G}$ is the \emph{derived} self-intersection of the identity in $G$, quotiented by $G$. This has dimension $-2\dim G$. The classical truncation of $\cM^{\Betti}_{G}$ is simply $\rM^{\Betti}_{G} \cong [\Spec \CC]/G$, reflecting that there is a unique local system on $\PP^1$ (the trivial one). A more na\"{i}ve attempt to upgrade $\rM^{\Betti}_{G}$ to a derived moduli space, with tangent complex $\rR\Gamma(\pi_1(C), \mf{g}_L[1])$, would have simply produced $[\Spec \CC]/G$, which is still quasismooth but has the wrong dimension $- \dim G$.

Note that the description of the tangent complex of $\cM^{\Betti}_{G}$ as the perfect complex isomorphic to $\rR\Gamma_{\Betti}(C_R,  \mf{g}_L[1])$, functorially in $L \in \cM^{\Betti}_{G}(R)$, shows that $\cM^{\Betti}_{G}$ is quasismooth, since $C$ has Betti cohomological amplitude in $[0,2]$. 
\end{example}

\begin{example}[Derived global unrestricted Galois deformations]\label{ex: gal def 3} 
We maintain the notation of Example \ref{ex: gal def 2}. The calculations of Example \ref{ex: gal def 2} suggest that we should look for a derived stack $\cM^{\Gal}_{G}$ whose tangent complex at $\rho \in \cM^{\Gal}_{G}(R)$ should be isomorphic naturally in $R$ to $\rR\Gamma(\pi_1(\cO_F[1/S]), \mf{g}_{\rho}[1])$. Note that one might expect, based on Example \ref{ex: betti 3}, that this answer should be adjusted to $\rR\Gamma_{\et}( \cO_F[1/S], \mf{g}_{\ol{\rho}}[1])$. It may indeed be better to think of the latter as the ``right'' answer, but the ring of $S$-integers in a global field is an ``\'{e}tale $K(\pi, 1)$'', so that the two always coincide.\footnote{
However, when $F$ is the function field of a smooth projective curve $C/\F_q$, then it is also natural to study \emph{everywhere unramified} deformations of a local system on $C$. In this case the tangent complex is $\rR\Gamma_{\et}(C; \mf{g}_{\ol{\rho}}[1])$, and it can have nontrivial $H^3_{\et}$. Therefore, the Galois deformation space is \emph{not} quasismooth in this case. }

The desired $\cM^{\Gal}_{G}$ is the \emph{derived Galois deformation functor} of Galatius-Venkatesh. We will not give its construction right now; that requires a bit of digression into homotopy theory. Morally, one wants to define the functor of points to be \eqref{eq: deformation lift diagram} but with $A$ an ``Artinian animated commutative ring augmented over $\F_p$''. The resulting object should be a functor from the category of Artinian animated commutative rings augmented over $\F_p$ to the category of anima (also known as ``simplicial sets'', ``spaces'', or ``$\infty$-groupoids'').

The determination of the cohomology of number fields is part of class field theory and Poitou-Tate duality. The output is that if $p \neq 2$, then $\rR\Gamma(\pi_1(\cO_F[1/S]), \mf{g}_{\ol{\rho}})$ has tor-amplitude in $[0, 2]$, so that $\cM^{\Gal}_{G}$ is quasismooth.

\end{example}

\begin{example}[Derived Hitchin stack]\label{ex: hitchin 3}
Continuing with the notation of Example \ref{ex: hitchin 2}, we expect that the \emph{derived Hitchin stack} $\cM^{\Hitch}_{m,n}$ should be such that the pullback of its tangent complex to $(\cE_1, \cE_2, \cF, t_1, t_2) \in \cM^{\Hitch}_{m,n}(R)$ is naturally isomorphic to $\rR\Gamma(C_R, \msf{T}(\cE_1, \cE_2, \cF, t_1, t_2))$. This $\cM^{\Hitch}_{m,n}$ is the ``split'' special case of the derived Hitchin stacks that appear in \cite{FYZ2}. It can be derived as a ``derived mapping stack'' in the sense of \cite{TV08}. Alternatively, it may be viewed as the derived vector bundle (\S \ref{sssec: derived vector bundles}) associated to the perfect complex taking $(\cE_1, \cE_2, \cF) \in \Bun_{\GL_m} \times \Bun_{\GL_m}  \times \Bun_{\GL_n}(R)$ to $\RHom_{C_R}(\cE_1, \cF) \oplus \RHom_{C_R}(\cF, \cE_2^\vee)$.

\end{example}

\begin{remark}As we see from these examples, the justification for the hidden smoothness philosophy in spaces arising in number theory seems to be the low cohomological dimensions of number fields and related rings, which lead to ``short'' tangent complexes. 
\end{remark}

\begin{example}[Spectral Hitchin stacks and spectral periods] 
Non-quasismooth derived moduli spaces are also of interest in the Langlands program. Several examples are studied in \cite{FW, BZSV}, where it is argued that they geometrize a spectral (i.e., Galois-theoretic) analogue of \emph{periods} in automorphic representation theory. An example is the \emph{spectral Hitchin stack} $\rM_{m,n}^{\spec}$, which parametrizes $(E_1, E_2, F, t_1, t_2)$ analogous to the points of $\rM_{m,n}$ except with bundles replaced by local systems. For concreteness, let $C/\CC$ be a smooth projective curve (but it is possible to formulate $\ell$-adic variants for curves over finite fields). The $R$-points of $\rM_{m,n}^{\spec}$ are triples $(E_1, E_2, F,t_1, t_2)$ where
\begin{itemize}
\item $E_i$ is an $R$-family of rank $m$ local systems on $C$,
\item $F$ is an $R$-family of rank $n$ local systems on $C$, and 
\item $t_1 \in \Hom(E_1,F)$ is a flat morphism from $E_1$ to $F$, and $t_2 \in \Hom(F, E_2^*)$ is a flat morphism from $F$ to $E_2^*$. 
\end{itemize}
This has a derived enhancement $\cM_{m,n}^{\spec}$ whose tangent complex is the \emph{de Rham} cohomology of $C$ with coefficients in an analogous complex $\msf{T}(E_1,E_2, F, t_1, t_2)$ to that of Example \ref{ex: hitchin 2}. However, since the de Rham cohomological dimension of $C$ is $2$ while the coherent cohomological dimension of $C$ is $1$, this tangent complex actually has non-zero cohomology in degree $2$, so that $\cM_{m,n}^{\spec}$ is \emph{not} quasismooth. In fact, the failure of $\cM_{m,n}^{\spec}$ to be quasismooth is an important facet of \cite{FW}, which studies the relation between automorphic periods and spectral periods. One of the main points of \cite{FW} is that the \emph{dualizing sheaves} of spectral Hitchin stacks, which in some sense measure the failure of quasismoothness, are the spectral counterpart to automorphic periods. 
\end{example}

\section{Derived special cycles and higher theta functions}\label{sec: derived special cycles}
We will explain applications of the derived Hitchin stacks introduced in Example \ref{ex: hitchin 3} towards the theory of special cycles and higher arithmetic theta series, developed in papers of Feng-Yun-Zhang \cite{FYZ, FYZ2, FYZ3}. 

The context for this discussion is the Kudla program, as discussed in the article of Chao Li \cite{Li24}. This concerns the so-called \emph{arithmetic theta functions}, which are Fourier series assembled out of special cycles. These were first studied on orthogonal or unitary Shimura varieties, but we will focus on the function-field context, as in the work of Feng-Yun-Zhang. We refer to \cite{Li24} for the classical background and motivations. 

\subsection{Hermitian shtukas}We recall some of the definitions from \cite[\S 6,7]{FYZ}.

Let $X$ be a smooth projective curve over $\F_q$. Let $\nu \co X' \rightarrow X$ be a finite \'{e}tale double cover, with the non-trivial automorphism of $X'$ over $X$ denoted $\sigma \co X' \rightarrow X'$. 

For each $r \geq 0$ and $n > 0$, \cite[\S 6]{FYZ} and \cite[\S 2]{FYZ2} construct moduli spaces of rank $n$ unitary shtukas, denoted $\Sht_{U(n)}^r$. These are variants of the spaces discussed in the lectures of Zhiwei Yun. There is a map $\pi \co \Sht_{U(n)}^r \rightarrow (X')^r$, called the ``leg map''. For $r=1$, $\pi$ is to be thought of as roughly analogous to the structure map for the integral model of a unitary Shimura variety.\footnote{As explained by Yun, a closer analogy for the latter is the fiber of the $\pi$ for $r=2$ over $X' \times \{\infty\}$ for a fixed point $\infty \in X'$.} We will now outline these constructions.

For a test scheme $S$, a rank $n$ \emph{Hermitian bundle} on $X \times S$, with respect to $\nu \co X'\to X$, is a vector bundle $\Cal{F}$ of rank $n$ on $X' \times S$, equipped with an isomorphism $h \co \Cal{F} \xrightarrow{\sim} \sigma^* \Cal{F}^{\vee}$ such that $\sigma^* h^{\vee} = h$. Here, $\cF^\vee = \cF^* \otimes_{\cO_{X'}} \omega_{X'}$ is the Serre dual of $\cF$. We refer to $h$ as the ``Hermitian structure''. 

We denote by $\Bun_{U(n)}$  the moduli stack of rank $n$ Hermitian bundles on $X$, which sends a test scheme $S$ to the groupoid of rank $n$ vector bundles on $X' \times S$ equipped with a Hermitian structure.

\begin{defn}\label{defn: Hk U(n)}
Let $r \geq 0$ be an integer. The \emph{Hecke stack} $\Hk_{U(n)}^{r}$ has as $S$-points the groupoid of the following data: 
\begin{enumerate}
\item $x_i' \in X'(S)$ for $i = 1, \ldots, r$, with graphs denoted by $\Gamma_{x'_i} \subset X' \times S$.
\item A sequence of vector bundles $\Cal{F}_0, \ldots, \Cal{F}_r$ of rank $n$ on $X' \times S$, each equipped with Hermitian structure $h_i \co \Cal{F}_i \xrightarrow{\sim} \sigma^* \Cal{F}_i^{\vee}$. 
\item Isomorphisms $f_i \co \Cal{F}_{i-1}|_{X' \times S - \Gamma_{x'_i} - \Gamma_{\sigma(x'_i)}} \xrightarrow{\sim} \Cal{F}_i|_{X' \times S - \Gamma_{x'_i}-\Gamma_{\sigma(x'_i)}}$, for $1 \leq i \leq r$, compatible with the Hermitian structures, with the following property: there exists a rank $n$ vector bundle $\Cal{F}^{\flat}_{i-1/2}$ on $X' \times S$ and a diagram of vector bundles
\begin{equation}\label{eq: modification diagram 1}
\begin{tikzcd}
&\Cal{F}_{i-1/2}^{\flat} \ar[dl, "f_i^{\leftarrow}"', hook']  \ar[dr, "f_i^{\rightarrow}", hook ]  &  \\
\Cal{F}_{i-1} & & \Cal{F}_i
\end{tikzcd}
\end{equation}
such that $\coker(f_i^{\leftarrow})$ is locally free of rank $1$ over $\Gamma_{x'_i}$, and $\coker(f_i^{\rightarrow})$ is locally free of rank $1$ over $\Gamma_{\sigma(x'_i)}$. In particular, $f_i^{\leftarrow}$ and $f_i^{\rightarrow}$ are invertible upon restriction to $X' \times S - \Gamma_{x'_i}-\Gamma_{\sigma(x'_i)}$, and the composition 
\[
\Cal{F}_{i-1}|_{X' \times S - \Gamma_{x'_i}-\Gamma_{\sigma(x'_i)}} \xrightarrow{(f_i^{\leftarrow})^{-1}} \Cal{F}_{i-1/2}^{\flat}|_{X' \times S - \Gamma_{x'_i}-\Gamma_{\sigma(x'_i)}} \xrightarrow{f_i^{\rightarrow}} \Cal{F}_{i} |_{X' \times S - \Gamma_{x'_i}-\Gamma_{\sigma(x'_i)}} 
\]
agrees with $f_i$. 
\end{enumerate}

\end{defn}

For a vector bundle $\Cal{F}$ on $X' \times S$, we denote by $\ft \Cal{F} := (\Id_{X'} \times \Frob_S)^* \Cal{F}$ its Frobenius twist. If $\Cal{F}$ has a Hermitian structure $h \co \Cal{F} \xrightarrow{\sim} \sigma^* \Cal{F}^{\vee}$, then $\ft \Cal{F}$ is equipped with the Hermitian structure $\ft h$; we may suppress this notation when we speak of the ``Hermitian bundle'' $\ft \Cal{F}$. Viewing $\cF \in \Bun_{U(n)}(S)$, $\ft \cF$ is the image of $\cF$ under the map $\Frob \co \Bun_{U(n)} \rightarrow \Bun_{U(n)}$. 

\begin{defn}
Let $r \geq 0$ be an integer. We define $\Sht_{U(n)}^{r}$ by the Cartesian diagram
\[
\begin{tikzcd}
\Sht_{U(n)}^{r} \ar[r] \ar[d] & \Hk_{U(n)}^{r} \ar[d] \\
\Bun_{U(n)} \ar[r, "{(\Id, \Frob})"]  & \Bun_{U(n)} \times \Bun_{U(n)}
\end{tikzcd}
\]
A point of $\Sht_{U(n)}^{r}$ will be called a ``Hermitian shtuka (of rank $n$)''.

Concretely, the $S$-points of $\Sht_{U(n)}^{r}$ are given by the groupoid of the following data: 
\begin{enumerate}
\item $x'_i \in X'(S)$ for $i = 1, \ldots, r$, with graphs denoted $\Gamma_{x'_i} \subset X' \times S$. These are called the \emph{legs} of the shtuka. 
\item A sequence of vector bundles $\Cal{F}_0, \ldots, \Cal{F}_r$ of rank $n$ on $X' \times S$, each equipped with a Hermitian structure $h_i \co \Cal{F}_i \xrightarrow{\sim} \sigma^* \Cal{F}_i^{\vee}$. 
\item Isomorphisms $f_i \co \Cal{F}_{i-1}|_{X' \times S - \Gamma_{x'_i} - \Gamma_{\s(x'_i)}} \xrightarrow{\sim} \Cal{F}_i|_{X' \times S - \Gamma_{x'_i}-\Gamma_{\s(x'_i)}}$ compatible with the Hermitian structure, which as modifications of the underlying vector bundles on $X' \times S$ are as in the definition of the Hecke stack $\Hk_{U(n)}^r$.

\item An isomorphism $\varphi \co \Cal{F}_r \cong \ft \Cal{F}_0$ compatible with the Hermitian structure. 
\end{enumerate}
\end{defn}

\begin{remark}One can show that $\Sht_{U(n)}^r$ is a Deligne-Mumford stack locally of finite type. The map $\Sht_{U(n)}^r \rightarrow (X')^r$ is smooth, separated, equidimensional of relative dimension $r(n-1)$. For these and more geometric properties, see \cite[\S 6.4]{FYZ}. 
\end{remark}

\subsection{Special cycles}
Now we will define special cycles on $\Sht_{U(n)}^r$, which are analogous to the special cycles constructed by Kudla-Rapoport on unitary Shimura varieties \cite{KRII}.

\begin{defn}\label{def: Z} Let $\Cal{E}$ be a rank $m$ vector bundle on $X'$. We define the stack $\rZ_{\Cal{E}}^{r}$ whose $S$-points are given by the groupoid of the following data: 
 \begin{itemize}
 \item A rank $n$ Hermitian shtuka $(\{x'_1, \ldots, x'_r\}, \{\Cal{F}_0, \ldots, \Cal{F}_r\}, \{f_1, \ldots, f_r\}, \varphi) \in  \Sht_{U(n)}^{r}(S)$.
 \item Maps of coherent sheaves  $t_i \co \Cal{E} \boxtimes \cO_{S} \rightarrow \Cal{F}_i$ on $X'\times S$ such that the isomorphism $\varphi \co \Cal{F}_r \cong \ft \Cal{F}_0$ intertwines $t_r$ with $\ft t_0$, and the maps $t_{i-1}, t_{i}$ are intertwined by the modification $f_i \co  \Cal{F}_{i-1} \dashrightarrow \Cal{F}_{i}$ for each $i = 1, \ldots, r$, i.e., the diagram below commutes:
\[
\begin{tikzcd}
\Cal{E}\boxtimes\cO_{S} \ar[d, "t_0"] \ar[r, equals] &  \Cal{E}\boxtimes\cO_{S} \ar[r, equals]  \ar[d, "t_1"]& \ldots \ar[d] \ar[r, equals] & \Cal{E}\boxtimes\cO_{S} \ar[r, "\sim"] \ar[d, "t_{r}"] & \ft (\Cal{E} \boxtimes\cO_{S})\ar[d, "\ft t_0"] \\
\Cal{F}_0 \ar[r, dashed, "f_1"] & \Cal{F}_1 \ar[r, dashed, "f_2"]  &  \ldots \ar[r, dashed, "f_r"] &  \Cal{F}_{r} \ar[r, "\sim"] & \ft \Cal{F}_0  
\end{tikzcd}
\]
\end{itemize}
In the sequel, when writing such diagrams we will usually just omit the ``$\boxtimes \cO_S$'' factor from the notation.

\end{defn}

\begin{remark}There is an evident map $\rZ_{\Cal{E}}^{r} \rightarrow \Sht_{U(n)}^r$ projecting to the data in the first bullet point. This map is finite by \cite[Proposition 7.5]{FYZ}, hence induces a pushforward map on Chow groups. 
\end{remark}

\begin{defn}Let $\rA_{\cE}(\F_q)$ be the $\F_q$-vector space of Hermitian maps $a \co \Cal{E} \rightarrow \sigma^* \Cal{E}^{\vee}$ such that $\sigma^* a^{\vee} =a$. 
\end{defn}

Let $(\{x'_i\}, \{\Cal{F}_i\}, \{f_i\}, \varphi, \{t_i\}) \in \rZ_{\Cal{E}}^r(S)$. By the compatibilities between the $t_i$ in the definition of $\rZ_{\Cal{E}}^r$, the compositions of maps in the sequence  
\begin{equation}\label{eq: hitchin map formula}
\Cal{E} \boxtimes \cO_S \xrightarrow{t_i} \Cal{F}_i \xrightarrow{h_i} \sigma^* \Cal{F}_i^{\vee} \xrightarrow{\sigma^* t^{\vee}_{i}} \sigma^* \Cal{E}^{\vee} \boxtimes \cO_S
\end{equation}
agree for each $i$, and \eqref{eq: hitchin map formula} for $i=r$ also agrees with the Frobenius twist of \eqref{eq: hitchin map formula} for $i=0$. Hence the composite map \eqref{eq: hitchin map formula} gives the same point of $ \rA_{\Cal{E}}(S)$ for every $i$, which is moreover fixed by Frobenius, hence must come from $\rA_{\Cal{E}}(\F_q)$. This defines a map $\rZ_{\Cal{E}}^{r} \rightarrow \rA_{\Cal{E}}(\F_q)$. For $a \in \rA_{\Cal{E}}(\F_q)$, we denote by $\rZ^{r}_{\Cal{E}}(a)$ the fiber of $\rZ^{r}_{\Cal{E}}$ over $a$. We have
\begin{equation*}
\rZ_{\cE}^{r}=\coprod_{a\in \rA_{\cE}(\F_q)} \rZ^{r}_{\cE}(a).
\end{equation*}
The $\rZ_{\cE}^r$ or $\rZ_{\cE}^r(a)$ are called \emph{special cycles} of corank $m$ (with $r$ legs).

\subsection{Virtual fundamental classes of special cycles} 

Motivated by the Kudla program, we expect to attach to each special cycle $\rZ_{\cE}^r(a)$, $a \in \rA_{\cE}(\F_q)$, a class $[\rZ_{\cE}^r(a)]$ in the Chow group $\Ch_*(\Sht_{U(n)}^r)$ such that the Fourier series with Fourier coefficients $[\rZ_{\cE}^r(a)]$ is ``automorphic''. We call such an automorphic form a \emph{higher theta function}, since in the case $r=0$ it recovers the definition of classical theta functions. 

However, a cursory inspection reveals that for $r>0$, the cycles $\rZ_{\cE}^r(a)$ necessarily have varying dimensions as $a$ varies. For example, if $a = 0$ then $\rZ_{\cE}^r(0)$ clearly surjects onto $\Sht_{U(n)}^r$. On the other hand, it is not hard to see that if $r>0$ and $a \neq 0$ then $\rZ_{\cE}^r(a)$ has smaller dimension than $\Sht_{U(n)}^r$. Therefore, we should not use the na\"{i}ve cycle classes $[\rZ_{\cE}^r(a)]^{\naive}$ in the definition of higher theta functions; we must construct \emph{virtual} fundamental cycles $[\rZ_{\cE}^r(a)]^{\vir}$, which should all lie in $\Ch_{(n-m)r}(\Sht_{U(n)}^r)$ where $m = \rank \cE$. 

\subsubsection{The non-singular terms} There are various heuristics which guide such a construction. A key input is that when $m = \rank \cE = 1$, and when $a$ is \emph{non-singular}, meaning that $a \co \cE \rightarrow \sigma^* \cE^\vee$ is injective as a map of coherent sheaves (equivalently, an isomorphism over the generic point of $X'$), then $\rZ_{\cE}^r(a)$ is LCI of the ``correct'' codimension, which is $r$. In this case, heuristics (justified a posteriori by \S \ref{ssec: classicality}) suggest that one may take $[\rZ_{\cE}^r(a)]^{\vir} = [\rZ_{\cE}^r(a)]^{\naive}$. More generally, using this observation one can define $[\rZ_{\cE}^r(a)]^{\vir}$ for any non-singular $a$ (the point being to allow $m >1$) using the original method of \cite{KRII} of presenting $\rZ_{\cE}^r(a)$ as an open-closed component of the ``derived intersection'' of special cycles labeled by $m=1$ and non-singular coefficients. The definition of $[\rZ_{\cE}^r(a)]^{\vir}$ in this case is carried out in \cite[\S 7]{FYZ}. 

\subsubsection{The singular terms} The definition of $[\rZ_{\cE}^r(a)]^{\vir}$ for singular $a$ requires further heuristics. Over a number field, the analogous problem has a relatively simple answer. In that case, $a$ can be represented as a matrix which is of the form $a_{ns} \oplus 0$, where $a_{ns}$ is non-singular of rank $m' \leq m$. Then the corresponding virtual fundamental class is the product of the virtual fundamental class for $a_{ns}$ times the $(m-m')$th power of the Chern class of a certain ``tautological'' line bundle. 

But in the function field case at hand, the answer is much more complicated: it is in general a sum of infinitely many terms, each involving the top Chern class of a different ``tautological'' bundle. We can already illustrate why in the special case $a = 0$. Then $\rZ_{\cE}^r(0)$ is stratified by the kernel of the map $t_i \co \cE \boxtimes \cO_S \rightarrow \cF$, and each such stratum contributes a piece of the form described above (a Chern class times a virtual class which can be constructed via derived intersections as in the non-singular case). The analogous number field situation has the special property that the governing Hermitian form is \emph{positive definite} on the generic fiber, which rules out all but one stratum there. This is discussed in more detail in \cite[\S 4]{FYZ2} and \cite[Introduction]{FYZ3}. 

The upshot is that $[\rZ_{\cE}^r(a)]^{\vir}$ can be constructed explicitly for all terms $a \in \rA_{\cE}(\F_q)$, but the answer is quite complicated. 

Alternatively, it was observed in \cite{FYZ2} that the definition of the special cycles $\rZ_{\cE}^r$ can be naturally promoted, according to the pattern of \S \ref{sec: derived moduli}, to \emph{derived special cycles} $\cZ_{\cE}^r$ which are quasismooth of the correct dimension, and whose classical truncations recover $\rZ_{\cE}^r$. According to \cite{KhanI}, any quasismooth derived stack $\cZ$ has an intrinsic notion of fundamental class $[\cZ]^{\naive}$.\footnote{This rough statement had long been a folklore theorem of derived algebraic geometry, announced for example in \cite{Lur04}, and its importance for Gromov-Witten or Donaldson-Thomas flavors of enumerative geometry was long understood. However, Khan's construction of $[\cZ]^{\vir}$ is the most general to appear in the literature, and \cite{FYZ2} uses the full scope of its generality.} In \cite[Theorem 6.5]{FYZ2} it is proved that $[\cZ_{\cE}^r(a)]^{\naive}$ coincides with the previously defined (elementary) construction of $[\rZ_{\cE}^r(a)]^{\vir}$. 

To sketch how $\cZ_{\cE}^r$ is defined, the key point is to realize that $\sqcup_{\cE} \rZ_{\cE}^r$ admits a description as a fibered product 
\begin{equation}\label{eq: Z fibered product}
\begin{tikzcd}
\sqcup_{\cE} \rZ_{\cE}^r \ar[r] \ar[d] & \Hk_{\rM}^r  \ar[d] \\
\rM \ar[r, "{(\Id, \Frob)}"] & \rM \times \rM
\end{tikzcd}
\end{equation}
for $\rM$ a unitary generalization of the Hitchin stack from Example \ref{ex: hitchin 1}, and $\Hk_{\rM}^r$ a certain Hecke correspondence for $\rM$. Now, $\rM$ and $\Hk_{\rM}^r$ have natural enhancements to (quasismooth) derived stacks $\cM$ and $\Hk_{\cM}^r$; for example $\cM$ is a unitary generalization of Example \ref{ex: hitchin 3}. Actually, the important point for the quasismoothness of $\cZ_{\cE}^r$ is that the maps $\Hk_{\cM}^r \rightarrow \cM$ are quasismooth. Then, we define $\sqcup_{\cE} \cZ_{\cE}^r$ as a derived fibered product 
\[
\begin{tikzcd}
\sqcup_{\cE} \cZ_{\cE}^r \ar[r] \ar[d] & \Hk_{\cM}^r  \ar[d] \\
\cM \ar[r, "{(\Id, \Frob)}"] & \cM \times \cM
\end{tikzcd}
\]
whose classical truncation recovers $\rZ_{\cE}^r$ by \eqref{eq: Z fibered product}. There is a map $\cZ_{\cE}^r \rightarrow \rA_{\cE}(\F_q)$ as before, and we define $\cZ_{\cE}^r(a)$ to be the fiber of the map $\cZ_{\cE}^r \rightarrow \rA_{\cE}(\F_q)$ over $a$.

\begin{remark}
Recently, Madapusi has defined derived special cycles on Shimura varieties in \cite{Mad23}. A key difference of that situation is that, unlike in the function field case, his construction does not yet have an accompanying moduli-theoretic interpretation. Indeed, the correct moduli problem should be in terms of ``global shtukas'', a notion which has not been defined in the setting of number fields, although recent work of Scholze et al. \cite{Sch18} gives tantalizing hints for its existence. 
\end{remark}

\subsection{Higher theta series}\label{ssec: higher theta series}
Let $n \geq 1$ and $n \geq m \geq 1$. Then we can assemble $[\rZ_{\cE}^r(a)]^{\vir}$ into a Fourier series (with Fourier parameter $a \in \rA_{\cE}(\F_q)$), called a \emph{higher theta series} in \cite{FYZ2}.

Let $m \in \Z_{\geq 1}$. The stack $\Bun_{U^{-}(2m)}$ parametrizes triples $(\cG, h)$, where $\cG$ is a family of rank $2m$ vector bundles on $X'$, and $h: \cG\isom \s^{*}\cG^*$ is a {\em skew-Hermitian} structure (i.e., $\s^{*}h^{*}=-h$).

Let $\Bun_{\wt P_{m}}$ be the moduli stack of quadruples $(\cG,h,\cE)$ where $(\cG, h)\in \Bun_{U^{-}(2m)}$, and $\cE \subset \cG$ is a Lagrangian sub-bundle (i.e., $\cE$ has rank $m$ and the composition $\cE\subset \cG\xr{h}\s^{*}\cG^{*}\to \s^{*}\cE^{*}$ is zero). Thus $\wt{P}_m$ corresponds to the Siegel parabolic of $U^-(2m)$.

Let $\cE' :=\cG/\cE$. The skew-Hermitian form on $\cG$ induces a perfect pairing $\cE\times \s^{*}\cE' \rightarrow \cO_{X'}$, which identifies $\cE'$ with $\s^{*}\cE^{*}$. We thus have a short exact sequence
\begin{equation*}
\xymatrix{0\ar[r] & \cE\ar[r] & \cG\ar[r] & \s^{*}\cE^{*}\ar[r] & 0}
\end{equation*}
and we denote by $e_{\cG,\cE}\in \Ext^{1}(\s^{*}\cE^{*}, \cE)$ its extension class. Thanks to the Hermitian structure on $\cG$, this extension class lies in a subspace $\Ext^{1}_{\Herm}(\s^{*}\cE^{*}, \cE) \subset \Ext^{1}(\s^{*}\cE^{*}, \cE)$. There is a Serre duality pairing between $\Ext^{1}_{\Herm}(\s^{*}\cE^{*}, \cE)$ and $\Hom(\cE, \s^{*}\cE^{\vee})  = \rA_{\cE}(\F_q)$, which we denote by $\tw{e_{\cG, \cE}, a} \in \F_q$. 

Choose a nontrivial additive character $\psi:\F_{q}\to \ol{\Q}^{\times}$. The \emph{higher theta series} is a function of $(\cG,h,\cE) \in \Bun_{\wt P_m}(\F_q)$, assigning to such a triple the class 
\begin{equation}
\Theta^r (\cG,h,\cE) := \chi(\det\cE)q^{n(\deg \cE-\deg\om_{X})/2}\sum_{a\in \rA_{\cE}(\F_q)}\psi(\j{e_{\cG,\cE},a})[\rZ^{r}_{\cE}(a)]^{\vir} \in \CH_{(n-m)r}(\Sht_{U(n)}^r)
\end{equation}
where $\chi: \Pic_{X'}(\F_q)\to \ol{\Q}^{\times}$ is a character whose restriction to $\Pic_{X}(\F_q)$ is the $n$th power of the quadratic character $\Pic_{X}(\F_q)\to \{\pm1\}$ corresponding to the double cover $X'/X$ by class field theory.

\subsection{The modularity conjecture}

By analogy with conjectures of Kudla explained in \cite{Kud04} (whose formulations are themselves partly conjectural), \cite[Conjecture 4.16]{FYZ2} predicts that these higher theta series are \emph{modular} in the sense of automorphic forms. Concretely, this means the following: 

\begin{conj}[Modularity Conjecture]\label{conj: modularity theta}
The function $\Theta^r(\cG, h, \cE)$ in \S \ref{ssec: higher theta series} is independent of the choice of Lagrangian sub-bundle $\cE \subset \cG$.
\end{conj}

\begin{remark}
The reason that this is called a ``Modularity Conjecture'' is that it implies that $\Theta^r(\cG, h, \cE)$ descends to a function on $\Bun_{U^-(2m)}(\F_q)$, which can then be interpreted as a function on the adelic double coset space associated to $U^-(2m)$ using Weil's uniformization. Thus, it implies that $\Theta^r$ defines an automorphic form for $U^-(2m)$, valued in $\CH_*(\Sht_{U(n)}^r)$. 
\end{remark}

When $r=0$, Conjecture \ref{conj: modularity theta} amounts to the modularity of classical theta functions, which has long been known. When $r=1$, it is parallel to the conjectural modularity of arithmetic theta series on orthogonal and unitary Shimura varieties, which has seen much recent progress, as explained in the article of Chao Li \cite{Li24}. But for $r>1$, it is a new phenomenon with no parallel in the classical theory.

We remark that the derived construction of $[\cZ_{\cE}^r(a)]^{\naive} = [\rZ_{\cE}^r(a)]^{\vir}$ already allowed \cite{FYZ2} to prove certain expected properties that were not accessible using the elementary definition. This is not an issue for the virtual classes of special cycles on Shimura varieties; the difference is because of the more complicated nature of singular terms in the function field situation, as described previously. The power of the description $[\cZ_{\cE}^r(a)]^{\naive}$ is that it is more uniform and conceptual, and ``only'' depends on the underlying derived stack $\cZ_{\cE}^r(a)$, whereas $[\rZ_{\cE}^r(a)]^{\vir}$ depends not only on $\rZ_{\cE}^r(a)$ but also on an auxiliary presentation of it. For this reason, it was felt at the time of writing \cite{FYZ2} that the derived perspective would be more useful for proving modularity, although it was not clear how. This view has been at least partially affirmed in \cite{FYZ3} and \cite{FK}, which prove that the $\ell$-adic realization of higher theta series after restriction to the generic fiber of $\Sht_{U(n)}^r \rightarrow (X')^r$ is modular. 

Indeed, \cite{FYZ3} and \cite{FK} are heavily steeped in the formalism of derived algebraic geometry, although the derived construction of the $[\rZ_{\cE}^r(a)]^{\vir}$ is only a small ingredient in the proofs. The modularity on the generic fiber comes from calculations within new forms of Fourier analysis, including a so-called \emph{arithmetic Fourier transform} on Borel-Moore homology groups that restricts to the usual finite Fourier transform in the case $r=0$ (where it pertains to the modularity of the classical theta function), and a \emph{derived Fourier transform} on the derived category of (motivic) sheaves on derived vector bundles. These situations are linked by a higher version of the trace construction, which gives a ``sheaf-cycle correspondence'' generalizing the classical sheaf-function correspondence. Even sketching the argument is beyond the scope of this discussion, so we just refer to \cite[\S 1, \S2]{FYZ3} for an outline. 

\begin{remark}
We anticipate that a statement which will be needed to go beyond modularity in the generic fiber is the \emph{Trace Conjecture} formulated later in \S \ref{ssec: trace conjecture}. It gives yet another approach to the definition of the virtual fundamental classes of special cycles. This conjecture can also be formulated purely within classical algebraic geometry. 
\end{remark}

\section{Derived Galois deformation rings}\label{DGdr}

In this section we will give an overview of the study of derived Galois deformation rings in \cite{GV}. We will elide many technical details in order to give the ``big picture''. Readers craving more details could consult the paper \cite{GV}; another resource is the Oberwolfach Arbeitsgemeinschaft Report \cite{Arb21}. Readers should be familiar with the classical theory of Galois deformation rings before reading this section.

\subsection{Overview}\label{ssec: deformation ring overview}

To orient the reader, we begin with an overview of the contents of this section. 

\emph{Step One.} First, we will explain how to construct the derived global unrestricted Galois deformation ring $\cR$ which represents the functor $\cM_{G}^{\Gal}$ from Example \ref{ex: gal def 3}. \footnote{This requires some assumptions on $\ol{\rho}$, and as literally formulated in Example \ref{ex: gal def 3}, can only exist  if $G$ is of adjoint type, although this can be generalized in a standard manner by accounting for determinants and centers.}

\begin{remark} According to a folklore conjecture of Mazur, $\cR$ should actually be isomorphic to a classical commutative ring, i.e., have vanishing higher homotopy groups. However, this conjecture is wide open, and treating $\cR$ as a derived ring allows us to ``circumvent'' the conjecture for some purposes. Indeed, the key fact we need about $\cR$ is the description of its tangent complex in terms of Galois cohomology, which we can calculate unconditionally for the derived ring $\cR$, for the reasons explained in \S \ref{ssec: geometric interpretation}, while knowing the tangent complex of $\pi_0(\cR)$ is equivalent to Mazur's conjecture. 
\end{remark}

\emph{Step Two.} The deformation ring $\cR$ from Step One does not capture local conditions which we expect to be satisfied for the Galois representations corresponding to automorphic forms. The second step is to cut down the deformation ring by asking for local conditions that ``match'' those seen from automorphic forms. The subtlest aspect of these local conditions concerns the decomposition group at $p$, where the local conditions are governed by $p$-adic Hodge theory.

To do this, one considers analogous derived deformation functors for local Galois groups, and then imposes local conditions by forming the \emph{derived} fibered product against the global deformation functor. We note that the local deformation functors are usually discrete, so that it is the step of forming the derived fibered product that introduces genuine derived structure. 

\begin{remark}
For the local condition at $p$, we only consider what may be considered the simplest situation from a modern perspective on integral $p$-adic Hodge theory: the crystalline Galois deformation functor in the Fontaine-Laffaille range. This should be regarded as a first test case; going beyond it seems to be the first obvious barrier to generalization. 
\end{remark}

This step results in an object $\cR^{\crys}$ called the \emph{derived crystalline global Galois deformation ring}. In many situations it has virtual dimension $-\delta$, where $\delta >0$ is the ``defect'' and therefore cannot be classical. 

We ``know'' the tangent complex of $\cR^{\crys}$ in terms of Galois cohomology, but the structure of its homotopy groups is mysterious a priori. It will be determined in the next step. 

\emph{Step Three.} As a derived generalization of $R=T$ theorems, one expects that $\cR^{\crys}$ acts on the cohomology of locally symmetric spaces, in a way that makes said cohomology ``free'' over $\cR^{\crys}$, at least generically. This action is constructed at the level of homotopy groups in \cite{GV} in tandem with the determination of $\pi_*( \cR^{\crys})$, using the Calegari-Geraghty modification of the Taylor-Wiles method. Thus, the calculation of $\pi_*( \cR^{\crys})$, which is in principle a purely Galois-theoretic problem, uses the correspondence with automorphic forms. Also, this calculation offers an explanation for the interesting numerical pattern in the multiplicities in the cohomology of locally symmetric spaces. This numerical pattern was one of the original motivations for \cite{GV} and \cite{V19}.

\subsection{Global unrestricted derived deformation ring}

Now we will circle back to the beginning, and discuss how to construct the derived Galois deformation functor in Step One of \S \ref{ssec: deformation ring overview}.

\subsubsection{Classical deformation ring} 

Let $k$ be a finite field. For the moment let $\Gamma$ be a discrete group and $G$ a split algebraic group over the Witt vectors $W(k)$.

Let $\ol{\rho}$ be a representation of $\Gamma$ in $G(k)$. Let $\rM_{G}^{\Gamma}$ be the functor on the category of Artinian local $W(k)$-algebras $\rA$ augmented over $k$, sending $\rA$ to the set of lifts modulo conjugacy, 
\begin{equation}\label{eq: classical deformation problem}
\left\{\begin{tikzcd}
& G(\rA)\ar[d] \\
\Gamma \ar[r, "\ol{\rho}"']  \ar[ur, "\rho", dashed] & G(k)
\end{tikzcd} \right\}/\sim
\end{equation}
If $\ol{\rho}$ is absolutely irreducible and $\Gamma$ satisfies certain finite conditions on its group cohomology, then Schlessinger's criterion implies that $\cM_{G}^{\Gamma}$ is representable. These observations, as well as the original idea to consider moduli of Galois deformations, are due to Mazur \cite{Maz89}. The story can be extended to profinite groups $\Gamma$, by treating $\rM_{G}^{\Gamma}$ as the (filtered) colimit of the deformation functors obtained from the finite quotients of $\Gamma$. 

Let $\cO_F$ be the ring of integers in a global field $F$ and $S$ is a finite set of places of $F$. Then $\Gamma = \pi_1(\Spec \cO_F[1/S])$ satisfies the requisite finiteness conditions by class field theory. We are actually going to rename $\rM^{\ol{\rho}}_{\cO_F[1/S]} :=  \rM_{G}^{\Gamma}$ in this case, because this notation emphasizes the parameters that we will focus on in this section. The pro-representing ring is the \emph{(classical) global unrestricted Galois deformation ring}, denoted $\rR_{\cO_F[1/S]}^{\ol{\rho}}$.

\subsubsection{The work of Galatius-Venkatesh}
The idea of Galatius-Venkatesh is to construct a \emph{derived} version of the Galois deformation ring by repeating this story at a simplicial level. The paper \cite{GV} is written in the language of simplicial commutative rings and model categories, rather than the language of animated commutative rings and $\infty$-categories. We will take this opportunity to discuss some of the subtleties that this comparison exposes. 

Morally, we want to upgrade $\rM^{\ol{\rho}}_{\cO_F[1/S]} $ to a functor $\cM^{\ol{\rho}}_{\cO_F[1/S]} $ on the category of \emph{animated Artinian local $W(k)$-algebras $\cA$ augmented over $k$} (to be appropriately defined) which ``sends'' such an $\cA$ to the \emph{anima} of lifts modulo conjugacy, 
\begin{equation}
\left\{\begin{tikzcd}
& G(\cA)\ar[d] \\
\Gamma \ar[r, "\ol{\rho}"']  \ar[ur, "\rho", dashed] & G(k)
\end{tikzcd} \right\}/\sim
\end{equation}

Because of the elaborate nature of the definition of functors between  $\infty$-categories, at least in the formalism of \cite{Lur09}, it can be difficult to rigorously construct functors. In particular, a functor of $\infty$-categories \emph{cannot} be specified by saying where objects and morphisms go (a point which is sometimes handled sloppily in the literature), and so the discussion of the previous paragraph does not really define a valid functor. 

Galatius-Venkatesh chose to use the \emph{model category of simplicial commutative rings} to address this difficulty. One can think of this as a rigidification of the $\infty$-category of animated commutative rings. It is an ordinary 1-category, so a functor can be specified in the usual way on objects and morphisms. However, part of the structure of a model category specifies a notion of \emph{weak equivalence}, alias ``homotopy equivalence'', and we are only interested in functors that send homotopy equivalences to homotopy equivalences. We call such functors \emph{homotopy invariant}. 

To relate this to the other perspective: the $\infty$-category of animated commutative rings is the localization of simplicial commutative rings at the homotopy equivalences, so a homotopy-invariant functor out of simplicial commutative rings should be equivalent to a functor out of animated commutative rings. Such statements can be made precise, but the approach of \cite{GV} is to use the model-theoretic framework throughout.

\subsubsection{Derived deformation ring}\label{sssec: derived global deformation functor} 

Now we summarize the construction of $\cR^{\ol{\rho}}_{\cO_F[1/S]}$. We begin by deriving the functor $\rM^{\ol{\rho}}_{\cO_F[1/S]}$. First we need to define the notion of an \emph{Artinian local simplicial commutative ring}. 

\begin{defn}
A simplicial commutative ring $\cA$ is \emph{Artinian local} if $\pi_0(\cA)$ is Artinian local and $\pi_*(\cA)$ is finitely generated as a module over $\pi_0(\cA)$. 
\end{defn}

Next, we roughly want to define $\cM^{\ol{\rho}}_{\cO_F[1/S]}(\cA)$ to be the \emph{simplicial set} of lifts of $\ol{\rho} \co \Gamma \rightarrow G(k)$ to $\rho \co \Gamma \rightarrow G(\cA)$, up to equivalence. We need to define this simplicial set, and ensure that its construction depends on $\cA$ in a homotopy-invariant way. For this, it is useful to adopt a different perspective on representations. 

Given a discrete group $H$, there is a \emph{classifying space} $B H$. The classifying space $BH$ can be defined as the geometric realization of the simplicial set $N_{\bu}(H)$ obtained as the nerve of $H$ viewed as a category (consisting of a single object, with endomorphisms given by $H$). The simplicial set $N_{\bu}(H)$ takes the form 
\begin{equation}\label{eq: classifying space}
	* \leftleftarrows H  \leftthreearrow H \times H \ldots
\end{equation}
where explicit formulas for face and degeneracy maps can be found in \cite[p.128]{May99}. In particular, $BH$ is equipped with a distinguished basepoint $\pt$. 

If $X$ is a ``nice'' space (e.g., the geometric realization of a simplicial set), then homotopy classes of maps from $X$ to $BH$ are in bijection with ``$H$-local systems'' on $X$. If $X$ is connected, then after choosing a basepoint $x \in X$ these are in bijection with the set of homomorphisms $\pi_1(X,x) \rightarrow H$ modulo conjugacy. Therefore, the \emph{space} of maps from $X$ to $BH$ gives a ``\emph{space} of local systems of $X$''. (See \cite[C.4]{Arb21} for more discussion.)

More generally, the formula \eqref{eq: classifying space} can be used to extend the definition of classifying spaces to groups that are not discrete. If $H$ is a simplicial group, then \eqref{eq: classifying space} is a bi-simplicial set, and we define the simplicial set $N_{\bu}(H)$ to be the diagonal bi-simplicial set. This generalization will be important in defining derived group representation functors, where $H$ will arise by evaluating an algebraic group $G$ on a simplicial commutative ring $\cR$ -- it is the passage from discrete $R$ to simplicial commutative $\cR$ which introduces a non-trivial simplicial structure on $H$.

\begin{example} Suppose $\Gamma$ is a discrete group, and shift perspectives so that $B \Gamma$ and $BH$ are viewed as simplicial sets. Then by general categorical considerations, there is a natural simplicial set of morphisms $B \Gamma \rightarrow B H$, which should be thought of as the simplicial set of representations of $\Gamma$ in $H$. Note that this formalism already incorporates ``modding out by conjugation''. If we want to model \emph{homorphisms} rather than representations, then we take the \emph{pointed} morphisms $(B \Gamma, \pt) \rightarrow (B H, \pt)$. (This is justified by \cite[\S 5.2.6.10, 5.2.6.13]{HA}, which says that the classifying space functor from simplicial groups to pointed spaces is fully faithful.) 
\end{example}

We have now defined a natural candidate for the ``simplicial set of representations $\Gamma \rightarrow G(\cA)$'' (modulo conjugation): it should be the simplicial set of morphisms $B \Gamma \rightarrow BG(\cA)$, where if $\Gamma \cong \limit_i \Gamma_i$ is profinite then we treat $B \Gamma$ as the formal inverse limit of $B\Gamma_i$. 

The actual approach of \cite{GV} is slightly different. To explain it, we recall that from a scheme $X$ one can extract, following Artin-Mazur \cite{AM86} and Friedlander \cite{Fr82}, an \emph{\'{e}tale homotopy type} $\Et(X)$. This is a pro-simplicial set whose purpose is to capture the \'{e}tale topology of $X$. For $\cO_F[1/S]$ the ring of $S$-integers in a global field $F$, it turns out that $\Et(\cO_F[1/S]) \cong B \Gamma$ where $\Gamma$ is the \'{e}tale fundamental group of $\Spec \cO_F[1/S]$ with respect to some basepoint. (Thus $\Spec \cO_F[1/S]$ is an ``\'{e}tale $K(\pi, 1)$'', at least with $p$-adic coefficients where $p>2$.) We prefer $\Et(\cO_F[1/S]) $ to $B \Gamma$, even though they are equivalent, because the former does not refer to a choice of basepoint.\footnote{There is a similarity to the discussion of the derived moduli space of Betti local systems in \S \ref{sec: derived moduli}, in which $\Gamma$ is analogous to $\pi_1(C)$, and $\Et(\cO_F[1/S])$ is analogous to $C$ itself.} By the preceding discussion, the moduli description of $\rM_{G}^{\Gamma}(A)$ can be reformulated as the set of lifts
\[
\begin{tikzcd} 
& BG(A) \ar[d] \\
\Et(\cO_F[1/S]) \ar[r, "\ol{\rho}"']  \ar[ur, "\rho", dashed] & BG(k)
\end{tikzcd}
\]

We then define $\cM^{\ol{\rho}}_{\cO_F[1/S]}$ to be the functor on the category of  Artinian local simplicial commutative rings $\cA$ equipped with an augmentation over $k$, which sends $\cA$ to the simplicial set of ``lifts'' 
\[
\begin{tikzcd} 
& BG(c(\cA)) \ar[d] \\
\Et(\cO_F[1/S]) \ar[r, "\ol{\rho}"']  \ar[ur, "\rho", dashed] & BG(k)
\end{tikzcd}
\]
where $c(\cA)$ is a functorial cofibrant replacement of $\cA$. This notion of ``cofibrant replacement'' is another piece of the model category structure, analogous to projective resolutions in homological algebra, and is what ensures that the functor is actually homotopy invariant.\footnote{A subtle point is that $c(-)$ is not monoidal, so by $BG(c(\cA))$ we mean the cosimplicial commutative ring 
\[
W(k) \rightrightarrows c(\mathcal{O}_G) \rightthreearrow c(\mathcal{O}_G \otimes \mathcal{O}_G) \ldots,
\]
where the maps are dual to those in \eqref{eq: classifying space}.} Here by ``lift'' we mean by definition the \emph{homotopy} fibered product
\[
\Map(\Et(\cO_F[1/S]), BG(c(\cA))) \htimes_{\Map(\Et(\cO_F[1/S]), BG(k))} \{\ol{\rho}\}
\]
which is computed by taking fibered product after forming fibrant replacements in the model category of simplicial sets. 

\subsubsection{Representability}
Lurie's thesis \cite{Lur04} establishes a derived version of Schlessinger's representability criterion. To check the criterion, one needs to compute the tangent complex; the answer was mentioned in Example \ref{ex: gal def 3} and will be revisited later. 

Using Lurie's derived Schlessinger criterion, Galatius-Venkatesh show that $\cM^{\ol{\rho}}_{\cO_F[1/S]}$ is representable by a pro Artinian animated ring $\cR^{\ol{\rho}}_{\cO_F[1/S]}$. Checking the criterion is not the most interesting aspect of \cite{GV}, in our opinion, and we will skip it entirely.

\subsection{Desiderata}
We abstract out the properties of $\cM^{\ol{\rho}}_{\cO_F[1/S]}$ and $\cR^{\ol{\rho}}_{\cO_F[1/S]}$ that we will actually need in the future.

\subsubsection{Tangent complex}

As already discussed in Example \ref{ex: gal def 3}, the tangent complex of $\cM^{\ol{\rho}}_{\cO_F[1/S]}$ at $\ol{\rho}$ is $\rR\Gamma(\Spec \cO_F[S^{-1}], \mf{g}_{\ol{\rho}}[1])$. This is explained in \cite[Example 4.38]{GV}; note the close similarity to derived moduli of Betti local systems, discussed in \S \ref{sec: derived moduli}.

\subsubsection{Comparison to the classical deformation functor}\label{sssec: comparison to classical} We have the following compatibility with the classical theory. 

As a variant of \S \ref{ssec: visualizing derived schemes}, one has a fully faithful embedding of the category of Artinian local rings $A$ augmented over $k$ into the category of Artinian local animated rings augmented over $k$. By construction, we have a commutative diagram
\[
\begin{tikzcd}[column sep = large]
\left\{\substack{\text{Artinian local augmented} \\ \text{commutative rings $A \rightarrow k$}}  \right\} \ar[d, hook] \ar[r, "\rM^{\ol{\rho}}_{\cO_F[1/S]}"] & \msf{Set} \\
\left\{\substack{\text{Artinian local augmented simplicial } \\ \text{commutative rings $\cA \rightarrow k$}}\right\}   \ar[r, "\cM^{\ol{\rho}}_{\cO_F[1/S]}"] & \msf{sSet} \ar[u, "\pi_0"]
\end{tikzcd}
\]
At the level of representing rings, the commutativity of the diagram is equivalent to an isomorphism 
\[
\pi_0(\cR^{\ol{\rho}}_{\cO_F[1/S]}) \cong \rR^{\ol{\rho}}_{\cO_F[1/S]}.
\]

\begin{remark}
As mentioned previously, it is actually a folklore Conjecture, attributed to Mazur, that $\rR^{\ol{\rho}}_{\cO_F[1/S]}$ is LCI of dimension $\chi(\rR\Gamma(\cO_F[1/S], \mf{g}_{\ol{\rho}}[1]))$. By \S \ref{ssec: classicality}, this is equivalent to the statement that $\cR^{\ol{\rho}}_{\cO_F[1/S]}$ is actually homotopy discrete, i.e., $\pi_i(\cR^{\ol{\rho}}_{\cO_F[1/S]}) = 0$ for $i>0$. 

Why then bother with all the effort to define $\cR^{\ol{\rho}}_{\cO_F[1/S]}$, when we could have simply used $\rR^{\ol{\rho}}_{\cO_F[1/S]}$? The most important fact we will use going forward is the description of the tangent complex of $\cR^{\ol{\rho}}_{\cO_F[1/S]}$, which comes relatively easily. Without knowing Mazur's conjecture, we don't know anything about the tangent complex of $\rR^{\ol{\rho}}_{\cO_F[1/S]}$ other than its 0th cohomology group. 
\end{remark}

\subsection{Imposing local conditions} Now we begin the consideration of imposing local conditions. We will use this in two ways: 
\begin{enumerate}
\item To impose conditions of $p$-adic Hodge-theoretic nature, which are necessary in order to cut down to Galois representations that could possibly ``match'' the automorphic side. 
\item To implement the Taylor-Wiles method, wherein one needs to consider the effect of adding and removing ramification at auxiliary primes. 
\end{enumerate}

\subsubsection{Local deformation functors} 
For a local field $F_v$ and a representation $\ol{\rho}$ of $\pi_1(F_v)$, one defines the \emph{local (unrestricted) derived deformation functor} $\cM^{\ol{\rho}}_{F_v}$ analogously to the global case, sending 
\[
\cA \mapsto \Map(\Et(F_v), BG(c(\cA))) \htimes_{\Map(\Et(F_v), BG(k))} \{\ol{\rho}\}.
\]

If $v$ is such that $\ol{\rho}|_{\pi_1(F_v)}$ is unramified, then it is also useful to have the variant $\cM^{\ol{\rho}}_{\cO_{F_v}}$ parametrizing \emph{unramified} deformations, which sends  
\[
\cA \mapsto \Map(\Et(\cO_{F_v}), BG(c(\cA))) \htimes_{\Map(\Et(\cO_{F_v}), BG(k))} \{\ol{\rho}\}.
\]
When these functors come up, they will typically not be representable, since $\ol{\rho}$ will typically not be irreducible. 

\subsubsection{Cutting out local conditions}
Suppose $v$ is a place of $F$. Then we have a map 
\[
\cM^{\ol{\rho}}_{\cO_F[1/S]} \rightarrow \cM^{\ol{\rho}}_{F_v} 
\]
obtained by pulling back a local system along $\Spec F_v \rightarrow \Spec \cO_F[1/S]$. If $v \notin S$, then the above map factors naturally through $\cM^{\ol{\rho}}_{\cO_F[1/S]} \rightarrow \cM^{\ol{\rho}}_{\cO_{F_v}} $. 

Now suppose that we are given some functor $\cM^{\ol{\rho}, \rD_v}_{F_v} \rightarrow \cM^{\ol{\rho}}_{F_v}$, where the $\rD$ stands for some deformation problem. Typically $\cM^{\ol{\rho}, \rD_v}_{F_v} \rightarrow \cM^{\ol{\rho}}_{F_v}$ should be a closed embedding, meaning a closed embedding on classical truncations. An example could be $\cM^{\ol{\rho}}_{\cO_{F_v}} \inj \cM^{\ol{\rho}}_{F_v}$, which informally speaking cuts out the deformations which are ``unramified at $v$''. Then we can form the derived fibered product of the diagram 
\[
\begin{tikzcd}
\cM^{\ol{\rho}, \{D_v\}}_{\cO_F[1/S]} \ar[r] \ar[d] & \cM^{\ol{\rho}}_{\cO_F[1/S]}  \ar[d] \\
\cM^{\ol{\rho}, \rD_v}_{F_v} \ar[r] &  \cM^{\ol{\rho}}_{F_v}
\end{tikzcd}
\]

\begin{example}[Adding ramification]\label{example: adding ramification}
Suppose $v \notin S$, and let $S' = S \sqcup \{ v\}$. The morphism $\Spec \cO_F[1/S'] \rightarrow \Spec \cO_F[1/S]$ induces a map $\cM^{\ol{\rho}}_{\cO_F[1/S]} \rightarrow \cM^{\ol{\rho}}_{\cO_F[1/S']}$. Informally, this morphism sends a deformation unramified outside $S$ to the same deformation regarded as unramified outside $S'$ (i.e., forgetting that it is unramified at $v$).  

In the theory of classical deformation rings one has that 
\[
\rM_{\cO_F[1/S']}^{\ol{\rho}} \times_{\rM_{F_v}^{\ol{\rho}}}  \rM_{\cO_{F_v}}^{\ol{\rho}} \cong \rM_{\cO_F[1/S]}^{\ol{\rho}}.
\]
Informally this amounts to the statement that ``a deformation of local system on $\cO_F[1/S']$ that is unramified at $v$ comes uniquely from a deformation of a local system on $\cO_F[1/S]$'', which is obvious. However, we single out this statement because (i) it is used crucially in the Taylor-Wiles method to descend from patched deformation rings, and (ii) its derived version is non-trivial. 

To follow up on point (ii): it turns out that the analogous identity holds at the level of derived deformation functors:
\begin{equation}\label{eq: derived remove ramification}
\cM_{\cO_F[1/S']}^{\ol{\rho}} \htimes_{\cM_{F_v}^{\ol{\rho}}}  \cM_{\cO_{F_v}}^{\ol{\rho}} \cong \cM_{\cO_F[1/S]}^{\ol{\rho}}.
\end{equation}
This is proved in \cite{GV} by calculating the map of tangent complexes. Roughly speaking, the content of \eqref{eq: derived remove ramification} is the statement that 
\[
\begin{tikzcd}
\Spec F_v \ar[d] \ar[r] & \Spec \cO_F[1/S'] \ar[d]  \\
\Spec \cO_{F_v} \ar[r] & \Spec \cO_F[1/S]
\end{tikzcd}
\]
is a homotopy pushout at the level of \'{e}tale homotopy types. This is similar to \eqref{eq: C homotopy pushout}, where we encountered a presentation with the special property that it was not just a pushout but a homotopy pushout.\footnote{To appreciate this point, it may be instructive to observe that if one uses a presentation of the fundamental group of a genus $g$ Riemann surface which is different from that in \eqref{eq: pi_1 presentation}, then one arrives at an a priori different moduli space of Betti local systems. Why is the presentation in \eqref{eq: pi_1 presentation} the ``right'' one?}

\end{example}

\subsubsection{Crystalline conditions}

The next step is to cut out a derived Galois deformation ring with a $p$-adic Hodge theory condition imposed at $p$. For concreteness, we focus our attention on the ``crystalline'' example, which should correspond to motives with good reduction at $p$, or from another perspective, automorphic forms with no level at $p$. 

We have a map $\cM_{\cO_F[1/S]}^{\ol{\rho}} \rightarrow \cM_{F_p}^{\ol{\rho}}$. (We now know that $\cM_{F_p}^{\ol{\rho}}  = \rM_{F_p}^{\ol{\rho}} $ is classical for $G = \GL_n$ by work \cite{BIP}, and it is expected to be true in general, although we will not logically use this.) We want to define a local functor $\cM^{\ol{\rho},\crys}_{F_p} \rightarrow \cM_{F_p}^{\ol{\rho}}$ that ``cuts out'' the ``crystalline'' deformations. This is something we do not know how to do in general, since we cannot formulate a moduli-theoretic notion of crystallinity over sufficiently general rings (this will be discussed further in \S \ref{ssec: local conditions}).

The only solution to this problem that we know of (at the present time) is to restrict ourselves to local conditions which ``should'' actually be classical, meaning ``$\cM_{\cO_F[1/S]}^{\ol{\rho}, \crys} = \rM_{\cO_F[1/S]}^{\ol{\rho},\crys}$'' where the right hand side is some classically studied object of $p$-adic Hodge theory, and take as our local condition the composition of the maps
\[
\rM_{\cO_F[1/S]}^{\ol{\rho},\crys} \rightarrow \rM_{\cO_F[1/S]}^{\ol{\rho}} \inj \cM_{\cO_F[1/S]}^{\ol{\rho}}.
\]
The way that one can recognize the deformation conditions which ``should'' be classical is, according to \S \ref{ssec: classicality}, that they are LCI of the ``correct'' dimension. In practice, to analyze the resulting objects one also needs to know the tangent complexes, which is typically hard to calculate unless the deformation functor is smooth, so that the tangent complex is just the tangent space. 

A general class of conditions where everything works is the \emph{Fontaine-Laffaille range}, which is what \cite{GV} studies. Here the problems all go away: 
\begin{itemize}
\item There is a complete moduli-theoretic description, even integrally. 
\item The obstructions vanish, so that the deformation functor is smooth and its tangent complex can be described explicitly. 
\end{itemize}
For the rest of this section, we put ourselves in the same situation. We assume we are in the Fontaine-Laffaille range, so that in particular $p$ is ``large enough'' compared to the Hodge-Tate weights that we are considering. This gives a closed embedding $\rM_{\cO_F[1/S]}^{\ol{\rho},\crys} \inj \cM_{\cO_F[1/S]}^{\ol{\rho}}$ as above. We then define $\cR^{\ol{\rho}, \crys}_{\cO_F[1/S]}$ to represent the derived fibered product $\cM^{\ol{\rho}, \crys}_{\cO_F[1/S]}$ of the diagram 
\[
\begin{tikzcd}
\cM^{\ol{\rho},\crys}_{\cO_F[1/S]} \ar[d] \ar[r] & \cM_{\cO_F[1/S]}^{\ol{\rho}} \ar[d] \\
\rM_{\cO_F[1/S]}^{\ol{\rho},\crys}\ar[r, hookrightarrow] & \rM_{\cO_F[1/S]}^{\ol{\rho}}
\end{tikzcd}
\]
Concretely, $\cR^{\ol{\rho}, \crys}_{\cO_F[1/S]}$ is calculated using the derived tensor product (cf. \S \ref{ssec: derived tensor}). 

A key point is that by taking the derived fibered product, one has control over the resulting tangent complex. Its Euler characteristic can be computed by using Poitou-Tate duality for Galois cohomology of global fields, and the answer is $-\delta$, where $\delta \geq 0$ is a certain quantity called the ``defect'' that we are not going to explicate (for our purposes, its definition can be taken to be the Euler characteristic, but it is more typical in the literature to take a different starting point). The appearance of this quantity $\delta$ is crucial; the fact that the same $\delta$ also manifests on the ``automorphic side'' is the linchpin of the Calegari-Geraghty method. 

\subsection{Relation to the Calegari-Geraghty method}\label{ssec: calegari-geraghty} We are now going to discuss the main results of \cite{GV}, which give a reinterpretation and elaboration of the Calegari-Geraghty method \cite{CG18}, which is itself a modification of the Taylor-Wiles method for ``$\delta>0$'' situations. The reader should be familiar with these methods at least at the level of \cite{C20} in order to get something out of this subsection.

\subsubsection{Conjectural picture}

We are going to abbreviate $\cR_0 := \cR^{\ol{\rho}, \crys}_{\cO_F[1/S]}$ and $\rR_0 := \pi_0(\cR_0)$. Let $\chG$ be the Langlands dual group of $G$. There is a \emph{complex} $\cM_0$ of ``automorphic forms for $\chG$ of level $S$", assembled out of the homology of locally symmetric spaces associated to $\chG$. We will index $\cM_0$ to be in degrees $[-\delta, 0]$, although it might be later (in \S \ref{sec: LSS}) reindexed to be in degrees ``$[-\delta-q_0, -q_0]$''. The crucial point is that this complex has length $\delta$ -- the \emph{same} $\delta$ which appears as the negative of the dimension of $\cR_0$. Finally, we write $M_0 := H_0(\cM_0)$. 

The ``automorphic to Galois'' direction of the Langlands correspondence, attaching Galois representations to automorphic forms, equips $M_0$ with the structure of an $\rR_0$-module. (This may be conjectural, depending on what $\chG$ is.) The expectation is (roughly): 

\begin{quote}
There is an action of $\cR_0$ on $\cM_0$ with the following properties. 
\begin{itemize}
\item It recovers the given action of $\rR_0$ on $M_0$ by taking homotopy/homology groups, as in the diagram below:
\[
\begin{tikzcd}[column sep = tiny]
\cR_0 \ar[d, twoheadrightarrow] & \acts & \cM_0 \ar[d, twoheadrightarrow] \\
\rR_0 & \acts & M_0
\end{tikzcd}
\]
\item It makes $\pi_*(\cM_0)[1/p]$ generically free over $\pi_*(\cR_0)[1/p]$.\footnote{technically, we have defined $\cR_0$ as a pro-ring, so that $\pi_*(\cR_0)$ is a pro-group. Some care must be taken to ``invert $p$'' on such an object.}
\end{itemize}
\end{quote}

The main result of \cite{GV} is the construction of an action of $\pi_*(\cR_0)$ on $\pi_*(\cM_0)$, under certain conjectures and certain Taylor-Wiles assumptions, with the properties above. 

\begin{example}[The case of $\GL_1$]\label{ex: gl_1}
At present, the above picture can be realized completely only for $G = \GL_1$, and even this situation is non-trivial.\footnote{So far we have focused on the case where $G$ is of adjoint type, so our definitions have to be modified slightly in order to accommodate this case.} If $\delta=0$, which happens when $F =\Q$ or a quadratic imaginary field, then all rings are classical, and the statement is a consequence of class field theory. The point is that the classical deformation ring is essentially the (completed) group ring of the abelianization of $\pi_1(\cO_F[1/S])$. Class field theory identifies this abelianization with the class group of $\cO_F[1/S]$, whose adelic description is the corresponding ``locally symmetric space'' (a finite set of points). 

However, in all other cases we have $\delta >0$, which implies:
\begin{itemize}
\item the derived Galois deformation ring is non-classical, 
\item the locally symmetric space is no longer a discrete set of points, 
\end{itemize}
and the hypothesized action goes beyond classical Galois deformation theory. Nevertheless, the desired action of $\cR_0$ on $\cM_0$ can still be constructed, which will be explained in the forthcoming paper \cite{FHMR}. The point is that $\cR_0$ admits a description as a (completed) ``derived group ring'' of the ``derived abelianization'' of $\pi_1(\cO_F[1/S])$. To elaborate on what this means: 
\begin{itemize}
\item The abelianization of a group controls its \emph{set} of morphisms to abelian groups. The \emph{derived abelianization} of a group controls its \emph{anima} of morphisms to animated abelian groups. One formulates this precisely in terms of universal properties, as a left adjoint to the forgetful functor from (animated) abelian groups to (animated) groups. 
\item The group ring of an abelian group can also be characterized in terms of a universal property: it is left adjoint to the functor from commutative rings to abelian groups given by extracting units. Analogously, the ``derived group ring'' construction can be characterized as left adjoint to the an analogous functor from animated rings to animated abelian groups. 
\end{itemize}
The punchline of \cite{FHMR} is that the derived abelianization of $\pi_1(\cO_F[1/S])$ ``is'' the animated abelian group obtained from the corresponding locally symmetric space for $G$ (which has a topological group structure in the special case $G = \GL_1$), and its derived group ring ``is'' the homology chains on this animated abelian group, with its induced ring structure. 

Meanwhile, $\cM_0$ is also the homology chains on the locally symmetric space, and the sought-for action is the tautological one. In particular, $\cM_0$ is a free module of rank one over $\cR_0$ for $G = \GL_1$. 
\end{example}

\subsubsection{Results of Galatius-Venkatesh}\label{sssec: Galatius-Venkatesh}
An informal statement of the main result of \cite{GV} is that, assuming ``standard'' conjectures and hypotheses related to automorphy lifting, and under ``no congruences'' and Fontaine-Laffaille assumptions,
\begin{enumerate}
\item $\pi_* (\cR_0)$ is supported in degrees $0 \leq * \leq \delta$, and is moreover an exterior algebra on $\pi_1(\cR_0)$, which in turn is free of rank $\delta$ over $\rR_0$. 
\item $H_*(\cM_{0})$ carries the structure of a free graded module over $\pi_*(\cR_0)$, extending the usual structure for $*=0$. 
\end{enumerate}

The conjectures alluded to above are about the existence of Galois representations attached to automorphic forms, and local-global compatibility for such representations. The hypotheses are (a stringent form of)\footnote{There should be wide scope for optimization of \cite{GV} from a technical perspective, by incorporating Kisin's methods.} the usual ``Taylor-Wiles'' type conditions on the residual representation. See \cite[\S 6.6, \S 10]{GV} for the precise formulations.

\subsubsection{Connection to the Taylor-Wiles method}
The proof of the main result of Galatius-Venkatesh is based on the Taylor-Wiles method, incorporating the modifications introduced by Calegari-Geraghty. It is an interesting feature that the determination of $\pi_* (\cR_0)$ seems to be closely bound to the Taylor-Wiles method, and we will give a brief and extremely impressionistic sketch of how this works. 

First, we orient ourselves psychologically. We believe there should be an action of $\cR_0$ on $\cM_0$, but here and throughout we will only be able to construct (at least at first) the ``classical shadow'' of this action, as depicted in the diagram below:
\[
\begin{tikzcd}[column sep = tiny]
\cR_0 \ar[d, twoheadrightarrow] & \stackrel{?}\acts & \cM_0 \ar[d, twoheadrightarrow] \\
\rR_0 & \acts & M_0
\end{tikzcd}
\]
Here $\stackrel{?}\acts$ stands for an action that we believe exists, but don't know how to construct. 

The starting point of the Taylor-Wiles method is to consider an auxiliary family of similar situations with additional level structure added at a well-chosen collection of primes $Q_n$. Indeed, letting $Q_n$ be such a set, we define
\begin{itemize}
\item $\cR_n := \cR^{\ol{\rho}, \crys}_{\cO_F[1/SQ_n]}$ and $\rR_n := \pi_0(\cR^{\ol{\rho}, \crys}_{\cO_F[1/SQ_n]})$, which by \S \ref{sssec: comparison to classical} is the classical crystalline deformation ring. The derived ring $\cR_n $ has virtual dimension $-\delta$. 
\item $\cM_n$ for the corresponding module of ``automorphic forms for $G$ of level $SQ_n$'', which by assumption is a complex concentrated in degrees $[-\delta, 0]$, and $M_n := H_0(\cM_n)$.
\end{itemize}
Then the same story holds for each $\cR_n$: we want an action $\cR_n \qacts \cM_n$, but what we can actually construct is its classical shadow $\rR_n \acts M_n$. Therefore, we have a diagram 
\[
\begin{tikzcd}[column sep = tiny]\label{eq: TW level n diagram}
 \cR_n \ar[d, twoheadrightarrow] & \qacts  & \cM_n \ar[d, twoheadrightarrow] \\
  \rR_n & \acts & M_n
\end{tikzcd}
\]
At first, the picture looks the same as when $n=0$, so it seems that no progress has been made. But the point is that by choosing $Q_n$ artfully, one can \emph{morally} arrange the $\cR_n$ to ``limit'' to an object $\cR_{\infty}$ which is actually \emph{classical}, so that $\cR_\infty \xrightarrow{\sim} \rR_{\infty} :=  \pi_0(\cR_{\infty})$. Now, in reality we will not be constructing $\cR_{\infty}$ in this way, but we will construct an $\rR_{\infty}$ which will be seen to be ``smooth and of the correct dimension''. By the discussion in \S \ref{ssec: classicality}, this justifies setting $\cR_{\infty}$ to be equal to $\rR_{\infty}$. 

Furthermore, the $\cM_{n}$ ``limit'' to an object $\cM_{\infty}$ (which really does exist), and it is really true that $\cM_{\infty} \xrightarrow{\sim} M_{\infty}$. Also, one should imagine that the diagram \eqref{eq: TW level n diagram} ``limits'' to a diagram 
  \[
\begin{tikzcd}[column sep = tiny]\label{eq: TW level infty diagram}
 \cR_\infty \ar[d, twoheadrightarrow] & \qacts  & \cM_\infty \ar[d, twoheadrightarrow] \\
  \rR_\infty & \acts & M_\infty
\end{tikzcd}
\]
Since in this case $ \cR_\infty \surj \rR_{\infty}$ is actually an isomorphism of simplicial commutative rings, it automatically gives the desired action on $\cM_{\infty} \xrightarrow{\sim} M_{\infty}$.

To see how this helps with our original problem, let us contemplate the relation between $\cR_n$ and $\cR_0$. The difference is that $\cR_0$ only parametrizes deformations unramified at $Q_n$ (recall that $\ol{\rho}$ is unramified at the primes in $Q_n$), while $\cR_n$ allows deformations that ramify at $Q_n$. As explained in Example \ref{example: adding ramification}, $\cR_0$ can be recovered from $\cR_n$ by forming a derived tensor product that cuts out the unramified locus from among the deformations of local Galois groups over the primes in $Q_n$. More precisely, one chooses $Q_n$ so that deformations of $\ol{\rho}$ are automatically tame above $Q_n$. Then the maps from the tame inertia groups over $Q_n$ to $\pi_1(\cO_F[1/SQ_n])$ induce a ring homomorphism $\rS_n \rightarrow \cR_n$, where $\rS_n$ is certain deformation ring for tame inertia at the primes in $Q_n$. There is an augmentation $\rS_n \rightarrow \rS_0 = W(k)$ that corresponds to the deformations which are trivial on tame inertia. By similar conditions as in Example \ref{example: adding ramification}, one has an equivalence of simplicial commutative rings
\begin{equation}\label{eq: derived tensor down}
\cR_n \dotimes_{\rS_n} \rS_0 \xrightarrow{\sim} \cR_0.
\end{equation}
This identity \emph{approximately} ``limits'' to 
\begin{equation}\label{eq: R infinity to 0}
``\cR_\infty \dotimes_{\rS_\infty} \rS_0 \xrightarrow{\sim} \cR_0."
\end{equation}
More generally, one has \emph{approximately} 
\begin{equation}\label{eq: R infinity to n}
``\cR_\infty \dotimes_{\rS_\infty} \rS_n \xrightarrow{\sim} \cR_n."
\end{equation}
The reason we have put these statements in quotes is that they are not exactly true, but they are true up to a certain amount of error which shrinks to $0$ as $n \rightarrow \infty$, which makes them acceptable for our purposes. The source of this ``error'' lies in the construction of $\rR_{\infty}$: it is extracted by a compactness argument from an inverse system of \emph{finite} quotients of the $\rR_n$. 

There will also be an action of $\rS_n$ on $M_n$ by ``diamond operators'', for which one has 
\begin{equation}\label{eq: M infinity to 0}
\cM_{\infty} \dotimes_{\rS_\infty} \rS_0 \xrightarrow{\sim} \cM_0,
\end{equation}
Furthermore, local-global compatibility ensures that this action is compatible with the one induced by $\rS_n \rightarrow \rR_n$. Then, from the $\cR_\infty \xrightarrow{\sim} \rR_{\infty}$ action on $M_{\infty} \xleftarrow{\sim} \cM_{\infty}$, taking derived tensor products and using \eqref{eq: R infinity to 0} and \eqref{eq: M infinity to 0}, we have  
\[
``\cR_0 \cong \cR_\infty \dotimes_{\rS_\infty} \rS_0"  \acts \cM_{\infty}  \dotimes_{\rS_\infty} \rS_0   \cong \cM_0.
\]
As discussed above, the statement in quotes is not exactly true but can be made true on homotopy groups up to arbitrary precision. This is done by instead using \eqref{eq: R infinity to n} to ``go from $\infty$ to $n$'', which incurs arbitrarily small error as $n$ increases, and then using \eqref{eq: derived tensor down} to ``go from $n$ to $0$'', which does not incur any error. That is enough to make some isomorphism 
\[
\pi_i(\cR_0 ) \approx \pi_i( \cR_\infty \dotimes_{\rS_\infty} \rS_0) \cong \Tor^{\rS_\infty}_i(\cR_{\infty}, \rS_0),
\]
which we will then be able to compute explicitly (under favorable assumptions). Indeed, the argument also shows that $\rS_{\infty} \rightarrow \rR_{\infty}$ is a quotient by a regular sequence of length $\delta$, which allows to compute $\Tor^{\rS_\infty}_i(\cR_{\infty}, \rS_0)$ by a standard Koszul complex calculation. 

This completes our shamelessly vague and impressionistic sketch. We encourage the reader to consult the article of Caraiani-Shin \cite{CS23} for a more substantial discussion of the Calegari-Geraghty enhancement of the Taylor-Wiles method.

\section{Derived Hecke algebras}\label{sec: DHA}
\subsection{The local derived Hecke algebra}\label{ldHa}

We briefly review the theory of derived Hecke algebras from \cite{V19}. Let $G$ be a split reductive group over a local field $K$, say of residue field $\F_q$ having characteristic $p$. Let $U \subset G(K)$ be a compact open subgroup. (For our purposes, we are most interested in the maximal compact subgroup $U = G(\cO_K)$.) 

Let $\Lambda$ be a commutative ring. Define the universal module
$$\mf{U}_\Lambda(U) = \cInd_{U}^{G(K)} \Lambda$$
where $U$ acts trivially on $\Lambda$. We can present the usual Hecke algebra for the pair $(G(K),U)$ as 
\[
\rH(G(K),U; \Lambda) := \Hom_{G(K)}(\mf{U}_\Lambda(U),\mf{U}_\Lambda(U)).
\]
This presentation suggests the following generalization. 

\begin{defn}
The \emph{derived Hecke algebra} for $(G(K),U)$ with coefficients in $\Lambda$ is the graded ring
\[
\Cal{H}(G(K),U;\Lambda) := \Ext^*_{G(K)}(\mf{U}_\Lambda(U),\mf{U}_\Lambda(U)),
\]
where the $\Ext$ is formed in the category of smooth $G(K)$-representations. For $U= G(\cO_K)$, we abbreviate $\Cal{H}(G(K);\Lambda) := \Cal{H}(G(K), G(\cO_K);\Lambda)$. 
\end{defn}

\begin{remark}
The ring $\Cal{H}(G(K),U;\Lambda)$ is really just a graded (associative) ring, rather than any kind of derived (associative) ring. However, its description makes it clear that it arises as the cohomology groups of a differential graded algebra $\RHom_{G(K)}(\mf{U}_\Lambda(U),\mf{U}_\Lambda(U))$. It might be more puritanical to call the latter object the ``derived Hecke algebra,'' but we follow \cite{V19} in applying this name to the graded ring $\Cal{H}(G(K),U;\Lambda)$. We will instead call $\RHom_{G(K)}(\mf{U}_\Lambda(U),\mf{U}_\Lambda(U))$ the \emph{differential graded Hecke algebra}, and denote it by $\DGHA(G(K), U; \Lambda)$, or sometimes simply $\DGHA(G(K); \Lambda)$ if $U = G(\cO_K)$. It does not have the structure of an animated ring; a priori, it is just a differential graded (i.e., $\EE_1$) algebra, although we conjecture in \S \ref{ssec: DGHA} that it has more commutative structure.
\end{remark}

\begin{example}
The $0$th graded group of $\Cal{H}(G(K),U;\Lambda) $ is the classical Hecke algebra $\rH(G(K),U;\Lambda)$.
\end{example}

We next give a couple more concrete descriptions of the derived Hecke algebra, following \cite[\S 2]{V19}.

\subsubsection{Function-theoretic description}\label{sssec: DHA function-theoretic}

Let $x,y \in G(K)/U$ and $G_{x,y} \subset G(K)$ be the stabilizer of the pair $(x,y)$. We note that $G_{x,y}$ is a compact open subgroup of $G(K)$. We can think of $\Cal{H}(G(K),U; \Lambda)$ as the space of functions 
\[
  G(K)/U \times G(K)/U \ni (x,y) \mapsto h(x,y) \in H^*(G_{x,y};\Lambda)
\]
satisfying the following constraints: 
\begin{enumerate}
\item The function $h$ is ``$G(K)$-invariant'' on the left. More precisely, we have 
\[
[g]^* h(gx,gy) = h(x,y)
\]
for all $g \in G(K)$, where $[g]^* \co H^*(G_{gx,gy}; \Lambda) \rightarrow H^*(G_{x,y}; \Lambda)$ is pullback by $\Ad(g)$. 
\item The support of $h$ is finite modulo $G(K)$. 
\end{enumerate}
The multiplication is given by a convolution formula, where one uses the cup product to define multiplication on the codomain, and restriction/inflation to shift cohomology classes to the correct groups \cite[eqn. (22)]{V19}.

\subsubsection{Double coset description}\label{sssec: double coset description}

For $x \in G(K)/U$, let $U_x = \Stab_U(x)$. Explicitly, if $x = g_xU$ then $U_x :=  U \cap g_x Ug_x^{-1}$.

We can also describe $\Cal{H}(G(K),U;\Lambda)$ as functions 
\[
x \in U \bs G(K) / U \mapsto h(x) \in H^*(U_x; \Lambda)
\]
which are compactly supported, i.e., supported on finitely many double cosets. For the description of the algebra structure, it seems better to convert to the model of \eqref{sssec: DHA function-theoretic}.

\subsubsection{The derived Hecke algebra of a torus}
Let $T$ be a split torus over $K$. We can explicitly describe the derived Hecke algebra of the torus $T(K)$, using now the double coset model. Since $T$ is abelian we simply have $T(\cO_K)_x = T(\cO_K)$ for all $x \in T(K)/T(\cO_K)$. 

We have $T(K)/T(\cO_K) \cong X_*(T)$. Identify
\begin{equation}\label{eq: satake embed}
X_*(T) = T(K) / T(\cO_K) \hookrightarrow G(K)/G(\cO_K)
\end{equation}
by the map $X_*(T) \ni \chi  \mapsto \chi(\varpi) \in G(K)/G(\cO_K)$, where $\varpi$ is a uniformizer of $K$. Then $\Cal{H}(T(K); \Lambda)$ simply consists of compactly supported functions
\[
X_*(T) \rightarrow H^*(T(\cO_K); \Lambda)
\]
with the multiplication given by convolution; in other words,
\[
\Cal{H}(T(K); \Lambda) \cong \Lambda[X_*(T)] \otimes_\Lambda H^*(T(\cO_K); \Lambda).
\]

\begin{remark} If $p$ is invertible in $\Lambda$, then the structure of $ H^*(T(\cO_K); \Lambda)$ can be elucidated as follows. First of all, the reduction map $T(\cO_K) \rightarrow T(\F_q)$ has pro-$p$ kernel, which then has vanishing higher cohomology with coefficients in $\Lambda$ by the assumption, and therefore induces an isomorphism $H^*(T(\F_q); \Lambda) \rightarrow H^*(T(\cO_K); \Lambda)$. 

Furthermore, since $T$ is split we have $T(\F_q) \cong (\F_q^{\times})^r$ where $r = \rank(T)$. We are mainly interested in the case where $\Lambda$ is of the form $\cO/l^n$ where $\cO$ is the ring of integers in a number field and $l$ is a prime above $\ell \neq p$. In this situation, $H^*(T(\F_q); \Lambda)$ has non-vanishing higher cohomology exactly when $q \equiv 1 \pmod{\ell}$. 

\end{remark}

\subsubsection{The derived Satake isomorphism}

We assume in this subsection that $q \equiv 1 \in \Lambda$. Let $T \subset G$ be a maximal torus and $W$ the Weyl group of $T \subset G$. We consider an analogue of the classical Satake transform for the derived Hecke algebra $\Cal{H}(G(K);\Lambda)$, which takes the form 
\[
``\text{Derived Hecke algebra for $G$}  \xrightarrow{\sim} (\text{Derived Hecke algebra for }T)^W."
\]

More precisely, we define the \emph{derived Satake transform}
\begin{equation}\label{eq: satake transform}
  \Cal{H}(G(K);\Lambda) \rightarrow \Cal{H}(T(K); \Lambda)
\end{equation}
simply by \emph{restriction} (in the function-theoretic model \S \ref{sssec: DHA function-theoretic}) along the map 
\[
(T(K)/T(\cO_K))^2 \rightarrow (G(K)/G(\cO_K))^2 
\]
from \eqref{eq: satake embed}. In more detail, let $h \in \Cal{H}(G(K); \Lambda)$ be given by the function 
\[
(G(K)/G(\cO_K))^2 \ni (x,y) \mapsto h(x,y) \in H^*(G_{x,y};\Lambda).
\]
Then \eqref{eq: satake transform} takes $h$ to the composition 
\[
\begin{tikzcd}
(T(K)/T(\cO_K))^2  \ar[r, hook] & (G(K)/G(\cO_K))^2  \ar[r, "h"] & H^*(G_{x,y};\Lambda)  \ar[r, "\res"]  &  H^*(T_{x,y};\Lambda)
\end{tikzcd}
\] 

\begin{remark}It may be surprising that this is the right definition, since the analogous construction in characteristic 0, on the usual (underived) Hecke algebra,  is far from being the usual Satake transform. It is only because of our assumptions on the relation between the characteristics (namely, that $q \equiv 1 \in \Lambda$) that this ``na\"{i}ve'' construction turns out to behave well (otherwise, it would not even be a ring homomorphism). The construction seems to have been motivated by \cite{TV16}; see \cite[Lemma 6.5]{Fen24} for a proof of the relevant point. 
 
\end{remark}

\begin{thm}[{\cite[Theorem 3.3]{V19}}] \label{thm: derived sat isom} Let $W$ be the Weyl group of $T$ in $G$. If $|W|$ is invertible in $\Lambda$ and $q =1 \in \Lambda$, then the map \eqref{eq: satake transform} induces an isomorphism 
\[
\Cal{H}(G(K);\Lambda) \xrightarrow{\sim} \Cal{H}(T(K); \Lambda)^W.
\]
\end{thm}

\begin{remark}
Strictly speaking, Theorem \ref{thm: derived sat isom} is only proved when $K$ has characteristic zero in \cite{V19}, but essentially the same argument should work in general. 
\end{remark}

\subsubsection{Commutativity?}\label{sssec: commutativity}

A priori, $\DGHA(G(K), U; \Lambda)$ is just an associative differential graded algebra and $\Cal{H}(G(K),U;\Lambda)$ is just an associative graded algebra. It is interesting to ask what more structure or properties these algebras have. Indeed, $\rH^0(\sH(G(K),U;\Lambda) ) \cong \rH(G(K),U;\Lambda)$ famously turns out to be commutative, although this is non-trivial to prove. 

We think it is reasonable to conjecture that $\cH(G(K) ;\Lambda)$ is also commutative, at least if the characteristic of $\Lambda$ is not in ``bad'' position with respect to $G$ and the residue characteristic of $K$. (No counterexample is known at present, even for bad characteristics.) We list some cases in which the commutativity is known:
\begin{enumerate}
\item The commutativity is known for $\Lambda$ of characteristic $\ell \nmid \# W$ such that $q \equiv 1 \pmod{\ell}$ by Theorem \ref{thm: derived sat isom} (which is due to Venkatesh). 
\item More generally, the commutativity is known if $\# G/B(\F_q)$ and $p$ are invertible in $\Lambda$ by Gehrmann \cite{G20}.\footnote{This paper also investigates examples for $G = \SL_2$ which do not fall under these hypotheses.} 
\end{enumerate}
To say a bit about these arguments, we recall two proofs of the commutativity of the classical spherical Hecke algebra: one via the Satake transform, and one via ``Gelfand's trick''. Venkatesh's argument is a generalization of the method of the Satake transform, and Gehrmann's argument is a generalization of Gelfand's trick. We note that although Gehrmann proves commutativity in more generality, Venkatesh's proof gives more refined information which is needed for global applications.

There are subtler questions that one can formulate beyond the commutativity of $\cH(G(K); \Lambda)$, concerning the structure of the differential graded Hecke algebra $\DGHA(G(K); \Lambda)$. Some more speculative discussion of this problem will be given in \S \ref{ssec: DGHA}.

\subsubsection{Extensions} The literature on derived Hecke algebras is not as comprehensive as one would hope. It would be useful to:
\begin{itemize}
\item Extend the theory to non-split groups $G$. 
\item Extend the theory to incorporate non-trivial coefficients, i.e., replace $\mf{U}_\Lambda(U)$ with the compact-induction of a non-trivial representation of $U$. In a global setting, this would be related to considering non-trivial local coefficients on the corresponding locally symmetric space. 
\end{itemize}

\subsection{Actions of derived Hecke algebras}\label{actiondHa}
Maintain the preceding notation: let $G$ be a reductive group over a local field $K$ of residue characteristic $p$, $U \subset G(K)$ be an open compact subgroup, and $\Lambda$ be a commutative ring. 

\begin{fact}\label{DHAaction} If $C^\bullet$ is an object of the derived category of smooth $G(K)$-representations over $\Lambda$, then we have (formally) a \emph{right} action of the differential graded Hecke algebra $\DGHA(G(K),U; \Lambda)$ (resp. the derived Hecke algebra $\cH(G(K),U; \Lambda)$) on $\RHom_{\Lambda[U]}(\Lambda,C^{\bullet})$ (resp. $\Ext^*_{\Lambda[U]}(\Lambda,C^{\bullet})$). 
\end{fact}
Indeed, the derived endomorphism ring $\RHom_{G(K)}(\mf{U}_\Lambda(U),\mf{U}_\Lambda(U))$ tautologically has a right action on the functor
$\RHom^{\bullet}_{\Lambda[G(K)]}(\mf{U}_\Lambda(U),-)$. Then the Fact follows from  the Frobenius reciprocity isomorphism
 \begin{equation}\label{FR}
\RHom^\bullet_{\Lambda[U]}(\Lambda,C^{\bullet}) \isoarrow \RHom^{\bullet}_{\Lambda[G(K)]}(\mf{U}_\Lambda(U),C^{\bullet}).
\end{equation}

The functor $C^{\bu} \mapsto \RHom^\bullet_{\Lambda[U]}(\Lambda,C^{\bullet})$ is that of \emph{derived $U$-invariants}. Thus Fact \ref{DHAaction} may be captured more colloquially by the slogan:

\begin{quote} There is a canonical right action of the derived Hecke algebra for $(G,U)$ on the derived $U$-invariants of a $G(K)$-representation.
\end{quote}

\begin{remark}
If $U$ is pro-$p$, and $p$ is invertible in $\Lambda$, then formation of $U$-invariants is exact. The case where $p$ is not invertible in $\Lambda$ will be discussed in the next subsection. 
\end{remark}

\begin{example}
A source of examples comes from the cohomology of locally symmetric spaces. Taking $U$ to be the level structure at $p$, the cohomology complex of the locally symmetric space for $G$ is realized canonically as the derived $U$-invariants of a $G(K)$-representation, which is at least heuristically ``the cohomology complex of the locally symmetric space with infinite level structure at $p$''. The motivation for the local derived Hecke algebra in \cite{V19} was to analyze such examples. We will say more about this later in the subsequent sections.
\end{example}

\subsection{Mod $p$ derived Hecke algebras of $p$-adic groups}

Let $G$ be a reductive group  over $\Q_p$ -- we use the same notation for the algebraic group and its group of $\Q_p$-valued points -- and
let $k$ be an algebraically closed field of characteristic $p$.  Let $\Rep_k(G)$ denote the category of continuous representations of $G$ on
$k$-vector spaces, and let $\Rep^{sm}_k(G) \subset \Rep_k(G)$ denote the subcategory of smooth representations.  When $G = \GL_2(\Q_p)$
or a closely related group, the irreducible objects in $\Rep^{sm}_k(G)$ have been classified for some time, but for every other group even an approximate classification remains elusive, in spite of impressive efforts over the past 15 years including \cite{Br03, Ko, BP}.  

The abelian category $\Rep^{sm}_k(G)$ nevertheless has a simple structure in one respect:  every irreducible object admits a canonical space of surjective
homomorphisms from a compact projective generator, denoted $\mf{U}$ below.  This fact has been exploited by Schneider to define a derived equivalence between $\Rep^{sm}_k(G)$ and the derived category of dg-modules over the derived endomorphism algebra $R\End_G(\mf{U})^{\op}$. The properties of the latter, a mod $p$ analogue of the mod $\ell$ derived Hecke algebras introduced above, are still largely unexplored.  A simpler version, called the {\it derived diamond algebra}, has been shown by Khare and Ronchetti \cite{KR23} to play a significant role in the theory of $p$-adic modular forms.

\subsubsection{Schneider's theory}  Fix an Iwahori subgroup $I \subset G$, and let $I(1) \subset I$ be its maximal pro-$p$ normal subgroup; thus $I/I(1)$ is the group of $\mathbb{F}$-points of a torus over a finite field $\mathbb{F}$. It is easy to show that every irreducible $\pi \subset \Rep^{sm}_k(G)$ is generated by its subspace
$\pi^{I(1)}$. This defines a canonical functor
\begin{align}\label{tau}
\tau_k:  \Rep^{sm}_k(G) &\ra \Mod(\rH(G,I(1);k)^{\op}) \\
 \pi &\mapsto \pi^{I(1)} \nonumber
\end{align}
where $\rH(G,I(1);k)$ is the $k$-valued Hecke algebra of $G$ relative to $I(1)$. Here we caution about a notational inconsistency: Schneider's convention differs from ours (which follows Venkatesh) by formation of opposite algebra. That is, letting $\mf{U} := \cInd_{I(1)}^G k$ be the universal induced module, Schneider defines the Hecke algebra of $G$ relative to $I(1)$ to be $\End_G(\mf{U})^{\op}$ instead of $\End_G(\mf{U})$. Thus, when we translate his results into our conventions, an extra ``op'' will appear.  

A well-known theorem of Borel and Casselman \cite{Bor76} asserts that, if $k$ is replaced by $\CC$, then \eqref{tau} defines an equivalence of categories between the subcategory of the left-hand
side of representations generated by their $I(1)$-fixed vectors and $\Mod(\rH(G,I(1);k)^{\op})$.  This is also true for $k$ in characteristic $p$ when $G = \GL_2(\Q_p)$, but Ollivier \cite{Oll09} proved that this is almost never the case for other $G$, even though every irreducible $\pi$ is generated by $\pi^{I(1)}$.  Schneider found the appropriate
generalization of the Borel--Casselman theorem by replacing $\rH(G,I(1);k)$ by the differential graded algebra
\begin{equation}\label{derivedH}
\DGHA(G,I(1);k) = R\End_G(\mf{U}) := \RHom_G(\mf{U},\mf{U}).
\end{equation}
The replacement of the functor \eqref{tau} is
\begin{align}\label{taub}
\tau^\bullet_k:  D(G) & \ra D(\DGHA(G,I(1);k)^{\op}) \\
  \pi &\mapsto \RHom_{I(1)}(1,\pi) = \RHom_G(\mf{U},\pi) \nonumber
\end{align}
Here $D$ denotes the (unbounded) derived category, we write $D(G)$ instead of $D(\Rep^{sm}_k(G))$, 
and the final equality is Frobenius reciprocity.  Schneider's theorem is the following.

\begin{thm}[\cite{Sc15}, Theorem 9]\label{dgI}  Suppose $I(1)$ is a torsion free $p$-adic group.  Then $\mf{U}$ is a compact generator of $D(G)$ and the functor \eqref{taub} is
an equivalence of triangulated categories.
\end{thm}

A first description of the category $\Rep^{sm}_k(G)$ is contained in \cite{AHHV}.  Parabolic induction is defined as a functor from
$\Rep^{sm}_k(M)$ to $\Rep^{sm}_k(G)$ if $M$ is a Levi component of a parabolic subgroup $P$; in characteristic $p$ this induction is not normalized and the
result depends on $P$ as well as $M$.   The irreducible admissible representation $\pi$ of $G$ is {\it supercuspidal} if it does not occur as a subquotient of a
representation induced from an admissible irreducible representation of a proper parabolic subgroup of $G$.  The main theorem of \cite{AHHV} reduces the classification
of $\Rep^{sm}_k(G)$ to the classification of supercuspidal representations of Levi subgroups.  This corresponds to a similar classification of modules for 
the (underived) Hecke algebra $\rH(G,I(1);k)^{\op}$:  the Hecke algebra modules corresponding to supercuspidal $\pi$ are called {\it supersingular} and have a
simple characterization in terms of $\rH(G,I(1);k)^{\op}$ (see \cite[p. 498]{AHHV}).   

\begin{remark}This suggests that the structure of $D(\DGHA(G,I(1);k)^{\op})$ can also be
reduced to the study of triangulated subcategories of supersingular modules, but no one seems to have looked into this question.  
\end{remark}

When $I(1)$ is not torsion free, it can be replaced in the definition of $\mf{U}$ and in Theorem \ref{dgI} by an appropriate
torsion free subgroup of finite index, but the theory of the underived Hecke algebras for general open compact subgroups is
much less clear.

As in the $\ell \neq p$ case, we will call the cohomology ring 
\[
\cH(G,U; \Lambda) := H^*(\DGHA(G, U; \Lambda))
\]
the \emph{derived Hecke algebra}; this is (a priori) just a graded associative ring. The thesis of N. Ronchetti \cite{Ron} studies the derived {\it spherical} Hecke algebra $\Cal{H}(G,G(\cO);\Lambda)$ and defines a Satake homomorphism
$$\Cal{H}(G,G(\cO);\Lambda) \ra \mathcal{H}(T,T(\cO);\Lambda) \isoarrow \Lambda[X_*(T)]\otimes H^*(T(\cO);\Lambda).$$  
In this case, with $T(\cO)$ $p$-torsion free, $H^*(T(\cO),\Lambda)$ is an exterior algebra on $H^1$.   This algebra is much more manageable
than Schneider's $\DGHA(G,I(1);k)$; however, the compact induction of the trivial representation
of $G(\cO)$ does not generate $D(G)$, so Ronchetti's algebra cannot account fully for the mod $p$ representation theory of $G$.  However, simpler
$p$-adic derived Hecke algebras do seem to play a role in the global theory of $p$-adic modular forms; see \S \ref{dda} below.

\section{Cohomology of locally symmetric spaces}\label{sec: LSS}

\subsection{Locally symmetric spaces}
In this section $G$ is a connected reductive group over $\QQ$, with center $Z = Z_G$.  We fix a subgroup $K_\infty \subset G(\RR)$ containing $Z(\RR)$ and a maximal
compact subgroup of the identity component $G(\RR)^0 \subset G(\RR)$, so that $K_\infty/Z(\RR)$ is compact.  Then $X = G(\RR)/K_\infty$ is a finite union of copies 
of the symmetric space for the derived subgroup $G(\RR)^{\rm der,0}$.  Let $K_f \subset G(\mathbf{A}_f)$ be an open
compact subgroup; the corresponding {\it congruence subgroup} of $G(\RR)$ is the subgroup $\Gamma = \Gamma_{K_f} \subset G(\QQ)$ given by the intersection $G(\QQ) \cap K_f$
in $G(\mathbf{A}_f)$.  Let $K'_f \subset K_f$ be a torsion-free normal subgroup of finite index, with quotient $Q$; then
$$X(\Gamma_{K'_f}) := \Gamma_{K'_f}\backslash X$$
is a $C^\infty$ manifold, and 
$$X(\Gamma_{K_f}) = \Gamma_{K_f}\backslash X = Q\backslash X(\Gamma_{K'_f})$$
is the locally symmetric space attached to $K_f$.  For most purposes the finite group $Q$ and the singularities of $X(\Gamma_{K_f})$ will play no role in what follows.

We will be most concerned with the adelic locally symmetric spaces
$$S_{K_f}(G,X) = G(\QQ)\backslash X \times G(\mathbf{A}_f)/K_f$$
and with the limit
$$S(G,X) = \varprojlim_{K_f}  S_{K_f}(G,X),$$
the  limit taken over all open compact subgroups $K_f \subset G(\mathbf{A}_f)$.  
For fixed $K_f$, the adelic locally symmetric space can be identified with the disjoint union of a finite set of discrete arithmetic quotients of $X$:
by reduction theory we can write $G(\aA_f) = \coprod_j G(\QQ) \alpha_j K_f$ for a finite set of $\alpha_j$, and then 
$$S_{K_f}(G,X) = \coprod_j X(\Gamma_j),$$
where $\Gamma_j = G(\QQ)\cap \alpha_j K_f \alpha_j^{-1}$.

Unless otherwise indicated, for any cohomology theory $H^*$, usually
with twisted coefficients, the
space $H^*(S(G,X))$ will be understood to be the filtered colimit 
$$\varinjlim_{K_f}H^*(S_{K_f}(G,X)).$$

\subsection{Review of the de Rham theory}\label{derham}

Let $\cA(G)$ denote the (complex) vector space of automorphic forms on $G(\QQ)\backslash G(\aA)$, and let
$\cA_0(G) \subset \cA(G)$ denote the subspace of cusp forms.  Let $Z = Z_G$ denote the center of $G$, $K_\infty \subset G(\RR)$
a closed subgroup of the identity component that contains $Z(\RR)$ and whose image in the adjoint group $G(\RR)^{\ad}$ is compact.
Depending on the circumstances, $K_\infty$ may or may not be maximal compact modulo center.  For bookkeeping purposes, we let
$K^c_\infty$ denote the maximal compact subgroup of $K_\infty$.  We write $\fg = \Lie(G)$,
$\fk = \Lie(K_\infty)$; unless otherwise indicated these are taken to be the Lie algebras of the complex points.   Then we have the
Cartan decomposition
\begin{equation}\label{Cartan}   \fg = \fk \oplus \fp.
\end{equation}
The (possibly disconnected) symmetric space $X = G(\RR)/K_\infty$ has an invariant hermitian structure precisely when there is a 
homomorphism
$$h:  \CC^\times = R_{\CC/\RR} \BG_{m,\CC} \ra G_\RR$$
with image in the center of $K_\infty$ that satisfies the axioms of a Shimura datum (cf. Morel's article \cite{Mor23} in this volume). In that case,
\eqref{Cartan} can be refined into the $\fk$-invariant Harish-Chandra decomposition:
\begin{equation}\label{CartanHC}   \fg = \fk \oplus \fp^+ \oplus \fp^{-},
\end{equation}
Here, under the map of tangent spaces 
$$\fg = T_{G,e} \ra T_{X,h},$$
where $T_{G,e}$ (resp. $T_{X,h}$) is the tangent space at the identity (resp. the tangent space at the $K_\infty$-fixed point $h \in X$), 
 $\fp^+$  (resp. $\fp^-$) is taken to the holomorphic (resp. antiholomorphic) tangent space; they are respectively of Hodge type $(-1,1)$ (resp. $(1,-1)$) for the Hodge structure on $\fg$ defined by $h$.
 
We sometimes write $\fp^+ = \fp^+_h$, $\fk = \fk_h$, in order to emphasize the choice of $h$.   

The symmetric space $X$ has a $G(\RR)$-invariant metric, unique up to multiplication by a positive real scalar, that descends to
$S(G,X)$.   If $G/Z$ is anisotropic, then $S(G,X)$ is compact and its de Rham cohomology can be computed by means of harmonic forms with respect to the invariant metric.
Matsushima's formula reinterprets this computation in terms of relative Lie algebra cohomology of the space $\cA(G)$ of automorphic forms.

Let $\rho:  G \ra \Aut(V)$ be a finite-dimensional linear representation, defined over a subfield $E \subset \CC$.  Let $\tilde{V}$ be the corresponding local 
system over $S(G,X)$, with coefficients in $E$:
$$\tilde{V} =  \varprojlim_{K_f} G(\QQ)\backslash X \times G(\aA_f) \times V/K_f.$$
In order for $\tilde{V}$ to be a local coefficient system, we need to assume that the Zariski closure of  $Z(\QQ) \cap K_f$ in $Z$, for all sufficiently small $K_f$,
is in the kernel of $\rho$; this assumption will be made without comment in what follows.  We write $V_\CC = V\otimes_E \CC$ and define
$\tilde{V}_\CC$ in the same way.

By an {\it admissible irreducible representation} of $G(\aA)$ we always mean an irreducible $(\fg,K_\infty) \times G(\aA_f)$-representation, with the chosen $K_\infty$.

\begin{thm}[Matsushima's formula]\label{matsu}   Assume $G/Z$ is anisotropic.  Then there is a canonical isomorphism of $G(\aA_f)$-representations
$$H^* (S(G,X),\tilde{V}_\CC) \isom \oplus_{\pi}  m(\pi)H^* (\fg,\fk; \pi_\infty \otimes V_\CC)\otimes \pi_f.$$
Here $\pi$ runs through admissible irreducible representations of $G(\aA)$, and there is a countable direct sum decomposition
\begin{equation}\label{discreteAG} \cA(G) \isom \oplus_\pi m(\pi)\pi 
\end{equation}
where $m(\pi)$ is an (integer) multiplicity and $\pi \isom \pi_\infty \otimes \pi_f$ is the factorization with $\pi_\infty$ (resp. $\pi_f$) an irreducible  $(\fg,K_\infty)$-module (resp. $G(\aA_f)$-representation).
The notation $H^* (\fg,\fk; -)$ denotes relative Lie algebra cohomology.
\end{thm}

The standard reference for this result is the book \cite{BW} of Borel and Wallach, to which we refer for definitions and proofs.  The  interest of Matsushima's formula
is that it separates the calculation of the cohomology into a global part, corresponding to the multiplicities $m(\pi)$, and a  local part that depends only on $\pi_\infty$.  
We will only consider $\pi$ for which $\pi_\infty$ is tempered; in this case the cohomology $H^* (\fg,\fk; \pi_\infty \otimes V_\CC)$ is completely calculated
in \cite{BW}, and we copy the answer in the following section.

The result is more complicated when $G/Z$ is not anisotropic, and the proofs are considerably more difficult.  In particular, there is no longer a direct sum decomposition
as in \eqref{discreteAG}.  Nevertheless, it was proved by Franke \cite{franke98}  (building on earlier work of Borel) that the following analogue of Theorem \ref{matsu} holds:
\begin{equation}\label{franke}
H^* (S(G,X),\tilde{V}_\CC) \isom H^* (\fg,\fk; \cA(G) \otimes V_\CC).
\end{equation}

\subsubsection{Cuspidal cohomology}

We define the {\it cuspidal cohomology}  
$$H^* _0(S(G,X),\tilde{V}_\CC) \subset H^* (S(G,X),\tilde{V}_\CC)$$ to be the image
of $H^* (\fg,\fk; \cA_0(G) \otimes V_\CC)$ in $H^* (\fg,\fk; \cA(G) \otimes V_\CC)$, with respect to the isomorphism \eqref{franke}.  Borel proved
that the map from the relative Lie algebra cohomology of cusp forms to cohomology is injective.  Thus if we write
\begin{equation}\label{cuspAG} \cA_0(G) \isom \oplus_\pi m_0(\pi)\pi,
\end{equation}
where $m_0(\pi)$ denotes the multiplicity of $\pi$ in the cuspidal spectrum, we have the analogue of Matsushima's formula:
\begin{equation}\label{matsucusp}
H^* _0(S(G,X),\tilde{V}_\CC) \isom \oplus_{\pi}  m_0(\pi)H^* (\fg,\fk; \pi_\infty \otimes V_\CC)\otimes \pi_f.
\end{equation}

\begin{remark}  In \eqref{matsucusp} and in Theorem \ref{matsu}, the relative Lie algebra cohomology should really be replaced by $(\fg,K^c_\infty)$-cohomology,
where $K^c_\infty$ is now a maximal compact subgroup.  The version stated here is only correct when $K^c_\infty$ is connected and the center of $G$ is finite.  
The discussion here can easily be adapted to handle complications introduced by disconnectedness of $K^c_\infty$ or by the center of $G$.
\end{remark}

 We fix a representation $\pi_f$ of $G(\aA_f)$ and let $[\pi_f]_\infty$ denote the set of irreducible representations
 of $G(\RR)$ such that $\pi = \pi_\infty\otimes \pi_f \subset \cA_0(G)$.  We make the hypothesis
\begin{hyp}\label{temp}  Every $\pi_\infty \in [\pi_f]_\infty$ is tempered.
\end{hyp}
We can define
$$H^* _0(S(G,X),\tilde{V}_\CC)[\pi_f] = \Hom_{G(\aA_f)}(\pi_f, H^* _0(S(G,X),\tilde{V}_\CC)).$$
By \eqref{matsucusp} we then have 
\begin{equation}\label{arch1}
H^* _0(S(G,X),\tilde{V}_\CC)[\pi_f] \isoarrow \oplus_{\pi_\infty \in [\pi_f]_\infty} m_0(\pi) H^* (\fg,\fk; \pi_\infty \otimes V_\CC).
\end{equation}

\subsubsection{Relative Lie algebra cohomology in the tempered case}

The set of tempered $\pi_\infty$ with non-trivial relative Lie algebra cohomology with coefficients in a given $V_\CC$ is determined in \cite[\S 3]{BW},
especially in Theorems 3.3 and Theorem 5.1.  Every such $\pi_\infty$ is isomorphic to a parabolically induced representation of the form
$I_{P,\sigma,0}$, in the notation of \cite{BW}, where $P$ is a fundamental parabolic subgroup of $G(\RR)$ with Langlands decomposition $^0MAN$ and
$\sigma$ is a discrete series representation of $^0M$, whose infinitesimal character is determined by the highest weight of $V$:
$$\chi_\sigma = -\chi_{-s(\rho + \lambda(V))|_{\mathfrak{b}}}$$
(see \cite[Theorem 3.3 (2)]{BW}) where $\lambda(V)$ is the highest weight of $V$, $\mathfrak{b}$ is a compact Cartan subalgebra of $\Lie(^0M)$,
and $s$ is a (uniquely determined) element of the Weyl group of length $\ell(s) = \dim{N}/2$.     Moreover, the relative Lie algebra cohomology 
$H^* (\fg,\fk; \pi_\infty \otimes V_\CC)$
is determined explicitly in the two theorems quoted.  We let 
$$q(G) = \frac{\dim(G) - \dim(K^c_\infty)}{2}; ~ \ell_0(G) =  {\rm rank}(G) - {\rm rank}(K^c_\infty)$$
$$q_0(G) = q(G) - \frac{\ell_0(G)}{2}.$$
The invariant $\ell_0$ coincides with $\dim A$ where $A$ is the split part of the maximally split torus in $G(\R)$. Then $q_0(G)$ is an integer (cf. \cite[\S 4.3]{BW}) and we have 
\begin{equation}\label{q0}
H^{q_0(G) + j}(\fg,K^c_\infty; \pi_\infty \otimes V_\CC) \isoarrow  H^{q(^0M)}(\mathfrak{^0m},K_P;\sigma\otimes W_{s(\lambda + \rho) - \rho})\otimes \wedge^j(\mathfrak{a}^*)
\end{equation}
Here $\mathfrak{a}^*$ is the linear dual of the Lie algebra $\mathfrak{a}$ 
of the split component $A$, $W_{s(\lambda + \rho) - \rho}$ is the irreducible finite-dimensional representation of $^0M$ with the
indicated highest weight, and $K_P$ is the maximal compact subgroup $K^c_\infty\cap P(\RR) = K^c_\infty \cap ^0M(\RR)$.   Moreover, $^0M$ is a group with discrete series
and $q(^0M)$ is the dimension of its associated symmetric space.  

In particular, since $\ell_0(G) = \dim A$, we have
\begin{equation}\label{range}
H^q(\fg,K^c_\infty; \pi_\infty \otimes V_\CC) \neq 0 \Leftrightarrow q \in [q_0(G),q_0(G) + \ell_0(G)]
\end{equation}

\subsubsection{The exterior algebra action}\label{extaction}  A basic fact about discrete series, due to Schmid, is that the space
$ H^{q(^0M)}(\mathfrak{^0m},K_P;\sigma\otimes W_{s(\lambda + \rho) - \rho})$ that appears in
\eqref{q0} is {\it one-dimensional}.  It follows that
\begin{cor}\label{extactarch}  $H^{q_0(G) + * }(\fg,K^c_\infty; \pi_\infty \otimes V_\CC)$ is naturally a free rank one differential graded
module over the exterior algebra $\wedge^{\bullet}(\mathfrak{a}^*).$
\end{cor}

Combining this with \eqref{arch1} we obtain the following global corollary, which is the starting point of Venkatesh's
motivic conjectures.
\begin{cor}\label{extactglobal}  Suppose $\pi_f$ has the property that every $\pi_\infty \in [\pi_f]_\infty$ is tempered.
Then 
$$H^{q_0(G)+\bullet}_0(S(G,X),\tilde{V}_\CC)[\pi_f]$$
 is a free differential graded module over $\wedge^{\bullet}(\mathfrak{a}^*)$ of rank 
$$m_0(\pi_f) = \sum_{\pi_\infty \in [\pi_f]_\infty} m_0(\pi).$$
\end{cor}

The starting point of Venkatesh's motivic conjectures is the observation that the dimension $\ell_0$ of $\mathfrak{a}^*$ coincides with a
rank of a motivic cohomology group hypothetically attached to $\pi_f$.  This is explained clearly in \cite{PV,V19}, and we recall the explanation
here.   To a pair consisting of cohomological cuspidal automorphic representation  $\pi$ of $G$ and a finite-dimensional representation $\tau$ of the
Langlands $L$-group ${}^LG$, the Langlands correspondence hypothetically attaches a collection of motives $M(\tau)_\pi$, obtained by composing $\tau$ with
the parameter of $\pi$ viewed as a homomorphism of a motivic version of the absolute Galois group of $\QQ$ to the algebraic group ${}^LG$.  For the purposes
of this and the following sections, we only need the case of the dual adjoint representation $\tau = \Ad^*$; then we write $M_\pi = M(Ad)_\pi$.  In general the motives have coefficients in a number field $E_\pi$, which we ignore temporarily.  

Attached to the hypothetical motive is the group $H^1_{\rm mot}(\QQ,M_\pi(1))$ of motivic cohomology.  This is supposed to be a direct summand of a certain
$K$-theory group tensored with $\QQ$, and is conjecturally a finite-dimensional $\QQ$-vector space (more generally an $E_\pi$-vector space);  
we denote this space  $\text{\dn{v}}_\pi$ (pronounced ``Va," in honor of Venkatesh).   While
both $M_\pi$ and its motivic cohomology are inaccessible, if $M_\pi$ is an actual piece of the cohomology of a smooth projective algebraic variety over $\QQ$
cut out by correspondences, then there are regulator maps comparing $\text{\dn{v}}_\pi\otimes C$ and certain  spaces of cohomological realizations, where
$C$ runs through possible fields of coefficients.  The key point
is\footnote{Normally one defines $\dnv_\pi$ as the space of classes that extend to an appropriate integral model.  In the applications Venkatesh explains how this condition can be ignored.  See for example \cite[p. 4]{PV}.}
\begin{conj}[Beilinson, Bloch-Kato]\label{conj-bbk}  (i)  The space $\dnv_\pi$ is of dimension $\ell_0$.

(ii)  The Beilinson regulator defines an isomorphism
$$\dnv_\pi \otimes_{\QQ} \CC \isoarrow H^1_\cD(M_{\pi,HdR},1),$$
where $H^1_\cD(*,i)$ is Deligne cohomology and $M_{\pi,HdR}$ denotes the Hodge-de Rham realization.

(iii)  For any prime $p$ at which $M_\pi$ has good reduction, the $p$-adic regulator defines an isomorphism
$$\dnv_\pi \otimes_{\QQ} \QQ_p \isoarrow H^1_f(\QQ, M_{\pi}(1)),$$
where $H^1_f$ is the Bloch-Kato Selmer group.
\end{conj}

Each conjecture also asserts that the space on the right-hand side in the isomorphism has dimension $\ell_0$.  This can often be proved or deduced from other conjectures.  Most importantly for our purposes, these conjectures provide rational structures on the exterior algebra $\wedge^{\bullet}(\mathfrak{a}^*)$ that can be compared by means of \eqref{arch1} to the rational structure on the cohomology space 
$H^*_0(S(G,X),\tilde{V}_\CC)[\pi_f]$, induced
by topological cohomology.  Specifically, we have

\begin{fact}[\cite{PV}, (5.1.1)]\label{regLie}  Conjecture \ref{conj-bbk} (i) and (ii) give rise to a canonical isomorphism
$$\dnv_\pi \otimes_{\QQ} \CC \isoarrow  \mathfrak{a}^*.$$
Thus Conjecture \ref{conj-bbk} endows $\mathfrak{a}^*$ with a canonical $\QQ$-rational structure.
\end{fact}

The conjectures of \cite{PV,V19,GV} all refer back to this $\QQ$-rational structure, and in particular to the following conjecture.
Since $\pi_f$ is realized in cohomology, it has a model over a number field $E(\pi_f).$\footnote{If $G = GL(n)$ it follows from the multiplicity one theorem that this can identified with the field of rationality of $\pi_f$, that is the fixed field in $Gal(\ol{\Q}/\Q)$
of the subgroup that fixes the isomorphism class of $\pi_f$.  In general $E(\pi_f)$ may be slightly larger.}  The $\pi_f$ isotypic component of $H^*_0(S(G,X),\tilde{V}_\CC)$ is then a direct summand of the $E(\pi_f)$-rational cohomology
$H^*_0(S(G,X),\tilde{V}_{E(\pi_f)})$.  We state a slight generalization of the main conjecture of Prasanna and Venkatesh:
\begin{conj}[\cite{PV}, Conjecture 1.2.1]\label{rational}   The action on $H^*_0(S(G,X),\tilde{V}_\CC)[\pi_f]$ of $\wedge^{\bu} \dnv_\pi$ preserves the $E(\pi_f)$-rational structure.
\end{conj}

We describe the  actions of the $p$-adic realizations of $\dnv_\pi$ in greater detail in \S \ref{venk}.

\subsection{Coherent cohomology of Shimura varieties}

We return to the notation of \S \ref{derham} but now we suppose that $X$ has $G(\RR)$-invariant hermitian structure.  We assume throughout that
$Z_G$ is trivial; this is certainly unnecessary but it simplifies the exposition.  We define 
$$\frP = \frP_h = \fk \oplus \fp^-,$$
in the notation of  \eqref{CartanHC}.  This is a maximal parabolic subalgebra with unipotent radical $\fp^-$.  The locally symmetric space $S(G,X)$ is then
a Shimura variety, and has a canonical model over a number field $E(G,X)$ that is preserved under the right action of $G(\aA_f)$.  This induces a canonical model
on the individual finite-level quotients $S{}_{K_f}(G,X)$.

Let $\tau: K_\infty \ra \GL(W_\tau)$ be an irreducible representation, and define the holomorphic vector bundle
\begin{equation}\label{autovb}
\cE_\tau = \varprojlim_{K_f} G(\QQ)\backslash G(\aA) \times W_\tau/K_\infty \times K_f
\end{equation}
where now $G(\QQ)$ acts on $G(\aA)$ on the left, $K_f$ acts on $G(\aA_f)$ on the right, and $K_\infty$
acts on $G(\RR)$ on the right and on $W_\tau$ on the left.   If $K_f$ is sufficiently small the bundle $\cE_\tau$ descends to a vector bundle, also denoted $\cE_\tau$,
on $S{}_{K_f}(G,X)$, which has a canonical model over a finite extension $E(G,X,\tau)$ of $E(G,X)$.  
In practice we fix a level subgroup $K_f$ such that $S{}_{K_f}(G,X)$ has a {\it canonical integral model}, as in the work of Kisin, Shin, and Zhu \cite{KSZ}; this 
presupposes that $S(G,X)$ is a Shimura variety of {\it abelian type}, which is closely related to moduli of abelian varieties with additional structure.

Suppose for the moment that $G$ is anisotropic over $\QQ$, so that $S(G,X)$ is (pro)-projective.  The  cohomology of
$S(G,X)$ with coefficients in (the coherent sheaf attached to) the vector bundle $\cE_\tau$ then admits an expression in terms of differential forms (Dolbeault cohomology),
which in turn can be written in terms of relative Lie algebra cohomology, analogously to Matsushima's formula:
\begin{equation}\label{dolbeault}
H^*(S(G,X),\cE_\tau) \isom H^*(\frP,\fk; \cA(G) \otimes W_\tau) \isom \oplus_{\pi}  m(\pi)H^*(\frP,\fk; \pi_\infty \otimes W_\tau)\otimes \pi_f.
\end{equation}
(Again, one should replace $(\frP,\fk)$-cohomology by $(\frP,K_h)$-cohomology, where $K_h \supset K_\infty$ is the stabilizer of $h$ in $G(\RR)$ and is not necessarily
connected.)
When $S(G,X)$ is not projective, one needs to replace the left-hand side of \eqref{dolbeault} by $H^*(S(G,X)_\Sigma,\cE_\tau^{can})$,
where $S(G,X)_\Sigma$ is a well-chosen (smooth projective) toroidal compactification and $\cE_\tau^{can}$ is Mumford's (or Deligne's) canonical
extension of $\cE_\tau$ to a vector bundle over $S(G,X)_\Sigma$; the result is independent of the choice of $\Sigma$.  
The analogue of the cuspidal cohomology \eqref{matsucusp} 
\begin{equation}\label{dolbeaultcusp}
H^*_0(S(G,X),\cE_\tau) \isom \oplus_{\pi \subset \cA_0(G)}  m(\pi)H^*(\frP,\fk; \pi_\infty \otimes W_\tau)\otimes \pi_f.
\end{equation}
was studied in \cite{H88}, where in most cases the left hand side is defined algebraically in terms of two distinct extensions of $\cE_\tau$ to $S(G,X)_\Sigma$.  The analogue of Franke's theorem \eqref{franke} holds for coherent cohomology and was proved recently in the thesis of Jun Su \cite{junsu}.  

Suppose the infinitesimal character of $\tau$, as a representation of the enveloping algebra of $\fk$, restricts on the enveloping algebra of $\fg$ to
the infinitesimal character of a finite-dimensional irreducible representation.  Then Schmid's theorem asserts that there is a unique $\pi_\infty$ and a unique $q(\tau)
\in \{0,1, \dots, \dim X\}$ such that
$H^{q(\tau)}(\frP,\fk; \pi_\infty \otimes W_\tau) \neq 0$; moreover, the cohomology space is then of dimension $1$, and $\pi_\infty$ is a discrete series.  This theory is described
at length in \cite{H88}.   In particular, the global coherent cohomology $H^q_0(S(G,X),\cE_\tau)$ is non-trivial only when $q = q(\tau)$.  

The version of Venkatesh's motivic conjectures for coherent cohomology, to be discussed in \S \ref{coherent}, concerns $\tau$ for which
$H^*_0(S(G,X),\cE_\tau)$ has cohomology in several degrees, and the corresponding $\pi_\infty$ belong to the {\it nondegenerate limit of discrete series} (NLDS).   It is shown in \cite{BHR}, following work of Floyd Williams, that any given NLDS $\pi_\infty$ contributes to $H^q(\frP,\fk; \pi_\infty \otimes W_\tau)$ for a unique $q$.  Let $P(\tau)$ denote
the set of NLDS $\sigma$ such that $H^q(\frP,\fk; \sigma \otimes W_\tau) \neq 0$, and define 
$$\Pi(\tau) = \oplus_{\sigma \in P(\tau)} \sigma.$$
For the purposes of this definition, it is most convenient to treat the different components of the restriction of any NLDS of $G(\RR)$ to its identity component as separate representations $\sigma$.

\begin{question}  Is there an abelian Lie algebra $\fra \subset \frp$ such that 
$$H^*(\frP,\fk; \Pi(\tau) \otimes W_\tau)$$
is a free module over $\wedge^\bullet(\fra^*)$?
\end{question}

The answer is affirmative in the admittedly unrepresentative case where $G = \GL_2$, but the question does not seem to have been explored more generally.

\subsection{Action of $q$-adic derived Hecke algebras on cohomology}\label{DHAcoh}

We return to the topological setting of \S \ref{derham}, and specifically the Prasanna-Venkatesh Conjecture \ref{regLie}.  Thus let
$\pi_f$ be as in that section; in particular, we assume all $\pi_\infty \in [\pi_f]_\infty$ are tempered.  We assume $E(\pi_f) = \QQ$;  this hypothesis can easily be relaxed, at the cost of complicating the notation.  We similarly assume throughout that the representation $V$ is defined over $\QQ$, and write $\tilde{V}$ for the corresponding local system with coefficients $\QQ$.   

Let $p$ be a prime number such that the local component  $\pi_p$ is unramified.  We fix a level subgroup $K_0 = \prod_v K_v \subset G(\aA_f)$, with $K_p$ hyperspecial maximal, such that $\pi_f^{K_0} \neq 0$.  We will be working with $\QQ_p$-valued cohomology of the locally symmetric space $S_{K_0}(G,X)$, with coefficients in $\tilde{V}\otimes \QQ_p$.   More precisely,
we choose a $\ZZ_p$-local system $\Lambda_V \subset \tilde{V}$, and define the Hecke algebra $\TT_{K_0}$, as in \cite[\S 6.1]{V19}, to be the $\ZZ_p$-algebra generated by Hecke operators at good primes, acting as endomorphisms of the chain complex of $\Lambda_V$, viewed as an object in the derived category.   
Then $\pi_f$ determines a maximal ideal $\fm = \fm_\pi \subset \TT_{K_0}$,
which we reserve for later.  

Write $S_{K_0} = S_{K_0}(G,X)$.  For any $r = 1, 2, \dots$ we let 
 $\Lambda_{V,r} = \Lambda_V/p^r\Lambda_V$, $A_r = \ZZ/p^r\ZZ$.
We consider the object 
$$\rR\Gamma(S_{K_0},\Lambda_{V,r})$$
as a module over $\TT_{K_0}$.
Let $q$ be a prime such that $K_q$ is hyperspecial maximal compact.  Consider the family $\cF_q$ of open compact subgroups $U \subset K_0$
that contain $\prod_{v \neq q} K_v$, and for each $U \in \cF_q$, the chain complex $C^\bullet(S_U(G,X),\Lambda_{V,r})$ as an object
in the derived category of $A_r$-modules.  Then $G_q = G(\QQ_q)$ acts on 
$C_q^\bullet(A_r) = \varinjlim C^\bullet(S_U(G,X),\Lambda_{V,r})$.  Let $\mf{U}_q$ be the universal $A_r[G_q]$-module
$\mf{U}_{A_r}(K_q)$.  Then just as in \S \ref{actiondHa},  Frobenius reciprocity identifies
\begin{equation}
\RHom_{A_r[K_q]}(A_r,C_q^\bullet(A_r)) \isoarrow \RHom_{A_r[G_q]}(\mf{U}_{A_r}(K_q),C_q^\bullet(A_r))
\end{equation}
and thus by Fact \ref{DHAaction} we have an action of the derived Hecke algebra
$\cH^*_{A_r}(G_q,K_q)$ on $H^*(K_q,C_q^\bullet(A_r))$.  

By the discussion in \cite[\S 7.5]{HV}, we can identify 
the hypercohomology $H^*(K_q,C_q^\bullet(A_r))$ with $H^*(S_{K_0},\Lambda_{V,r})$.   Thus formally we have 
\begin{prop}\label{DHactionq}  For any $r$ and any $q$, there is a canonical action of the derived Hecke algebra $\cH(G_q,K_q; {A_r})$
on $H^*(S_{K_0},\Lambda_{V,r})$.
\end{prop}

The action in Proposition \ref{DHactionq} is only of interest when $\cH_{A_r}(G_q,K_q)$ is non-trivial.  Thus we assume $q \equiv 1 \pmod{p^r}$.  For such $q$ we have seen in Theorem \ref{thm: derived sat isom} that
\begin{equation}\label{d-Satake}
 \cH^*(G_q,K_q; {A_r}) \xrightarrow{\sim}  \cH^*(T(\Q_q), T(\Z_q); {A_r})^W
 \xrightarrow{\sim} \left( A_r[X^*(T)] \otimes_{A_r} H^*(T_q; A_r) \right)
 ^W
\end{equation}
where $W$ denotes the Weyl group of $G$ and we have adapted the notation for the algebra on the right hand side.

\begin{remark}  When $q \equiv 1 \pmod{p^r}$,
the action of an individual derived Hecke operator in $\cH^1_{A_r}(G_q,K_q)$ on cohomology modulo $p^r$
can be made explicit in terms of the mod $p^r$ characteristic classes of cyclic unramified topological coverings of degree $p^r$
in the passage from Iwahori level  to  pro-$q$ Iwahori level at $q$.
They are thus seen to be derived versions of the diamond operators that play a central role in the Taylor-Wiles method.
This is explained in greater detail in \S \ref{coherent} in the setting of coherent cohomology of modular curves.  The topological version was developed in \cite{V19} and is completely analogous.
\end{remark}

Venkatesh explains in \cite[\S 2.13]{V19} how to take an increasing union of sets of $q$ as above, while passing to the limit on $r$, in order to define various versions of a {\it global} $p$-adic derived Hecke algebra  that acts on  $H^*(S_{K_0},\Lambda_{V})$.  We choose the {\it strict global derived Hecke algebra}, denoted $\tilde{\TT}'$ in \cite{V19}, and assume $p$ is prime to the order of the Weyl group.  Then $\tilde{\TT}'$ is graded commutative and contains  underived Hecke operators at good primes other than $p$.  These generate the usual (underived) Hecke algebra $\TT_{K_0} \subset \tilde{\TT}'$.  Moreover, the action of $\tilde{\TT}'$ preserves
the localization $H^*(S_{K_0},\Lambda_{V})_\fm = H^*(S_{K_0},\Lambda_{V})_{\fm_\pi}$ at the maximal ideal of $\TT_{K_0}$ we have been saving for this moment.

The first main result of \cite{V19} is
\begin{thm}\cite[Theorem 7.6]{V19}\label{Hecke-cyclic}  Under the assumptions listed below, $H^*(S_{K_0},\Lambda_{V})_\fm$ is generated as a graded module over $\tilde{\TT}'$ by its minimal degree component $H^{q_0(G)}(S_{K_0},\Lambda_{V})_\fm$.
\end{thm}

The proof is an adaptation of the ideas of Calegari-Geraghty, in the version developed by Khare and Thorne,
to the language of derived structures.  It depends on a list of technical assumptions \cite[\S 6.1]{V19} and one substantial assumption \cite[\S 6.2]{V19}.  The most significant of the technical assumptions is that the $\fm$-localized homology [N.B.] $H_*(S_{K_0},\Lambda_{V})_\fm$ is concentrated in the interval
$[q_0(G),q_0(G)+\ell_0(G)]$ that we have already seen in the calculation of relative Lie algebra cohomology of tempered representations.  Related to this is the assumption that $H_*(S_{K_0},\Lambda_{V})$ is torsion-free.  Finally, $G$ is assumed to be a split $\QQ$-group; this appears to rule out generalizations to simple groups over totally real fields, for example, though the analogous theorem must hold for such groups as well.

The more serious assumption is that $\pi$ has an associated $p$-adic Galois parameter.  This is a continuous homomorphism
\begin{equation}\label{rhopi}
\rho_\pi:  \Gal(\ol{\QQ}/\QQ) \ra \chG(\ol{\QQ_p}),
\end{equation}
where $\chG$ is the Langlands dual group to $G$, whose Frobenius conjugacy class at the unramified prime $q$ corresponds to the Satake parameter of $\pi_q$, and whose restriction to a decomposition group at $p$ is crystalline in an appropriate sense.  The existence of such a $\rho_\pi$ is known when $G$ is the general linear group over
a real or CM field, under our assumption that $\pi$ contributes to cohomology of some $S_{K_0}$.  
Venkatesh then makes assumptions on the image of the reduction $\ol{\rho}_\pi$ mod $p$ that are familiar from the Taylor-Wiles theory, in particular that the image of
$\ol{\rho}_\pi$ is not contained in a proper parabolic subgroup of $\chG$ and is not too small.

\begin{remark} Strictly speaking, Venkatesh only proves  Theorem \ref{Hecke-cyclic} when $V$ is the trivial representation.  Moreover, he assumes that $\pi$ is congruent to no other representations at level $K_0$, and indeed that
$\TT_{K_0,\fm} \isoarrow \Z_p$.  As far as we know, no one has written down a proof without these hypotheses, which presumably would require additional restrictions on $p$.
\end{remark}

The second main result is the construction of an action of the $p$-adic realization of motivic cohomology, as in Conjecture \ref{conj-bbk} (iii), on $H^*(S_{K_0},\Lambda_{V})_\fm$ as endomorphisms increasing the degree by $1$.  This will be discussed in the following section, which is devoted to Venkatesh's main conjectures regarding this action and to generalizations in other settings.

\section{Venkatesh's motivic conjectures}\label{venk}

\subsection{The action of motivic cohomology}\label{actmot}

We retain the assumptions of \S \ref{DHAcoh}.  Then Venkatesh's Theorem \ref{Hecke-cyclic} tells us that the 
$\ell$-adic cohomology of the locally symmetric space, localized at the maximal ideal corresponding to $\pi$,
is cyclic over the strict global derived Hecke algebra (DHA) .  On the other hand, the analytic theory tells us that, after tensoring
with $\QQ$, it is cyclic over the exterior algebra on a vector space of dimension $\ell_0(G)$.  This invariant  plays no role in the definition of the strict global DHA but  is linked to motivic cohomology by the Bloch-Kato Conjecture.   The relation
between the two appearances of $\ell_0(G)$ is provided by a canonical map from the dual adjoint Selmer group on the 
right hand side of  Conjecture \ref{conj-bbk} (iii) to the degree 1 part of the strict global DHA that is characterized by a local reciprocity law.  

The local reciprocity law is nevertheless stated in a global context.  We realize $\pi$ and the Galois
representation in \eqref{rhopi} over the ring of integers $\cO$ in
a $p$-adic field, with residue field $k$.
Choose a Taylor-Wiles system for $\pi$, a set $Q_r = (q_1,\dots,q_m)$ of primes $q_i \equiv 1 \pmod{p^r}$, as in \cite[\S 6.3]{V19}, to which we refer for the relevant properties.\footnote{See also the chapters of Caraiani-Shin  and Emerton-Gee-Hellmann in these proceedings.  In practice we can work over a base field other than $\QQ$, of course, and the $q_i$ will then be prime ideals.}   We write $\cH^*_{q,A_r}$ for
$ \cH^*(G_q,K_q; {A_r})$.
In particular $\pi$ is unramified at each $q \in Q_r$ and the Satake parameters of each $\pi_{q}$, $q \in Q_r$, are strongly regular. This means that if $\fm_q \subset \rH(T(\Q_q),T(\Z_q); {A_r})^W$ is the maximal ideal induced by the maximal ideal $\fm_\pi$, then the fiber over $\fm_q$ in
$$\Spec(\rH(T(\Q_q),T(\Z_q), A_r) = A_r[X^*(T^\vee)] = A_r[X_*(T)])$$ 
is a set of $|W|$ maximal ideals.  We choose one maximal ideal  in $\Spec(\cH^0(T(\Q_q),T(\Z_q);{A_r}))$ for each $q$; this corresponds to choosing an element of the Langlands dual torus $T^{\vee}(k)$ that is contained in the  conjugacy class of $\ol{\rho}_\pi(\Frob_q)$; here $\ol{\rho}$ denotes the reduction of
$\rho_\pi$ in $G^\vee(k)$.   Since $W$ acts on both factors in the tensor product in \eqref{d-Satake},
this choice together with the derived Satake isomorphism allows us to define a canonical injection \cite[\S 8.17]{V19}
\begin{equation}\label{mapi}
\iota_{q,r}:  H^1(T_q; A_r) \hookrightarrow  \cH^1_{q,A_r}
\end{equation}

We can define a $\cO$-lattice in the
dual adjoint Selmer group $H^1_f(\QQ, M_{\pi}(1))$.  With Conjecture \ref{conj-bbk} in mind, we let $\dnv_{\pi,p}^*$ denote the $\cO$-dual of this lattice.   
Under the running assumptions, Venkatesh proves in \cite[Lemma 8.9]{V19} that $\dnv_{\pi,p}$ is free of rank $\ell_0(G)$
over $\cO$.  Note that global duality identifies
$\dnv_{\pi,p}^*$ with a lattice in (the space of obstructions) $H^2_f(\QQ,Ad(\rho_\pi))$.

There is a second map 
\begin{equation}\label{mapf}
f_{q,r}:  H^1(T_q,A_r) \rightarrow \dnv^*_{\pi,p}/p^n
\end{equation}
defined by local duality of Galois cohomology (see p. 78 of \cite{V19}).   

Venkatesh's second theorem, after Theorem \ref{Hecke-cyclic}, is then a more precise version of the following:

\begin{thm}\cite[Theorem 8.5]{V19}\label{motaction}  There is an action of $\dnv^*_{\pi,p}$ on $H^*(S_{K_0},\Lambda_{V})_\fm$ by endomorphisms of degree $+1$.  This action is uniquely characterized by the property that, for any $r \geq 1$ and any Taylor-Wiles prime $q \in Q_{a(r)}$ for sufficiently large $a(r)$, the  two actions of $H^1(T_q,A_r)$ on $H^*(S_{K_0},\Lambda_{V})_\fm$ -- by derived Hecke operators (via $\iota_{q,r}$) or by
the action of  $\dnv^*_{\pi,p}$ (via $f_{q,r}$) coincide.

In particular, there is an embedding of the exterior algebra $\wedge^\bullet \dnv^*_{\pi,p}$ in the strict global DHA $\TT'$ of
endomorphisms of $H^*(S_{K_0},\Lambda_{V})_\fm$.
\end{thm}

\begin{remark}[The spectral Hecke algebra]\label{rem: SHA}
The proof of Theorem \ref{motaction} is an intense computation in \cite{V19}. Another perspective is given by \cite{F20}. There, the author constructs an object called the \emph{spectral Hecke algebra}; it is a derived Hecke algebra that acts on the spectral (i.e., Galois) side of the Langlands correspondence. The construction is motivated by Geometric Langlands, but we will not explain that motivation here; we refer to \cite[Introduction]{F20} for details.  

The usual spherical Hecke algebra for $G$ is the space of functions on $G(\Z_q) \bs G(\Q_q) / G(\Z_q)$ viewed as a double coset set. It is natural to view this double coset space as a \emph{groupoid}, i.e. remembering automorphisms. Then it has a non-trivial topological structure, and taking cohomology (instead of just functions) gives the derived Hecke algebra. For $G = \GL_n$, the groupoid may be interpreted as the space of tuples $(\cF_1, \cF_2, \tau)$ where $\cF_1, \cF_2$ are rank $n$ vector bundles on $\Spec \Z_q$ and $\tau$ is an identification of their restrictions to $\Spec \Q_q$. 

To make the spectral Hecke algebra, we instead consider the space of tuples $(F_1, F_2, \tau)$ where $F_1, F_2$ are rank $n$ local systems on $\Spec \Z_q$ and $\tau$ is an identification of their restrictions to $\Spec \Q_q$. Although the language of the two situations is parallel, the geometry is quite different. In the second case, the construction is only interesting if performed in a derived sense, because the space of rank $n$ local systems on $\Spec \Z_q$ embeds in the space of rank $n$ local systems on $\Spec \Q_q$ (so that the information of $F_2$ is redundant). We call this derived space the ``spectral Hecke stack'', and its functions are essentially the spectral Hecke algebra of \cite{F20}. 

There is a co-action of the spectral Hecke algebra on the derived Galois deformation ring, and \cite{F20} computes that it is dual to Venkatesh's reciprocity law for the action of the derived Hecke algebra; furthermore \cite{F20} explains that there is a co-multiplication on the spectral Hecke algebra making it dual to the derived Hecke algebra. Then \cite{F20} gives a conceptual interpretation of Venkatesh's reciprocity law as a ``derived local-global compatibility'' between the actions of the derived Hecke algebra and the spectral Hecke algebra. 

\end{remark}

Now we tensor the motivic cohomology and the cohomology of $S_{K_0}$ with $\QQ_p$ and assume the Beilinson-Bloch-Kato Conjecture \ref{conj-bbk} (iii).   We let $\dnv^*_{\pi,E(\pi_f)} \subset \dnv^*_{\pi,p}\otimes \QQ_p$ be the space of classes whose pairing with motivic cohomology lies in $E(\pi_f)$.  Let $H^*(S(G,X),\tilde{V})_\pi$ denote the $\pi_f$-isotypic subspace
of $H^*(S(G,X),\tilde{V})$. 
The $p$-adic version of Venkatesh's motivic conjecture is

\begin{conj}\cite[Conjecture 8.8]{V19}\label{motconj}  Under the action of $\wedge^\bullet \dnv^*_{\pi,p}$ on
$$[H^*(S(G,X),\tilde{V})_\pi]^{K_0} \otimes \QQ_p \subset H^*(S_{K_0},\Lambda_{V})_\fm\otimes \QQ_p$$
defined by Theorem \ref{motaction}, the action of $\dnv^*_{\pi,E(\pi_f)}$
preserves the $E(\pi_f)$ rational structure induced from $H^*(S(G,X),\Lambda_{V})_\pi$.
\end{conj}

The strongest evidence for this conjecture, apart from the results of \cite{DHRV} on a similar conjecture
for coherent cohomology of modular curves (see \S \ref{coherent}), is the substantial evidence provided in \cite{PV} for the analogous conjecture for complex cohomology.  In that conjecture the $p$-adic regulator with values in the Bloch-Kato Selmer group is replaced by the 
Beilinson regulator with values in Deligne cohomology, as in Conjecture \ref{conj-bbk} (ii).  Since derived structures are
not required for the statement of the Prasanna-Venkatesh conjecture, we say no more about it here.

\subsection{The action of the derived deformation ring}\label{actddr}After constructing the derived deformation ring in \cite{GV}, Galatius and Venkatesh use it to give an interpretation in derived geometry of the method of Calegari-Geraghty.  
The cuspidal cohomological representation $\pi$, the localized cohomology $H^*(S_{K_0},\Lambda_{V})_\fm$, the Galois parameter $\rho_\pi$, its reduction $\ol{\rho}_\pi$ modulo $p$, and the dual adjoint Selmer group $H^1_f(\QQ, M_{\pi}(1))$ are as in \S \ref{DHAcoh} and \S\ref{actmot}.   
Let $S$ be the set of places where $\pi$ and $K_0$ are ramified, and assume $p \notin S$.   
Then we can define the crystalline derived deformation ring $\cR^{\ol{\rho}, \crys}_{\cO_F[1/S]}$, with unrestricted ramification at $S$.  \footnote{The assumptions of \cite{V19} already admitted in the previous section remain in force in \cite{GV}; these include the hypothesis that the local unrestricted deformation spaces at primes in $S$ are well-behaved.}  

The Galatius-Venkatesh interpretation of the Calegari-Geraghty method is formulated in terms of {\it homology} rather than
{\it cohomology} of $S_{K_0}$.  To avoid complications in the duality between homology and cohomology
we assume $V$ is the trivial representation, so $\Lambda_V$ is just
the constant local system $\cO$; \cite{GV} does not treat the more general case. The graded homology group $H_*(S_{K_0}, \cO)_{\mf{m}}$ is $H_*(\cM_0)[q_0]$ in the notation of \S \ref{ssec: calegari-geraghty}. The proof of the following theorem, the main result of \cite{GV}, was sketched earlier: 

\begin{thm}\cite[Theorem 14.1]{GV}\label{GVtheorem}  Under the hypotheses of the previous sections the localized homology
$H_*(S_{K_0},\cO)_\fm$ carries the structure of a free graded module over the graded ring $\pi_*(\cR^{\ol{\rho}, \crys}_{\cO_F[1/S]})$. 
\end{thm}

\subsubsection{Duality of the derived Hecke and derived deformation actions}\label{duality}
The article \cite{GV} does not contain any new formulations of the motivic conjectures, unlike the article \cite{V19} devoted
to the derived Hecke action.  However, when the two actions can be defined they are related canonically.  

The main point is the determination of $\pi_1(\cR^{\ol{\rho}, \crys}_{\cO_F[1/S]})$ \cite[(15.7)]{GV}:  in our notation, we have

\begin{lemma}\cite[Lemma 15.3]{GV}\label{pi1rs}  There is an isomorphism
$$\pi_1(\cR^{\ol{\rho}, \crys}_{\cO_F[1/S]}) \isoarrow \dnv_{\pi,p}
$$
that induces an isomorphism of graded algebras
$$\pi_*(\cR^{\ol{\rho}, \crys}_{\cO_F[1/S]}) \isoarrow  \wedge^*\dnv_{\pi,p}
$$
\end{lemma}

The relation between the two actions cannot be completely straightforward, because the derived deformation ring acts on homology while the (strict global) derived Hecke algebra acts on cohomology.   
Under the duality between cohomology and homology, the action of  $\wedge^\bullet \dnv^*_{\pi,p}$
on $H^*(S_{K_0},\cO)_\fm$ defined in Theorem \ref{motaction} becomes an adjoint action on
$H_*(S_{K_0},\cO)_\fm$; $\dnv^*_{\pi,p}$ acts as operators of degree $-1$.  The compatibility between
the two actions takes the following form:
\begin{thm}\cite[Theorem 15.2]{GV}\label{GVcompatible}  The DHA action of $\wedge^\bullet \dnv^*_{\pi,p}$ and the derived deformation ring
action of $\pi_* (\cR^{\ol{\rho}, \crys}_{\cO_F[1/S]})$ are compatible in the following sense.  Under the identification
$$\pi_*(\cR^{\ol{\rho}, \crys}_{\cO_F[1/S]}) \isoarrow  \wedge^*\dnv_{\pi,p}
$$
of Lemma \ref{pi1rs}, let $v \in \dnv_{\pi,p} \in \pi_1 (\cR^{\ol{\rho}, \crys}_{\cO_F[1/S]})$, $v^* \in \dnv^*_{\pi,p}$, $m \in H_*(S_{K_0},\cO)_\fm$.  Then
$$v^*\cdot v\cdot m + v \cdot v^*\cdot m = \langle v,v^* \rangle m$$
where $\langle , \rangle$ is the natural $\cO$-valued pairing between $\dnv_{\pi,p}$ and $ \dnv^*_{\pi,p}$.
\end{thm}

\subsection{Coherent cohomology}\label{coherent}

\subsubsection{The derived Hecke action}
In this section we use the classical notation for modular curves.   If $d > 1$ is an integer we let $X_0(d)$ and
$X_1(d)$ be the familiar modular curves over $\ZZ[\frac{1}{d}]$ with level structure of type $\Gamma_0(d)$
and $\Gamma_1(d)$ respectively.  We think of $X_?(d)$, $? = 0, 1$, as the quotient
$$\Gamma_?(d)\backslash \mathfrak{H}; ~~ \mathfrak{H} = \GL_2(\RR)^0/\S\mathrm{O}_2\cdot \RR^\times,$$
where $\GL_2(\RR)^0$ is the identity component of $\GL_2(\RR)$.  
The line bundle over $X_?(d)$, $? = 0, 1$ whose sections are modular forms
of weight $1$ is denoted $\omega$; the isotropy representation of $\S\mathrm{O}_2$ on the fiber of $\omega$ at
its fixed point $h \in \mathfrak{H}$ is denoted $\omega_h$.   

The bundle $\omega$ is the only automorphic line bundle with non-trivial cuspidal
coherent cohomology, in the sense of \eqref{dolbeaultcusp}, in more than one degree -- in other words in degrees
$0$ and $1$.  Moreover, $H^0_0(X_?(d),\omega)$ and $H^1_0(X_?(d),\omega)$ are isomorphic 
as modules over the classical Hecke algebra, which we denote $\TT_{?d}$; this is the algebra generated by
standard Hecke operators at all primes not dividing $d$.  This remains true for  modular curves attached to all congruence
subgroups:  there are two representations (Harish-Chandra modules) of the identity component $\GL_2(\RR)^0$,
say $\pi_0$ and $\pi_1$, such that
$$\dim H^i(\frP,\fk;\pi_i\otimes \omega_h) = 1, ~ i = 0, 1$$
in the notation of \eqref{dolbeault}.  The representations $\pi_0$ and $\pi_1$ are  holomorphic and antiholomorphic limits of discrete series, respectively, and the direct sum $I(0) := \pi_0 \oplus \pi_1$ has a structure of irreducible 
$(\Lie(\GL_2), \mrm{O}_2(\RR)\cdot\RR^\times)$-module whose restriction to $(\Lie(\GL_2),\S\mathrm{O}_2\cdot\RR^\times)$ is the sum of the two Harish-Chandra modules for $\GL_2(\RR)^0$.  The notation $I(0)$ is meant to suggest (accurately) that, up to a harmless twist by a power of the absolute value of the determinant, it is obtained by parabolic induction from a unitary character of the Borel subgroup.

For any prime $q$ we let $G_q = \GL(2,\QQ_q)$, $K_q = \GL(2,\ZZ_q)$.  Suppose $q$ is relatively prime to $d$ and is congruent to $1$ $\pmod{p^r}$ for some $r > 0$.  We let $X_?(d)_r$ denote the base change of $X_?(d)$ to $\Spec(\ZZ/p^r\ZZ)$ and
define $\omega_r$ analogously. 
One defines  in \cite{HV} an action of the derived Hecke algebra $\cH(G_q,K_q; \Z/p^r \Z)$
on $H^*(X_?(d)_r,\omega_r)$.  The action can be defined formally  as in \S \ref{DHAcoh}, but it
also has a concrete definition in terms of cyclic coverings.  We let $X_{0?}(qd)$ denote the modular curve corresponding
to the congruence subgroup $\Gamma_?(d)\cap \Gamma_0(q)$, and define $X_{1?}(qd)$ similarly.   As long as $q \geq 5$ the natural map 
$X_{1?}(qd) \ra X_{0?}(qd)$
factors through an \'etale cover $X_\Delta(qd)/X_{0?}(qd)$ that is cyclic of degree $p^r$.  This cover defines a class
\begin{equation}\label{shim}  \frS \in H^1(X_{0?}(qd)_r; C_r)
\end{equation}
that is called the {\it Shimura class}.  Here $C_r$ is the cyclic quotient of $(\ZZ/q\ZZ)^*$ of degree $p^r$.  This can be identified
with $\ZZ/p^r\ZZ$ but not canonically; we choose an isomorphism (discrete logarithm)
\begin{equation}\label{dislog} {\rm log}: C_r \isoarrow \ZZ/p^r\ZZ.
\end{equation}

On the other hand, there are two maps $p_i:  X_{0?}(qd) \ra X_?(d)$.  If we think of $X_{0?}(qd)$ as the moduli space of triples
$(E,C,\phi)$ where $E$ is a (generalized) elliptic curve, $C$ is a cyclic subgroup of order $q$, and $\phi$ is
some level structure at $d$, then $\pi_1$ is the map that forgets $C$ and $\pi_2$ is the map that replaces $E$ by $E/C$.  
\begin{prop}\label{t1q} (i)  The composition
$$H^0_0(X_?(d)_r,\omega_r) \overset{\pi_1^*}\to H^0_0(X_{0?}(d)_r,\omega_r) \overset{\cup \frS}\to H^1_0(X_{0?}(d)_r,\omega_r)
\overset{\pi_{2,*}}\to H^1_0(X_{?}(d)_r,\omega_r),$$
where the middle arrow identifies $\omega_r\otimes C_r$ with $\omega_r$ by means of the discrete logarithm,
coincides with the action of a (precise non-trivial) element $T^1_q \in \cH^1(G_q,K_q)$.  

(ii)  The action of $\cH^1(G_q,K_q; \Z/p^r\Z)$ on $H^*(X_?(d)_r,\omega_r)$ commutes with the action of the classical
Hecke operators at primes not dividing $qd$.
\end{prop}

Part (i) is equation (3-1) in \cite{HV}, while part (ii) was observed in Mazur's Eisenstein ideal paper \cite{Maz77} (and is easy to prove).

Analogous operators have been constructed and studied in the Columbia thesis of S. Atanasov \cite{Ata22} for Shimura varieties
attached to unitary groups of signature $(1,n-1)$.  These groups also have pairs of (nondegenerate) limits of discrete series
$\pi_0, \pi_1$ that have $(\frP,\fk)$-cohomology with coefficients in the same irreducible representation of (a maximal compact) $K$.

\subsubsection{The conjecture}

The action of the global Hecke algebra $\TT_{?d}$ on the weight one coherent cohomology spaces
$H^i_0(X_?(d),\omega)\otimes \overline{\QQ}_p$, $i = 0, 1$, is diagonalizable and
the characters $\lambda:  \TT_{?d} \ra \overline{\QQ}_p$ that appear in the decomposition are in bijection
with the cuspidal automorphic representations $\pi$ of $\GL(2,\aA)$ whose archimedean component
$\pi_\infty$ is isomorphic to the irreducible representation $I(0)$ introduced above.  As already mentioned, the same characters occur in $H^0_0$ and $H^1_0$, with the same multiplicity.  
It is also convenient to identify each $\pi$ that contributes to weight one coherent cohomology with its normalized newform $g_\pi$, which is uniquely determined by its classical $q$-expansion.  In the statement of the conjecture we will just write $g = g_\pi$ and we always assume $g$ to be a newform of level $d$.  

The analogue of Venkatesh's conjecture is easier to state in this setting because the hypothetical motive $M_\pi$ actually exists. By the Deligne-Serre \cite{DS74} Theorem each $\pi$ (or $g_\pi$) corresponds, in the sense of the Langlands correspondence, to a $2$-dimensional Galois representation
$$\rho_\pi:  \Gal(L/\QQ) \ra \GL(2,\ol{\QQ})$$
for some finite extension $L/\QQ$, where $\ol{\QQ}$ is any algebraically closed field.  
Each such $\rho_\pi$ is {\it odd}:  the determinant of $\rho_\pi(c)$, for any complex conjugation $c \in \Gal(L/\QQ)$, equals $-1$.  
Conversely, the Artin conjecture asserts that every $2$-dimensional odd  representation of the Galois group of a finite extension of $\QQ$ is of the form $\rho_\pi$; this has been known since the work of Khare--Wintenberger \cite{KW09} on Serre's Conjecture, although it had been established in most cases by Buzzard--Dickinson--Shepherd-Barron--Taylor \cite{BDST01}.  
Moreover, when the image of $\rho_\pi$ is a dihedral group the result is essentially due to Hecke.

Thus to each $\pi$ we have the $2$-dimensional representation $\rho_\pi$, which can be viewed as the Galois realization of a rank $2$ Artin motive, with coefficients in a ring $\cO$ of algebraic integers, which is a direct factor of the $0$-dimensional motive $R_{L/\QQ} \Spec(L)$.  
The motive $M_\pi$ is then the direct factor of $R_{L/\QQ} \Spec(L)$ whose Galois realization is given by the $3$-dimensional Galois representation 
$$\Ad^0\rho_\pi: \Gal(L/\QQ) \ra \GL(3,C),$$
the trace free summand of the $4$-dimensional representation $\Ad\rho_\pi$.   (In general $\Ad^0\rho_\pi$ factors through the Galois group of a proper subfield of $L$, but this is unimportant.)
It is explained in \cite[\S 2.7]{HV} that the motivic cohomology group $H^1_{\rm mot}(\QQ,M_{\pi}(1))$, denoted $\dnv_\pi$ above, can be replaced in the formulation of the conjecture by a certain $1$-dimensional space of {\it Stark units}, denoted $U_\pi$ or
$U_g$ if we write $g = g_\pi$.  
This is defined as
\begin{equation}\label{Ug}
U_g = U_\pi \stackrel{\sim}{\rightarrow} \operatorname{Hom}_{\mathcal{O}\left[G_{L / Q}\right]}\left(\operatorname{Ad}^0\rho, \mathrm{U}_{L} \otimes \mathcal{O}\right)
\end{equation}
This becomes one-dimensional upon tensoring with $E = \Frac(\cO)$.   One introduces a natural map, denoted
$$\mathrm{red}_q:  U_g \ra \mathbb{F}_q^\times \otimes \ZZ/p^r\ZZ.$$

The conjecture stated in \cite{HV} relates the action of the derived Hecke operator $T^1_q$ with a rational multiple of the discrete $\pmod{p^r}$ logarithm of a generator of $U_g$.  Since rational multiplication does not sit well with reduction modulo $p^r$, one needs to define the inputs and the relations carefully.   The Shimura class in $H^1(X_{0?}(qd)_r;C_r)$ defines a non-trivial class $\frS_\cO$ in $H^1(X_{0?}(qd)_r;\cO_r)$ (where $\cO_r$ is here the structure sheaf of $X_{0?}(qd)_r$) through the composition of \eqref{dislog} with the inclusion of Zariski sheaves $\ZZ/p^r\ZZ \hookrightarrow \cO_r$.  A feature specific to the case of $\GL_2$ and its inner forms is that the class $\frS_\cO$ is the reduction modulo $p^r$ of the class in $H^1(X_{0?}(qd),\cO)$
of (the complex conjugate of) an {\it Eisenstein cusp form} $f_E^\#$, in the sense of \cite{Maz77}.\footnote{Atanasov has shown in \cite{Ata22} that for the Shimura varieties attached to unitary groups of signature $(1,n-1)$ with $n > 2$, the analogues of the Shimura classes do not lift to characteristic zero. }  At least if the space of Eisenstein cusp forms is one-dimensional, this allows us to assign an invariant that measures the action of $T^1_q$.   Define $g^*$ to be the weight $1$ newform corresponding to the automorphic representation $\pi^\vee$ dual to $\pi$.  We let $g^{*,(q)}$ be the ``old form" in the space $\pi^{\vee}$ which in classical notation is just the analytic function $z \mapsto g^*(qz)$ for $z$ in the upper half-plane.  By Proposition \ref{t1q}, $T^1_q(g) \in H^1(X_{0?}(d)_r,\omega_r)$ has the same Hecke eigenvalues as $g$.  Hence in the Serre duality pairing
$$\langle,\rangle:  H^0_0(X_{0?}(d)_r,\omega_r)\otimes H^1_0(X_{0?}(d)_r,\omega_r) \ra H^1_0(X_{0?}(d)_r,\omega_r^{\otimes 2}) \ra \cO_r,$$
where the last arrow is the trace map, the map $\bullet \mapsto \langle \bullet, T^1_q(g)\rangle$ factors through
projection on the $g = g_\pi$-eigenspace of level $qd$.  Let now $f_E^\# = a_1\cdot f_E$ where $f_E$ is the normalized Eisenstein cusp form
(with leading coefficient $1$; the constant $a_1$ is essentially the inverse of what is called the {\it Merel constant} in \cite{HV}).  Then $T^1_q(g)$ is determined by the pairing $\langle gg^{*,(q)},f_E \rangle$, which can be identified with the Petersson pairing of
the weight $2$ forms $gg^{*,(q)}$ and $f_E$; $g^*$ has to be replaced by the old form because $f_E$ has level $q$ (the pairing is taking place in the second group from the right in Proposition \ref{t1q})  It can also be related to the square root of the central value of the triple product $L$-function $L(s,\pi\times \pi^\vee \times \pi(f_E))$, where the notation $\pi(f_E)$ is self-explanatory.  

The motivic conjecture, as formulated in \cite{HV}, comes down to an equality up to an integer factor:
  \begin{conj} 
  \label{conj:HV}
There exists an integer $m =m_g \geq 1$ and
  $u_g \in U_g$ such that,     for all  primes $q$ and $p$ as above 
$$ m \cdot \langle gg^{*,(q)}, \mathfrak{S} \rangle = \log( \mathrm{red}_q(u_g)).$$
 \end{conj}

Here both sides depend linearly on the choice of  $\log$  as in \eqref{dislog}.  The independence of $m$ from $q$ and $p$ is essential in this
statement, which would otherwise be vacuous.

\subsubsection{Results}\label{DHRVresults}

Weight $1$ eigenforms $g = g_\pi$ can be classified by the images of the corresponding (odd) $2$-dimensional Galois representations $\rho_\pi$, into those whose image is dihedral -- induced from a character $\psi_1$ of a quadratic extension $K/\QQ$ -- and the ``exotic" cases where the projectivized image is  $A_4$, $S_4$,
or $A_5$.  In the dihedral case the Stark unit group $U_\pi$ can be identified.  When $K$ is imaginary -- say $g$ or $\pi$ is of CM type -- $U_\pi$ consists of elliptic units;
when $K$ is real -- then $g$ or $\pi$ is of RM type -- $U_\pi$ is generated by the fundamental unit of $K$ itself.  Let $D$ be the discriminant of $K$.
Using the explicit determination of $U_\pi$, the article \cite{DHRV}
proves the following theorem:

\begin{thm}\label{DHRVthm}  If $K$ is imaginary, assume that $D$ is an odd prime and that
 $\psi_1$ is unramified.  If $K$ is real assume that $D$ is odd and that $\psi_1$  has conductor  dividing the different of $K/\QQ$. 
 Then Conjecture \ref{conj:HV} is true for $g_\pi$.
\end{thm}

The restrictions on $D$ and the conductor of $\psi_1$ are certainly not necessary; they have been relaxed in the Columbia thesis of R. Zhang \cite{ZhaThesis}, whose results are close to optimal. The methods of proof of Theorem \ref{DHRVthm} depend crucially to the realization of $g_\pi$ as explicit theta functions.   D. Marcil has provided striking numerical confirmation of the conjecture in some $A_4$ cases \cite{Ma21}.  As matters stand, it seems unlikely that one can prove unconditional results in the exotic case, though it's conceivable that the conjecture can be deduced somehow as a consequence of Stark's conjecture \cite{Stark1, Stark2, Stark3, Stark4}.

\subsection{Venkatesh's conjectures, $p$-adic cohomology, and the invariant $\ell_0$}

Specialists in $p$-adic modular forms already encountered the invariant $\ell_0$  as early as the 1990s. Hida's theory of ordinary  modular forms, or automorphic cohomology classes, for the reductive group $G$ over $\QQ$,
constructs a $p$-adic analytic space whose points parametrize $p$-adic modular forms, which can be viewed as the generic
fiber $\cX_{ord,G}$ of the formal spectrum $\Spf(\TT)$, where $\TT$ is Hida's ``big Hecke algebra."  The ring $\TT$ is an algebra over the Iwasawa algebra $\Lambda$
of the (maximal compact subgroup of the $p$-adic points of the) maximal torus of $G$.  In Hida's original theory for modular curves, and its extension to Shimura varieties, 
$\TT$ is a finite flat $\Lambda$-algebra, and the points of $\cX_{ord}$ corresponding to eigenvalues of classical (complex) Hecke algebras are dense in $\cE_{ord}$.
In general Hida conjectured in \cite{Hida} that the codimension of $\cX_{ord,G}$ in the rigid generic fiber of $\Spf(\Lambda)$ is exactly $\ell_0$.  

Urban made the analogous conjecture in the paper \cite{U} in which he constructed $p$-adic families of modular forms of finite slope, parametrized by the {\it eigenvariety} $\cX_G$, one of several constructions of a higher-dimensional generalization of the Coleman-Mazur eigencurve for elliptic modular forms.    In the finite slope situation, the eigenvariety is purely a rigid analytic object, Hida's Hecke algebra is replaced by a noetherian subring, which we also denote $\TT$, of 
the ring of global sections of $\cO_{\cX_G}$, and the map of $\Spf(\TT)$ to $\Spf(\Lambda)$ is replaced by a morphism $\pr:  \cX_G \ra \cW$.  Here  the {\it weight space} $\cW$ is again a rigid torus of dimension equal to the absolute rank of $G$, the Iwasawa algebra is replaced by the completed (noetherian) local ring
$\Lambda_w$ at a point $w \in \cW$.  Fixing a point $x \in \cX_G$, with $\pr(x) = w$,  the map $\pr$ makes the completion $\TT_x$  into a  finite $\Lambda_w$-algebra.

Both Hida and Urban refer to their conjectures on the codimension as non-abelian analogues of the Leopoldt Conjecture on the $p$-adic multiplicative independence of 
global units in number fields.  The Leopoldt Conjecture has for decades been the Bermuda triangle of algebraic number theory, tempting many of the most distinguished specialists into losing precious years in the pursuit of promising but ultimately futile attempted proofs.  Best not to try to prove the codimension conjectures of Hida and Urban, in other words.   Nevertheless, it's natural to reinterpret these codimension conjectures in  Venkatesh's motivic framework.  We briefly discuss two such interpretations.   The article \cite{HT17} of Hansen and Thorne approaches the eigenvarieties as parameter spaces for Galois representations, whereas 
the article \cite{KR} of Khare and Ronchetti provides evidence that the derived diamond algebras introduced in \S \ref{dda}, acting on Hida's
ordinary families, are the wild $p$-adic analogues of  Venkatesh's derived Hecke algebras at (tame) Taylor-Wiles primes.  In both constructions the motivic cohomology group   $H^1_{\rm mot}(\QQ,M_\pi(1))$ that figures in Venkatesh's motivic conjecture appears through its $p$-adic realization as the dual Selmer group.

\subsubsection{The Hansen-Thorne construction}

The article \cite{HT17} is concerned exclusively with eigenvarieties attached to the group 
of $G = \GL(n)_\QQ$.  The eigenvarieties are constructed by Hansen's method in \cite{Ha17}, following the topological approach of Ash and Stevens. By work of Harris--Lan--Taylor--Thorne \cite{HLTT}, and also by work of Scholze \cite{Sch15}, it is known that the classical point $x \in \cX := \cX_{GL(n)}$ corresponding to the cuspidal cohomological
automorphic representation $\pi$
has an associated semisimple $p$-adic Galois representation
\begin{equation}\label{GLL}
\rho_\pi:  \Gal(\ol{\QQ}/\QQ) \ra GL(n,L)
\end{equation}
for some finite extension $L/\QQ_p$ over which $x$ is defined.   The point $x$ also depends on the secondary datum of a {\it refinement} $t_x$, or an ordering on the eigenvalues of crystalline Frobenius, that is assumed implicitly in what follows.  The representation $\rho_\pi$ does not depend on the choice of $t_x$.

We then define the dual adjoint Selmer group $H^1_f(\QQ,M_\pi(1)) = H^1_f(\QQ,\Ad\circ \rho_\pi(1))$ as in \S \ref{actmot}.   Note that this is a
$\QQ_p$-vector space; the constructions in \cite{HT17}, in contrast to those of \cite{GV, V19, HV}, are carried out entirely in characteristic zero.  

The main result of \cite{HT17} is that, at least under favorable conditions, and assuming a conjecture recalled in the statement of Theorem \ref{HTtheorem} below, $H^1_f(\QQ,M_\pi(1))$ controls the 
infinitesimal geometry of $\cE$ at $x$ with respect to the morphism $\pr$.  
Combining this with the universal coefficients spectral sequence constructed in \cite{Ha17}, the authors recover 
the analogue of  the Galatius-Venkatesh Theorem \ref{GVtheorem}, but without introducing the derived deformation ring.

The main hypothesis required for the final results of \cite{HT17} is a version of the isomorphism 
\begin{equation}\label{ReT} R_x \isoarrow \TT_x
\end{equation} 
of the Taylor-Wiles method.  Here $R_x$ is an {\it underived} deformation ring defined using the theory of $(\varphi,\Gamma)$-modules over the Robba ring.  In particular, $R_x$ is defined purely by characteristic zero methods, and (pro)-represents a functor on Artinian local $L$-algebras for the $p$-adic coefficient field $L$ of $\rho_\pi$ \eqref{GLL}.  More precisely, it is assumed that $\rho_\pi$ is absolutely irreducible; then $R_x$ (pro)-represents {\it trianguline} deformations of $\rho_\pi$ along with the refinement $t_x$. 
The property of being \emph{trianguline}\footnote{pronounced in English to rhyme with {\it fine}, {\it wine}, and {\it turpentine}} \cite{Ber11} is determined by the $(\varphi,\Gamma)$-module attached to the restriction of $\rho_\pi$ to decomposition groups at $p$; the relevant fact for us is that the Galois representations attached to points on the eigenvariety, including $\rho_\pi$ itself, necessarily have this property.  

We now assume $\pi$ contributes to $H^*(S_{K_0},\tilde{V})$, with notation as in \S \ref{DHAcoh} for $G = \GL(n)_\QQ$.  We change notation slightly,
and let $\mathfrak{m}_\pi$ denote the maximal ideal corresponding to $\pi$ in $\TT_{K_0}\otimes \QQ_p$;  $H^*(S_{K_0},\tilde{V})_{\fm_\pi}$
denote the localization.   Let $q_0 = q_0(\GL(n)_\QQ)$, $\ell_0 = \ell_0(\GL(n)_\QQ)$.   The following theorem summarizes the main results of Theorems 4.9 and 4.13 of \cite{HT17}.

\begin{thm}\label{HTtheorem}  Let $x \in \cX$ be a classical point as above, with $w = p(x)$,
and let $\lambda:  \TT \ra L$ be the corresponding character of the Hecke algebra.
Assume $\lambda$ is {\rm numerically noncritical}\footnote{This is a regularity condition, for which we refer to \cite[Theorem 4.7]{HT17}.}    Then 
\begin{enumerate}
\item  $\dim \TT_x \geq \dim \Lambda_w - \ell_0$;
\item  If $\dim \TT_x = \dim \Lambda_w - \ell_0$ let $V_x = \ker(\Lambda \ra \TT_x)\otimes_{\Lambda_w} L$, where the map from $\Lambda_w$ on $L$ is the augmentation.  Then $H^*(S_{K_0},\tilde{V})_{\fm_\pi}$ is free of rank $1$ as graded module over the exterior algebra $\wedge^*_L V_x$ and is generated  by its term in top degree $q_0 + \ell_0$.
\item  Assume $\rho_\pi$ is absolutely irreducible and crystalline at $p$, with the expected Hodge-Tate weights, and its associated Weil-Deligne representation at every prime $\ell$ (including $p$) corresponds to $\pi_\ell$ under the local Langlands correspondence.
Assume also the isomorphism \eqref{ReT}.  

Then $V_x$ can be identified with $H^1_f(\QQ,\Ad \rho_\pi(1)) = \dnv_{\pi,p}[\frac{1}{p}]$, and 
$H^*(S_{K_0},\tilde{V})_{\fm_\pi}$ is free of rank $1$ as graded module over the exterior algebra $\wedge^*_L \dnv_{\pi,p}$, by its term in top degree $q_0 + \ell_0$.

\end{enumerate}

\end{thm}

\begin{remark}  It is natural to assume that the action of $\wedge^*_L \dnv_{\pi,p}$ on $H^*(S_{K_0},\tilde{V}))_{\fm_\pi}$ coincides, up to duality, with the one defined in Theorem \ref{GVtheorem}.  However, the identification of $V_x$ with $\dnv_{\pi,p}[\frac{1}{p}]$ depends on the deformation ring $R_x$, which in turn depends on the choice of refinement $t_x$.  Hansen has conjectured that the action does not depend on this additional choice, but even this conjecture remains open.\footnote{We thank David Hansen for pointing this out.}
\end{remark}

\subsubsection{The Khare-Ronchetti construction}\label{dda}

The constructions in this section are based on the paper \cite{KR}. Although their derived diamond algebras are defined purely locally, the construction is bound up with the action on
the ordinary part of global cohomology.  It is not clear to us how the representation theory of the derived diamond algebras fits in a hypothetical $p$-adic
local Langlands program.

Let $p > 2$ be a prime number and let $G$ be the group $\PGL(n)$ over a finite extension $F/\Q$.  For any open compact subgroup $K \subset G$ and any ring $A$ we let $H(G,K;A)$ denote the (underived) Hecke algebra of $K$-biinvariant functions on $G$ with coefficients in $A$.
Let $B = T\cdot U$ be the upper-triangular Borel subgroup, $U$ its unipotent radical, and $T$ the diagonal maximal torus.   We restrict attention to $\PGL(n)$ in order to use explicit matrix groups, and to make use of the known construction of Galois parameters attached to cuspidal cohomology classes.  The article \cite{KR} treats the general case and the reader should be able to make the appropriate adjustments.

Let $U^-$ be the opposite unipotent (lower triangular) subgroup.  Let $v$ be a prime of $F$ dividing $p$, with integer ring $\cO = \cO_v$ and uniformizer $\varpi = \varpi_v$. 
For $1 \leq b \leq c$ we let 
\begin{equation}\label{Ibc}
I_{b,c,v} = U^-(\varpi^c\cO)\cdot T(1+\varpi^b\cO) \cdot U(\cO);
\end{equation}
 i.e. 
the generalization from $\GL_2$ to $G$ of the congruence subgroup
$$I_{b,c,v} = \begin{pmatrix} 1 + \varpi^b* &  * \\ \varpi^c* & 1+\varpi^b*\end{pmatrix} \subset \GL(2,\cO).$$
We let $I_{b,c} = \prod_{v \mid p} I_{b,c,v}$.

Let $K = K_p \cdot K^p \subset G(\A_f)$, where $K^p$ is a fixed level subgroup away from $p$
and $K_p$ varies over the subgroups $\prod_{v\mid p} I_{b,c,v}$ introduced in \eqref{Ibc}.   
We consider the cohomology of the locally symmetric spaces of dimension $d$
$$Y_{b,c} = G(\QQ)\backslash G(\A_f)/K_\infty\cdot K^p\cdot  I_{b,c}$$
where $K_\infty \subset G(\RR)$ is maximal compact and $K^p$ is a fixed level subgroup away from $p$,
so that $K^p\cdot I_{b,c} \subset G(\A_f)$ is open compact.  We assume $K^p$ is hyperspecial maximal
outside a finite set of primes $S$.

The cohomology of $Y_{b,c}$ is a module over the unramified Hecke algebra 
$$H^{un} = \prod_{q \notin S, q \neq p} H(G_q,K_q).$$
 It is also  a module over a derived Hecke algebra at $p$, which we proceed to define in the setting of Hida
 theory.
 
Following Hida, we define an operator $U_p = \prod_{\alpha \in \Sigma} U_\alpha$ acting as correspondences on each $Y(K(b,c))$, compatibly with
the natural maps $Y(K(b,c)) \ra Y(K(b',c'))$ for $b' < b, c' < c$.    (There are really operators $U_p(b,c)$ but we don't stress this.)
Roughly, with $G = \PGL(n)$  we have  
$$U_p = \prod_{j = 1}^{n-1} U_{j,p}, ~ U_{j,p} = U_{\alpha_j} = I_{b,c} \prod_{v \mid p}{\rm diag}(\varpi_v \Id_j, \Id_{n-j})I_{b,c}.$$  

\begin{defn}  The ordinary part $H^i(Y_{b,c},\Zn)_{ord} \subset H^i(Y_{b,c},\Zn)$ is the submodule (and direct summand) where the $U_p$ operator
acts invertibly.
\end{defn}

\begin{fact}\label{bbc} For all $c \geq b \geq 1$, $n \geq 1$, the pullback
$$H^i(Y_{b,b},\Zn)_{ord}  \ra H^i(Y_{b,c},\Zn)_{ord}$$
is an isomorphism.
\end{fact}
 
\begin{proof}  The proof in \cite{KR} is a calculation with double cosets acting on singular cochains.
\end{proof}

\begin{remark}  If Fact \ref{bbc} admits a generalization to overconvergent cohomology it will necessarily need to be formulated
in characteristic zero.
\end{remark}

The group $I_{1,c}/I_{b,c}$ is isomorphic to the abelian $p$-group 
$$T_b := \prod_{v \mid p}T(1+\varpi_v\cO_v)/T(1+\varpi^b_v\cO_v).$$  
  Since $I_{b,c}$ is normal in $I_{1,c}$, the elements
$\alpha \in I_{1,c}$ define {\it diamond operators} 
$$\da = I_{b,c}\alpha I_{b,c} \in H_{\Zm}(G,I_{b,c}).$$ 
for all $b,c$. 

The action of $\da$ on $H^i(Y_{b,c},\Zn)$ 
commutes with the maps $Y_{b,c} \ra Y_{b',c'}$ for $b' < b, c' < c$.  
Let $\Lambda_b = \ZZ_p[T_b]$, $\Lambda = \varprojlim_b \Lambda_b$, which is the Iwasawa algebra of $T(\ZZ_p)$.   
We also use $\Lambda_F = \Lambda\otimes_{\ZZ_p} F$ for any finite extension $F/\QQ_p$.

The $T_{b,c}$-covering $Y_{c,c}/Y_{b,c}$ corresponds to a morphism to the classifying space
$$\pi_{b,c}: Y_{b,c} \ra BT_{b,c}.$$  
Thus there is a ring homomorphism
$$\pi_{b,c}^*:  H^*(BT_{b,c},\Zm) = H^*(T_{b,c},\Zm) \ra H^*(Y_{b,c},\Zm)$$  
and thus an action of $H^*(T_{b,c},\Zm)$ on $H^*(Y_{b,c},\Zm)$ by {\it derived diamond operators}.    Since $I_{b,c}^{ab} = T_{b,c}$, 
we have an isomorphism
$$\iota:  H^1(T_{b,c},\Zm) \isoarrow H^1(I_{b,c},\Zm) \subset \Cal{H}_{\Zm}(G_p,I_{b,c}).$$ 

\begin{lemma}[\cite{KR}]  The derived diamond action of  $\phi \in H^1(T_{b,c},\Zm)$ coincides with the derived Hecke 
action of $\iota(\phi)$ (obtained as in Fact \ref{DHAaction}). 

In particular, when $b = 1$ this action commutes with the $U_{j,p}$-operators.
\end{lemma}
The claim about commutation follows from an explicit double coset calculation.   
It implies (when $Y_{b,c} = Y_{1,c}$) 
\begin{cor}  The derived degree $1$ diamond operators act on $H^*(Y_{1,c},\Zm)_{ord}$, preserving all
$U_{j,p}$-eigenspaces.
\end{cor}  
Another calculation shows:
\begin{prop}  As $c$ varies, the actions of $H^1(T_{1,c},\Zm)$ on $H^*(Y_{1,1},\Zm)_{ord}$ are compatible, and thus
give rise in the limit to an action of $H^1(T(1+p\Zp)),\Zm)$ on $H^*(Y_{1,1},\Zm)_{ord}$.  

Letting $m \to \infty$ we thus obtain a graded action of $H^*(T(1+p\Zp),\Zp)$ on $H^*(Y_{1,1},\Zp)_{ord}$.
\end{prop}  

We return to the setting of \S \ref{derham}.  Let $\pi$ be a cuspidal cohomological representation of $G$, as in that section.  Thus we
have, as in Corollary \ref{extactglobal}
$$H^{q_0(G)+\bullet}_0(Y_{1,1},\ol{\QQ}_p)[\pi_f]$$
 is a free differential graded module over ($p$-adic) $\wedge^{\bullet}(\mathfrak{a}^*)$ of rank $0$ or $1$;
 here we are using multiplicity one, and the existence of Whittaker models, for $\PGL(n)$.

The following theorem is an imprecise version of the main resuilt of \cite{KR}.    Let $E$ be an appropriate $p$-adic coefficient field for $\pi_f$.

\begin{thm}[\cite{KR}]\label{KRaction}  Under certain favorable circumstances, $H^*(Y_{1,1},E)_\pi$ is generated over $H^{q_0}(Y_{1,1},E)_\pi$ by the action of the derived diamond operators $H^*(T(1+p\Zp),E)$.

Moreover, there is a canonical map from $H^1(T(1+p\Zp),E)$ to the Selmer group $\dnv_{\pi,p}[\frac{1}{p}]^\vee$ attached to $\pi$.   And under even more favorable circumstances,
the action of $H^*(T(1+p\Zp),E)$ on $H^*(Y_{1,1},E)_\pi$ factors through this Selmer group.
\end{thm}

The ``favorable circumstances" to which the theorem alludes include the isomorphism \eqref{ReT}, in other words the non-abelian Leopoldt conjecture.  
As in Theorem \ref{GVcompatible}, the actions of Theorem \ref{KRaction} and Theorem \ref{HTtheorem} should be dual to each other.

\section{A selection of open problems}\label{sec: open problem}

We discuss some open problems at the interface of derived algebraic geometry and the Langlands correspondence.

\subsection{Derived local deformation conditions}\label{ssec: local conditions}
The context for this problem is \S \ref{sec: derived moduli} and \S \ref{DGdr}. In \S \ref{sec: derived moduli} we discussed examples of moduli spaces that we know how to upgrade to derived moduli spaces. However, there are also many interesting examples, arising in connection with concrete problems, for which we do not know how to construct any good derived enhancement, and for which it is unclear to the authors whether such an enhancement should even exist. The examples of Galois deformation functors \emph{with local conditions} are particularly interesting, since they would presumably have consequences for automorphy lifting. 

In practice, when implementing modularity lifting arguments, one wants to impose more general local conditions on Galois deformation functors, as explained for example in the articles of \cite{CS23} in these proceedings. This is done by constructing Galois deformation functors for local fields with local conditions, and it turns out that these can often fail to be LCI; see \cite[Theorem 3.19]{CG18} and \cite{Sn18} for some ``$\ell=p$'' examples, and \cite{BKM} for $\ell \neq p$ examples. This is an interesting phenomenon to understand, because Kisin's modification of the Taylor-Wiles method ``reduces'' automorphy lifting to the problem of having good control of \emph{local} deformation rings.

We do not know if these local deformation functors should actually be subject to the hidden smoothness philosophy, since these are typically not ``honest'' moduli functors, in the sense that their moduli descriptions are not known ``concretely''. Rather, they are typically constructed from ``honest'' moduli spaces by processes such as formation of Zariski closure or scheme-theoretic image, which do not admit good derived versions. 

For example, the local deformation functors at $p$ are often constructed by the following sort of procedure (which to our knowledge was first employed in \cite{CDT}): define a subset of $\ol{\Q}_p$-points of the unrestricted deformation functor using $p$-adic Hodge theory, and take their Zariski closure. (The point here is that the $p$-adic Hodge theory is understood most generally in the context of representations with characteristic zero coefficients.) Then the characteristic zero points admit some moduli theoretic description in terms of rational Galois representations with $p$-adic Hodge theoretic conditions, but it is a priori unclear how to interpret the $\Z_p$ or $\F_p$-points concretely in terms of Galois representations. 

There are some cases where the theory is well-behaved with integral coefficients, for example in the Fontaine-Laffaille range. In this case, there is a full moduli-theoretic description of the local deformation functor, which even makes clear that it is smooth. 

In general, the local deformation functors at $p$ are constructed by Kisin \cite{Kis08}, but their construction again involves processes such as formation of scheme-theoretic image or Zariski closure from the generic fiber, which we do not know how to perform in a derived way. This is why we do not know natural ways to derive these local deformation functors. However, in the past few years there have been major advances in our understanding of integral $p$-adic Hodge theory due to Bhatt-Morrow-Scholze and Bhatt-Scholze. In particular, Bhatt-Scholze have established in \cite{BS21} a new interpretation of lattices in crystalline Galois representations, in terms of prismatic $F$-crystals. While this itself does not suffice to give a notion of crystalline Galois representations with general coefficients, it is a promising step towards such a concept, and therefore towards a moduli-theoretic description of crystalline deformation functors, which could then be derived. 

We caution, however, that even if this is all possible, it is not quite clear at present what applications this would unlock. Global methods, such as the Taylor-Wiles method, are predicated upon the premise that the structure of global deformation rings and automorphic forms can be ``simplified'' in a suitable sense by adding level structure. However, this simplification occurs completely in the local-to-global aspect, and it is not clear that one can hope for the same type of simplification if local deformation rings are non-classical. Although it does not directly tackle this concern, the work \cite{IKM22} of Iyengar-Khare-Manning, which generalizes Wiles' \emph{numerical criterion} using derived algebra, is a promising step in this broad direction.

\subsection{Differential graded Hecke algebra}\label{ssec: DGHA}
The context of this discussion is \S \ref{ldHa}. In \S \ref{sssec: commutativity}, we explained that commutativity of the local derived Hecke algebra $\cH(G(K); \Lambda)$ is still unknown. We conjecture that it is true in general, and moreover should be a consequence of ``higher commutativity'' for the differential graded Hecke algebra $\DGHA(G(K); \Lambda)$. For differential graded algebras, commutativity is a structure rather than a property, and there is in fact a spectrum of structures that measures ``how commutative'' an algebra is, called $\EE_n$-structures. An associative differential graded algebra has an $\EE_1$-structure, while an $\EE_{\infty}$-structure is the homotopy-coherent analogue of commutativity. It is clear from the construction that $\DGHA(G(K); \Lambda)$ has an $\EE_1$-structure. If this $\EE_1$-structure can be promoted to an $\EE_2$-structure, then it would follow that $\cH(G(K); \Lambda)$ is commutative. We believe that this is the case, and we conjecture moreover that the $\EE_1$-structure on $\DGHA(G(K); \Lambda)$ can be promoted to an $\EE_3$-structure.

We explain the motivation for this conjecture. Unpublished work of Feng-Gaitsgory \cite{FG} implies (conditionally) that if $K$ is a function field, then $\cH(G(K); \Lambda)$ is commutative whenever $\Lambda$ is of characteristic $\ell$ which is not too large with respect to $G$. The strategy of \cite{FG} is to apply the \emph{categorical trace of Frobenius} to the \emph{modular derived Satake equivalence}, a result which has not appeared in the literature but which has been announced by Arinkin-Bezrukavnikov, and may well be possible to deduce directly from the recent paper \cite{BR24} of Bezrukavnikov-Riche. 

This strategy is not merely an elaboration of a proof of commutativity for the classical Hecke algebra, but rather proceeds from a categorical equivalence. It should also clarify the ``meaning'' of the derived Hecke algebra in terms of the dual group, but the answer is complicated to state, so we do not formulate it here. It seems likely that a variant of the strategy will also work for $p$-adic local fields $K$, although significant foundational groundwork would need to be laid for this to happen. 

Now, we explain why this suggests the possibility of an $\EE_3$-structure on $\DGHA(G(K);\Lambda)$. Within the above strategy, the differential graded Hecke algebra appears as the categorical trace of the derived Satake category. This latter object is expected to admit an $\EE_3$-structure, with one $\EE_1$-monoidal operation coming from convolution, and a compatible $\EE_2$-monoidal operation related to its structure as a ``factorization category'' over an algebraic curve, which is an algebraic analogue of an $\EE_2$-structure. 

\subsection{Action of the derived Galois deformation ring}\label{ssec: construct action} The context for this discussion is \S \ref{DGdr}. There we sketched the main result of \cite{GV}, which is the construction of an action of \emph{the homotopy groups of} a derived Galois deformation ring $\cR^{\ol{\rho}, \crys}_{\cO_F[1/S]}$ on the \emph{homology groups of} a locally symmetric space (localized at a corresponding maximal ideal of the Hecke algebra). 

The above action occurs at the level of classical algebra: a graded commutative ring acts on a graded abelian group. We expect that this action can be refined to an action of the derived Galois deformation ring $\cR^{\ol{\rho}, \crys}_{\cO_F[1/S]}$ on the localization of the homology chains (i.e., an action of an animated commutative ring on an animated abelian group), which induces the previous action after taking homotopy groups. 

The construction of this action presents an intriguing challenge. Recall that the Galatius-Venkatesh construction uses Taylor-Wiles patching, which is based on adding level structure \emph{away} from $p$. Morally, Galatius-Venkatesh find that the derived Galois deformation ring should become classical when patched, and thus deduce the derived action from the classical action, but this does not literally make sense since there is no patched derived Galois deformation ring. Instead, one could hope to add some kind of infinite level structure at $p$ and hope that the corresponding derived Galois deformation ring is classical. Indeed, the work of Khare-Ronchetti \cite{KR} suggests that this should happen when ascending the ordinary (Hida) tower, if one considers ordinary deformations. On the other hand, \cite{KR} also points out that this approach is intimately connected with $p$-adic transcendence problems which are known to be extremely difficult, such as the Leopoldt Conjecture, so perhaps one should look for some other angle. 

Although the construction of this refined action seems like a homotopy-theoretic problem, it would have concrete consequences for classical algebra. In fact, by considerations with the spectral Hecke algebra \cite{F20}, it would actually imply Venkatesh's reciprocity law \cite[Theorem 8.5]{V19} for the action of the derived Hecke algebra. This is very desirable, since in the current formulation of the Theorem from \emph{loc. cit.} there is no effective way to determine for which primes $q$ it applies.

\subsection{The Trace Conjecture for Hitchin stacks}\label{ssec: trace conjecture} The context for this discussion is \S \ref{sec: derived special cycles}. The following problem is extracted from \cite[\S 6.4]{FK}. Let
\begin{equation}\label{eq: general correspondence}
\begin{tikzcd}
\cM & \ar[l, "{c_0}"'] \cC \ar[r, "{c_1}"] &  \cM
\end{tikzcd}
\end{equation}
be a correspondence of (derived) stacks over $\F_q$. Define $\Sht_{\cM}$ by the derived Cartesian square 
\begin{equation}\label{eq: sht diagram}
\begin{tikzcd}
\Sht_{\cM} \ar[r, "\Delta'"] \ar[d, "c'"] & \cC \ar[d, "{(\Frob \circ c_0, c_1)}"] \\
\cM \ar[r, "\Delta"] & \cM \times \cM
\end{tikzcd}
\end{equation}

In \cite[Lemma 4.2.1]{FYZ3}, the following is proved: 

\begin{lemma}\label{lem: tangent complex of ShtM}
The tangent complex $\bT_{\Sht_{\cM}/\F_q}$ is the restriction of $\bT_{c_1}$. In particular, if $c_1$ is quasismooth, then $\Sht_{\cM}$ is quasismooth (over $\Spec \F_q$), of dimension equal to $d(c_1)$ (the relative virtual dimension of $c_1$). 
\end{lemma}

We will formulate the \emph{Trace Conjecture} for motivic sheaves. See \cite[\S 3]{FK} for the relevant definition of motivic sheaves; the reader unfamiliar with motivic sheaf theory may replace $\Q$ with the constant sheaf $\Ql$ in the following discussion without losing the main gist. 

For a motivic sheaf $\cK$ and $i \in \Z$, we write $\cK \tw{i} := \cK[2i](i)$. For a map of (derived) stacks $\cX \rightarrow \cY$, we denote $\CH_{i}(\cX/\cY) := H^{-2i}(\cX, f^! \Q_{\cY}(-i))$. When $f \co \cX \rightarrow \Spec k$ is the structure map to a field, we abbreviate $\CH_i(\cX) := \CH_i(\cX/\Spec k)$. If $f$ is quasismooth, then we write $d(f) = \chi(\bT_f)$ for the Euler characteristic of its tangent complex, and call it the ``relative (virtual) dimension'' of $f$. By \cite{KhanI}, in this situation there is a relative fundamental class $[f] \in \CH_{d(f)}(\cX/\cY)$. 

Assume that in the above setup $c_1$ is quasismooth, so that $[c_1] \in \CH_{d(c_1)}(\cC/ \cM)$ exists. Then Lemma \ref{lem: tangent complex of ShtM} implies that $\Sht_{\cM}$ is quasismooth (over $\Spec \F_q$), so the derived fundamental class $[\Sht_{\cM}]\in \CH_{d(c_1)}(\Sht_{\cM})$ exists. 

On the other hand, regarding $[c_1]$ as a map $c_1^* \Q_{\cM} \rightarrow c_1^! \Q_{\cM} \tw{-d(c_1)}$, we have a sequence of maps
\[
(\Frob \circ c_0)^* \Q_{\cM} = \Q_{\cC}  = c_1^* \Q_{\cM} \xrightarrow{[c_1]} c_1^! \Q_{\cM} \tw{-d(c_1)},
\]
whose composition we call $\cc_{\cM}$. This is a \emph{cohomological correspondence} between $\Q_{\cM}$ and $\Q_{\cM}\tw{-d(c_1)}$. We will define another class $\Tr(\cc_{\cM})  \in \CH_{d(c_1)}(\Sht_{\cM})$. Consider \eqref{eq: sht diagram}, and abbreviate $c$ for the right vertical map. It is explained in \cite[\S 3.8]{FK} that 
\begin{align}\label{eq: tr eq 1}
\cRHom((\Frob \circ c_0)^* \Q_{\cM}, c_1^! \Q_{\cM})  & \cong  c^! \cRHom(\pr_0^* \Q_{\cM}, \pr_1^! \Q_{\cM}) \nonumber \\
& \cong c^! (\DD(\Q_{\cM}) \boxtimes \Q_{\cM})
\end{align}
where $\DD$ is the Verdier duality functor. The evaluation map $\DD(\Q_{\cM}) \otimes \Q_{\cM}\tw{-d(c_1)}  \rightarrow \DD_{\cM} \tw{-d(c_1)}$ is adjoint to a map $\DD(\Q_{\cM}) \boxtimes \Q_{\cM} \rightarrow \Delta_* \DD_{\cM} \tw{-d(c_1)}$ where $ \DD_{\cM}$ is the dualizing complex of $\cM$. Composing this with \eqref{eq: tr eq 1} gives a map 
\begin{equation}\label{eq: tr eq 2}
\cRHom((\Frob \circ c_0)^* \Q_{\cM}, c_1^! \Q_{\cM})  \rightarrow c^! \Delta_* \DD_{\cM}\tw{-d(c_1)}.
\end{equation}
Finally, using proper base change, we have isomorphisms 
\begin{equation}\label{eq: defining tr 3}
c^! \Delta_* \DD_{\cM}\tw{-d(c_1)} \cong  \Delta'_* (c')^! \DD_{\cM} \tw{-d(c_1)}\cong   \Delta'_* \DD_{\Sht_{\cM}} \tw{-d(c_1)}. 
\end{equation}
We may regard $\cc_{\cM}$ as a global section of $\cRHom(c_0^* \Q_{\cM}, c_1^! \Q_{\cM} \tw{-d(c_1)})$. Then 
\[
\Tr(\cc_{\cM}) \in H^0( \cC, \Delta'_* \DD_{\Sht_{\cM}} \tw{-d(c_1)}) \cong \CH_{d(c_1)}(\Sht_{\cM})
\]
is defined as its image under the composition of \eqref{eq: tr eq 2}, and \eqref{eq: defining tr 3}.

At this point, we have two natural classes in $\CH_{d(c_1)}(\Sht_{\cM})$. One is the derived fundamental class $[\Sht_{\cM}]$, and the other is $\Tr(\cc_{\cM})$. It is then natural to ask if
\[
 \Tr(\cc_{\cM}) = [\Sht_{\cM}] \in \CH_{d(c_1)}(\Sht_{\cM}).
\]
We have no evidence to believe that it is true in the full stated generality. However, for $\cM$ a derived Hitchin stack of the type considered in Example \ref{ex: hitchin 3}, we feel that the equality is probably necessary for the Modularity Conjecture \cite[Conjecture 4.12]{FYZ2} to be true, so we conjecture the equality to be true in this case. 

\begin{conj}[Trace Conjecture]
For $\cM$ a derived Hitchin stack as in \cite[\S 5]{FYZ2}, we have
\[
 \Tr(\cc_{\cM}) = [\Sht_{\cM}] \in \CH_{d(c_1)}(\Sht_{\cM}).
 \] 
\end{conj}

Such $\cM$ have special properties which could conceivably be necessary for the proof; for example, they are quasismooth\footnote{Note that this quasismoothness condition is not relevant for the formulation of the Trace Conjecture.}. We expect that the Trace Conjecture should be an important ingredient for a proof of the Modularity Conjecture \ref{conj: modularity theta}.

\appendix 

\section{A crash course on simplicial commutative rings}\label{appendix}

This appendix is a sketchy primer on simplicial commutative rings, leading up to the cotangent complex. Its purpose is to provide more grounding for the readers to whom the ``black box'' approach of the main text is too vague, and to indicate references for those interested in delving further into the subject.  For short introductions we like \cite{Qu68, Qu70, Iy07}; longer textbooks include \cite{May67} and \cite{GJ09}.

In particular, we think it is helpful for readers to learn about simplicial commutative rings and model categories, even if their ultimate goal is work in the language of animated rings and $\infty$-categories.

Our treatment is extremely far from comprehensive. In fact, we view the sketchiness of our writeup as its main (and perhaps only) feature. Each of the subjects we touch upon is already well-documented in the literature, but simplicial homotopy theory is so vast and technical a subject that there is real risk in getting lost in the details. In keeping with the spirit of the main text, we prefer merely to sketch the skeletal structure of the subject, contenting ourselves with slogans and intuition, and deferring details to our favorite references.

\subsection{Why ``simplicial''?} We will introduce ``simplicial sets'' as a combinatorial model for topological spaces, and ``simplicial commutative rings'' as a combinatorial model for topological commutative rings. 

Why do we say ``simplicial'' instead of ``topological''? One answer is that the combinatorial nature of simplicial objects makes them better-behaved from a technical perspective -- it rules out ``pathological'' topological spaces. We do not have a good philosophical answer; the reader could consult \cite{58497} for more musings on this question. 

Furthermore, it turns out that the adjective ``simplicial'' is flexible enough to be applied to any object of any category, and provides a way to do abstract homotopy theory in many examples of interest. That is, if $\msf{C}$ is a category of ``widgets'', we will be able to define a category $\msf{sC}$ of ``simplicial widgets''.

\subsection{The simplex category}\label{ssec: simplex category}

The \emph{simplex category} $\Delta$ has as its objects the sets $[n] := \{0, 1, \ldots, n\}$ for $n \geq 0$, and $\Hom_{\Delta}([n],  [m])$ consists of all maps of sets $f \co [n] \rightarrow [m]$ that do not reverse order, i.e., if $i \leq j$ then $f(i) \leq f(j)$. 

We will name specific generating morphisms between $[n]$ and $[n+1]$. For $0 \leq i \leq  n+1$, the \emph{coface map} 
\[
\delta_i \co [n] \rightarrow [n+1]
\]
is the morphism ``skipping over $i$'', or more formally 
\[
\delta_i(j) =\begin{cases} j & j<i, \\ 
j+1 & j \geq i. \end{cases}
\]
For $0 \leq i \leq n$, the \emph{codegeneracy map} 
\[
\sigma_i \co [n+1] \rightarrow [n]
\]
is the morphism ``doubling up on $i$'', or more formally 
\[
\sigma_i(j) =\begin{cases} j & j \leq i, \\ 
j-1 & j \geq i+1. \end{cases}
\]

\begin{remark}\label{rem: simplex morphism generators}
The sense in which these morphisms ``generate'' is that any morphism in $\Delta$ can be written as a composition of codegeneracy maps followed by a composition of coface maps. 
\end{remark}

\subsubsection{Simplicial sets}

A \emph{simplicial set} is a functor $\ul{X} \co \Delta^{\op} \rightarrow \msf{Sets}$. Concretely, it is specified by a collection of sets $X_n  := \ul{X}([n])$, for $n \geq 0$, with maps between them indexed by the morphisms in $\Delta^{\op}$. In particular, we have \emph{face maps}
\[
d_i := \ul{X}(\delta_i) \co X_{n+1} \rightarrow  X_n 
\quad \text{for $0 \leq i \leq n+1$},
\]
and \emph{degeneracy maps}
\[
s_i := \ul{X}(\sigma_i) \co X_{n} \rightarrow  X_{n+1} \quad 
\text{for $0 \leq i \leq n$}.
\]
Because of Remark \ref{rem: simplex morphism generators}, to specify a simplicial set it suffices to specify the data of the $\{X_n\}_{n \geq 0}$ and the $d_i$ and $s_i$ for each $n$, satisfying the ``simplicial identities'' \cite[I.1]{GJ09}:
\begin{itemize}
\item $d_id_j = d_{j-1}d_i$ for $i<j$. 
\item $d_js_j = d_{j+1}s_j  = \Id$. 
\item $d_i s_j = s_{j-1} d_i$ for $i<j$. 
\item $d_is_j = s_j d_{i-1}$ for $i>j+1$. 
\item $s_i s_j  = s_{j+1} s_i$ for $i \leq j$. 
\end{itemize}

\begin{example}
Each element $[n] \in \Delta$ induces a representable functor on $\Delta^{\op}$, namely $\Hom_{\Delta}(-, [n])$, and so induces a simplicial set that we shall call $\Delta^n$. Intuitively, we think of it as corresponding to the topological space of the standard $n$-simplex, denoted 
\[
|\Delta^n| := \{ (x_0, \ldots, x_n) \in \R^{n+1} \co 0 \leq x_i \leq 1, \sum x_i = 1\}.
\]
\end{example}

\begin{example}
The definition of simplicial set is an axiomatization of the structure obtained from a topological space $X$ by setting $\ul{X}[n]$ to be the set of continuous maps from the $n$-simplex $|\Delta^n|$ to $X$, i.e. $X_n$ is the set of ``$n$-simplices of $X$''. The face maps $d_i$ are the usual boundary maps. The degeneracy maps account for ways to view lower-dimensional simplices as ``degenerate'' instances of higher-dimensional simplices; these do not play an important role in the topological theory, but are technically convenient in the simplicial theory. This defines the \emph{singular simplices} functor 
\[
\Sing \co \msf{Top} \rightarrow \msf{sSets}.
\]
Therefore, for a general simplicial set $\ul{X}$ we will refer to $\ul{X}[n]$ as the ``$n$-simplices of $\ul{X}$''. 
\end{example}

\subsubsection{Geometric realization}
Given a simplicial set $\ul{X}$, we can build a topological space $X := |\ul{X}|$ from $\ul{X}$ by viewing $\ul{X}$ as a recipe for assembling a simplicial complex. In formulas, the geometric realization is given by 
\[
|\ul{X}| := \frac{\coprod_n \ul{X}[n] \times |\Delta^n|}{(d_i x, u) \sim (x, \delta_i u), (s_i x, u) \sim (x, \sigma_i u)}.
\]
The topological space $|\ul{X}|$ is a CW complex \cite[\S 14]{May67}.

An alternative, concise way to express this is as follows: we declare $|\Delta^n|$ to be the standard $n$-simplex $\{ (x_0, \ldots, x_n) \in \R^{n+1} \co 0 \leq x_i \leq 1, \sum x_i = 1\}$ and then we define 
\[
|\ul{X}| = \colim_{\Delta^n \rightarrow \ul{X}} |\Delta^n|
\]
where the indexing category has as its objects the maps $\Delta^n \rightarrow \ul{X}$ for varying $n$, and as its morphisms the maps $\Delta^m \rightarrow \Delta^n$ respecting the given maps to $\ul{X}$. Note that the geometric realization of the simplicial set $\Delta^n$ agrees with the standard $n$-simplex $|\Delta^n|$, justifying the notation. 

\begin{prop}
Geometric realization is left adjoint to the singular simplices functor:
\[
\Hom_{\msf{Top}}(|\ul{X}|, Y) \cong \Hom_{\msf{sSet}}(\ul{X}, \Sing Y).
\]
\end{prop}

\begin{remark}\label{rem: quillen equivalence} These adjoint functors define a ``homotopy equivalence'' of categories in a suitable sense, which we are not yet equipped to define precisely. The technical statement is that the functors induce a Quillen equivalence of \emph{model categories}, with the usual Quillen model structures on each side. See \cite[\S 16]{May92} for the precise formulation and proof. This justifies the assertion that simplicial sets provide a combinatorial model of topological spaces. 
\end{remark}

\subsubsection{Homotopy groups of a simplicial set}\label{sssec: htpy groups}
A pointed simplicial set $(\ul{X}, \star)$ has \emph{homotopy groups} $\pi_i(\ul{X}, \star)$, which could be defined as the homotopy groups (in the usual sense) of the geometric realization $|\ul{X}|$ with respect to the basepoint $\star$. They could also be defined directly within the category of simplicial sets, basically by using a combinatorial model for the sphere as the simplicial set $\Delta^n$ modulo its boundary, but this involves the subtlety of ``resolving'' $\ul{X}$ by a \emph{Kan complex} (see Remark \ref{rem: quillen equivalence}). We remark that for a topological space $Y$, $\Sing Y$ is a Kan complex. The point is that in ``homotopy-theoretic'' categories (such as the category of simplicial sets), objects may not have ``enough'' maps in or out of them and so must be ``resolved'' by ones that do. Every topological space is \emph{fibrant}, which is why $|\ul{X}|$ does not have to be similarly ``resolved''. 

When we discuss simplicial (abelian) groups and simplicial commutative rings, we will always take $\star$ to be the point corresponding to the identity element, and omit the basepoint from the notation. 

\subsubsection{Internal Hom}\label{sssec: internal hom}

Given two simplicial sets $\ul{X}$ and $\ul{Y}$, we construct a simplicial set $\ul{\Map}(\ul{X}, \ul{Y})$ whose set of $0$-simplices  is identified with $\Hom_{\msf{sSet}}(\ul{X}, \ul{Y})$. Intuitively, we are trying to put a topological structure on the set of maps from $\ul{X}$ to $\ul{Y}$. This will make the category of simplicial sets into a \emph{simplicial category} -- that is, a category enriched over simplicial sets. 

We define 
\[
\ul{\Map}(\ul{X}, \ul{Y})([n]) := \Hom_{\msf{sSet}}(\ul{X} \times \Delta^n, \ul{Y}).
\]
Note that the product of functors is formed level-wise.

\begin{remark}
A simplicial category is a possible model for the notion of $\infty$-category (although not the main one used in \cite{Lur09}). Intuitively, this is a type of category enriched over \emph{spaces}, i.e., where the sets of morphisms are promoted to spaces. 
\end{remark}

\subsection{Simplicial widgets}\label{ssec: simplicial widget}

More generally, if $\msf{C}$ is any category, then we say that a \emph{simplicial object of $\msf{C}$} is a functor $\Delta^{\op} \rightarrow \msf{C}$. If $\msf{C}$ is the category of ``widgets'', then the functor category $\Hom(\Delta^{\mrm{op}}, \msf{C})$ will be called the category of ``simplicial widgets'', and denoted $\msf{sC}$. 

\begin{example}
The following examples will be of particular interest to us. 
\begin{itemize}
\item A \emph{simplicial abelian group} is a functor from $\Delta^{\op}$ to the category of abelian groups.
\item Given a commutative ring $R$ (for us this always means commutative and with unit), a \emph{simplicial $R$-module} is a functor from $\Delta^{\op}$ to the category of $R$-modules. 
\item A \emph{simplicial commutative ring} is a functor from $\Delta^{\op}$ to the category of commutative rings. 
\item Given a commutative ring $R$, a \emph{simplicial $R$-algebra} is a functor from $\Delta^{\op}$ to the category of $R$-algebras. 
\end{itemize}
\end{example}

The categories of such objects is, in each case, enriched over $\sSet$. To see this, we first observe that for a simplicial set $\ul{S}$ and a simplicial widget $\ul{X}$, we have a simplicial widget 
\[
(\ul{X}^{\otimes \ul{S}})[n] = \coprod_{S_n} X_n
\]
with the coproduct on the RHS formed in the category of widgets. Then we define
\[
\ul{\Hom}_{\msf{sC}}(\ul{X}, \ul{Y})[n] := \Hom_{\msf{sC}}( \ul{X}^{\otimes \Delta^n }, \ul{Y}).
\]

\begin{example}\label{ex: free sab and scr}
Let $\ul{X}$ be a simplicial set. Then there is a simplicial abelian group $\Z\langle \ul{X} \rangle$, obtained by forming the free abelian group level-wise: define $\Z\langle \ul{X} \rangle ([n]) := \Z\langle\ul{X}([n]) \rangle$ and letting the face and degeneracy maps be induced by those of $\ul{X}$. 

There is also a simplicial $\Z$-algebra  $\Z[\ul{X}]$ obtained by forming the free (polynomial) $\Z$-algebra level-wise, with the face and degeneracy maps induced by those of $\ul{X}$. 
\end{example}

\begin{remark}\label{rem: kan complex}
It turns out that simplicial groups are \emph{automatically} Kan complexes. Therefore, the homotopy groups of simplicial groups can be calculated as maps from simplicial spheres, in contrast to the situation for general simplicial sets as cautioned in \S \ref{sssec: htpy groups}. The same applies for simplicial abelian groups, simplicial commutative rings, etc. because of how the model category structures are defined in each of these cases. 
\end{remark}

\subsection{The Dold-Kan correspondence}

Recall that if $\msf{C}$ is an abelian category, then we can form the category of chain complexes of objects in $\msf{C}$. This is essential for example in homological algebra. It turns out that this is closely related to the category $\msf{sC}$.

\begin{thm}[Dold-Kan]\label{thm: Dold-Kan}
Let $R$ be a commutative ring. The category of simplicial $R$-modules is equivalent to the category of non-negatively graded chain complexes of $R$-modules. 
\end{thm}

We may view Theorem \ref{thm: Dold-Kan} as explaining that ``simplicial'' is a generalization of ``chain complex'' to non-abelian categories (where the notion ``chain complex'' does not exist).

We will indicate the functors used to define the equivalence of Theorem \ref{thm: Dold-Kan}, starting with the functor from simplicial $R$-modules to chain complexes of $R$-modules. A proof may be found in \cite[\S III.2]{GJ09}; we also like the treatment in \cite{MDK}, which we follow. Let $\ul{M}$ be a simplicial $R$-module. First we explain an auxiliary construction called the \emph{Moore complex} of $\ul{M}$. From $\ul{M}$ we form a chain complex $M_{*}$ by taking $M_n:= \ul{M}[n]$, and taking the differential $\partial_n \co M_n \rightarrow M_{n-1}$ to be the alternating sum of the face maps, 
\[
\partial_n := \sum_{i=0}^n (-1)^i d_i.
\]

\begin{example}
The singular chain complex of a topological space $X$ is by definition the Moore complex of the simplicial abelian group $\Z\langle\Sing X\rangle$. 
\end{example}

\begin{remark}The homology groups of the Moore complex $M_*$ coincide with the homotopy groups of the simplicial abelian group $\ul{M}$. 
\end{remark}

Let us analyze the structure on the Moore complex. There are degeneracy maps $s_i \co M_n \rightarrow M_{n+1}$ for $0 \leq i \leq n$. Topological intuition suggests that ``removing'' the degenerate simplices should not affect the homology. To implement this, define $DM_{*} \subset M_{*}$ so that $DM_{n+1}$ is the span of the image of the degeneracy maps $s_0, \ldots, s_n$. 

\begin{exercise} Check that $DM_{*}$ is a subcomplex of $M_*$, i.e. is preserved by the differential $\partial$. Furthermore, show that the map $M_* \rightarrow (M/DM)_*$ is a quasi-isomorphism. 
\end{exercise}

Basically, we want to instead consider the chain complex $(M/DM)_*$. However, it is convenient to use a different normalization of this, which we will call the \emph{normalized Moore complex}. Define $NM_n$ to be the kernel of \emph{all} the face maps $d_i$ for $i<n$ (but not $i=n$). Then $(-1)^n d_n$ defines a differential $NM_n \rightarrow NM_{n-1}$. 

\begin{exercise}Check that $NM_*$ is a chain complex. 
\end{exercise}

\begin{remark}The sum map 
\[
NM_n \oplus DM_n \rightarrow M_n
\]
is an isomorphism. Therefore, $NM_*$ maps isomorphically to $(M/DM)_*$. 
\end{remark}

The functor from simplicial $R$-modules to chain complexes of $R$-modules, that we will take in Theorem \ref{thm: Dold-Kan}, is $\ul{M} \mapsto NM_*$.  A key step in the proof of the Dold-Kan equivalence is to show that 
\[
\ul{M}[n] 	\cong \bigoplus_{[n] \surj [k]} NM_k
\]
functorially in $\ul{M}$. Here the index set is over maps $[n] \rightarrow [k]$ in the simplex category (i.e., order-preserving maps) which are surjective. This tells us how to define the inverse functor: given a chain complex $M_*$, we will define a simplicial $R$-module by 
\[
\ul{M}[n] := \bigoplus_{[n] \surj [k]} M_k.
\]
The details in checking that this defines an equivalence of categories are left to the references.

\subsection{Simplicial commutative rings}

A \emph{simplicial commutative ring} is a functor $\cR$ from $\Delta^{\op}$ to the category of commutative rings. Because we shall work with these a lot, we adopt a slightly more economical notation. Simplicial commutative rings will be denoted using calligraphic letters, and we abbreviate $\cR_n := \cR[n]$. Classical commutative rings will be denoted using Roman letters such as $R$. 

So, a simplicial commutative ring can be specified concretely by a collection of commutative rings $\cR_n$, for $n \geq 0$, with maps between them indexed by the morphisms in $\Delta^{\op}$, which can be specified by face maps $d_i \co \cR_{n+1} \rightarrow \cR_n$
and degeneracy maps $s_i \co \cR_n \rightarrow \cR_{n+1}$ satisfying the simplicial identities. This is an axiomatization of the structure that exists on the singular simplices of a topological commutative ring. The category of commutative rings is denoted $\msf{CR}$ and the category of simplicial commutative rings is denoted $\msf{SCR}$.

\begin{example}\label{2ex: constant SCR 1}
Let $R$ be a commutative ring. Then we may define a simplicial commutative ring $\ul{R}$ such that $\ul{R}[n] = R$, and all face and degeneracy maps are the identity map. Intuitively, $\ul{R}$ corresponds to the topological commutative ring which is $R$ equipped with the discrete topology. 
\end{example}

\subsubsection{Homotopy groups}\label{sssec: homotopy groups}
If $\Cal{R}$ is a simplicial commutative ring, then its homotopy groups $\pi_*(\Cal{R})$ form a graded-commutative ring. The underlying group of $\pi_*(\cR)$ coincides with that of the underlying simplicial abelian group of $\cR$, and can therefore be computed by the Moore complex 
\[
\ldots \rightarrow \cR_n \xrightarrow{\partial_n} \cR_{n-1} \xrightarrow{\partial_{n-1}} \ldots \xrightarrow{\partial_0} \cR_0.
\]

\begin{example}
For the constant simplicial commutative ring $\ul{R}$, the Moore complex reads
\[
\ldots  \xrightarrow{\Id} R  \xrightarrow{0} R \xrightarrow{\Id} R \xrightarrow{0} R 
\]
Hence we see $\pi_0(\ul{R}) = R$ and $\pi_i(\ul{R}) = 0$ for $i >0$, in accordance with the intuition expressed in Example \ref{2ex: constant SCR 1}.
\end{example}

\begin{remark}
It turns out that the augmentation ideal $\pi_{>0}(\cR)$ always has a \emph{divided power structure} in $\pi_*(\cR)$. We will not explain what this means or why it exists, except that for a $\Q$-algebra a divided power structure always exists and is unique, whereas in positive characteristic it is a rather special piece of structure. 
\end{remark}

\subsubsection{Classical truncation} Recall that in \S \ref{ssec: visualizing derived schemes} we asserted the existence of an adjunction 
\begin{equation}\label{1eq: CR SCR}
\begin{tikzcd}[column sep = huge]
\msf{CR} \ar[r, bend right, "R \mapsto \ul{R}"']  & \msf{SCR} \ar[l, bend right, "\pi_0(\cR) \mapsfrom \cR"'] 
\end{tikzcd}
\end{equation}
To justify this, we need to argue that 
\[
\Hom_{\msf{SCR}}(\cR, \ul{R}) \cong \Hom_{\msf{CR}}(\pi_0(\cR), R).
\]
A morphism $f \in \Hom_{\msf{SCR}}(\cR, \ul{R})$ amounts to a system of maps $f_n \co \cR_n \rightarrow R$ compatible with the face and degeneracy maps. Since all the face and degeneracy maps on $\ul{R}$ are the identity, it must therefore be the case that $f_n(r)$ can be calculated by composing with any of the $n+1$ distinct maps $\cR_n \rightarrow \cR_0$ and then applying $f_0 \co \cR_0 \rightarrow R$. For this to be well-defined, we see that $f_0$ must vanish on the image of $d_1 - d_0$, and thus factor through $\pi_0(\cR) \rightarrow R$. It remains to show that any such $f_0$ does induce a well-defined map of simplicial commutative rings, which is completed by the following exercise.

\begin{exercise}
Show that all of the $n+1$ maps $\cR_n \rightarrow \cR_0$ induced by the $n+1$ maps $[0] \rightarrow [n]$ have the same composition with the quotient $\cR_0 \surj \pi_0(\cR)$.  
\end{exercise}

\subsubsection{Enrichment over simplicial sets}\label{sssec: internal hom SCR}

Because of its importance, let us explicate the enrichment of simplicial commutative rings over simplicial sets, although it is a special case of the description in \S \ref{sssec: internal hom}. Given a set $S$ and a commutative ring $R$, we can form the commutative ring $R^{\otimes S}$. Then for a simplicial set $\ul{S}$ and a simplicial commutative $\cR$, we define $\cR^{\otimes \ul{S}}$ by $(\cR^{\otimes \ul{S}})_n:=  (\cR_n)^{\otimes \ul{S}[n]}$. Finally, we define a simplicial set $\ul{\Hom}_{\msf{SCR}}(\cR, \cR')$ with 
\[
\ul{\Hom}_{\msf{SCR}}(\cR, \cR')[n] := \Hom_{\msf{SCR}}( \cR^{\otimes \Delta^n}, \cR')
\]
and the natural face and degeneracy maps.

\begin{exercise}
For commutative rings $R$ and $R'$, calculate the homotopy groups of $\ul{\Hom}_{\msf{SCR}}(\ul{R}, \ul{R}')$. (Answer: $\pi_0 = \Hom_{\msf{CR}}(R, R')$, and $\pi_i$ vanish for $i>0$.) This justifies thinking of $\msf{CR}$ as being embedded fully faithfully, in the sense of \emph{simplicial categories}, in $\msf{SCR}$. 
\end{exercise}

\subsection{Simplicial resolutions} We now discuss a key construction, in which a morphism of simplicial commutative rings is ``resolved'' by a process analogous to formation of projective or injective resolutions.

\subsubsection{Derived functors}

Let us recall the paradigm of derived functors in homological algebra, which is probably familiar to you. If $F$ is a right-exact functor on an abelian category $\msf{C}$, then we define the ``higher derived functors'' of $F$ on $M \in \msf{C}$ by finding a ``free'' resolution 
\[
\ldots \rightarrow P_n \rightarrow P_{n-1} \rightarrow \ldots \rightarrow P_0 \rightarrow M \rightarrow 0
\]
and then applying $F$ instead to the complex $\ldots \rightarrow P_n \rightarrow P_{n-1} \rightarrow \ldots \rightarrow P_0$, which we view as a ``replacement'' for $M$. We are being deliberately vague about what the adjective ``free'' should mean. 

\begin{example}
Let $\msf{C}$ be the category of $R$-modules and $N$ an $R$-module, $F = N \otimes_R (-)$. The higher derived functors of $F$ are $\Tor^R_i(N, -)$. 
\end{example}

Similarly, if $F$ is a left-exact functor, then one performs an analogous process using \emph{injective} resolutions instead. 

\begin{example}
Let $\msf{C}$ be the category of $R$-modules and $N$ an $R$-module, $F = \Hom_R(N, -)$. The higher derived functors of $F$ are $\Ext_R^i(N, -)$. 
\end{example}

Simplicial commutative rings will play the role for commutative rings that chain complexes of $R$-modules play for $R$-modules. In particular, we will see that replacing a commutative ring by a type of ``free resolution'' will allow us to build a theory of ``derived functors'' on the category of simplicial commutative rings. We will illustrate this in the particular example of the cotangent complex. 

Recall that for a morphism $f \co A \rightarrow B$, we have the $B$-module of K\"{a}hler differentials $\Omega^1_f$, which is ``well-behaved'' if $f$ is smooth. For general $f$, we will ``resolve'' $B$ by a ``smooth'' simplicial $A$-algebra in order to define the cotangent complex $\bL_f$. For abelian categories, this meant replacing $B$ by a complex of objects which level-wise had good properties (e.g., ``free''), but that does not apply to the category of $A$-algebras, which is far from abelian. Instead, this theory of homological algebra will be replaced by a theory of \emph{homotopical algebra} developed by Quillen. The correct framework for this is that of Quillen's \emph{model categories}, but this would involve a significant digression to explain completely, so we will adopt an ad hoc presentation and then hint at the general theory later, in \S \ref{ssec: model cats}.

\begin{defn}\label{defn: weak equivalence}
A morphism of simplicial commutative rings $f \co \cR \rightarrow \cR'$ is a \emph{weak equivalence} if it induces an isomorphism of homotopy groups $\pi_*(f) \co \pi_*(\cR) \xrightarrow{\sim} \pi_* (\cR')$. 
\end{defn}

This notion is an analogue for simplicial commutative rings of the notion of \emph{quasi-isomorphism} for chain complexes. Indeed, we could make an analogous definition for simplicial $R$-modules, which would transport under the Dold-Kan Theorem to the usual notion of quasi-isomorphism.

\subsubsection{Free simplicial algebras}

Let $\cA$ be a simplicial commutative ring. A \emph{simplicial $\cA$-algebra} is a simplicial commutative ring $\cR$ equipped with a homomorphism of $\cA \rightarrow \cR$ of simplicial commutative rings. 

\begin{defn}
A \emph{free simplicial $\cA$-algebra} $\cR$ is a simplicial $\cA$-algebra of the following form: 
\begin{itemize}
\item There is a system of sets $X_n$ such that $\cR_n = \cA[X_n]$. 
\item The degeneracy maps send $s_j(X_n) \subset X_{n+1}$. (Note that there is \emph{no} condition on face maps!) 
\end{itemize}
\end{defn}

The following Lemma illustrates that a free simplicial $\cA$-algebra has good ``mapping out'' properties, similar to those of free resolutions of modules. 

\begin{lemma}\label{lem: lifting resolutions}
Consider a commutative diagram of simplicial algebras 
\[
\begin{tikzcd}
\cA \ar[r] \ar[d] & \cB \ar[d, "\sim"', "\Phi"] \\
\cR \ar[r] \ar[ur, dashed] & \cC
\end{tikzcd}
\]
where $\Phi$ is level-wise surjective and a weak equivalence, and $\cR$ is a free simplicial $\cA$-algebra. Then there exists a lift $\cR \rightarrow \cB$ making the diagram commutative. 
\end{lemma}

\begin{proof}
See \cite[(5.4)]{Goe07}. For the more general statements in the framework of model categories, see \cite[\S 4]{DS95}.
\end{proof}

\begin{prop}[Existence of resolutions]\label{prop: resolutions exist}
Let $f \co A \rightarrow B$ be a map of (classical) commutative rings. Then there exists a free simplicial $A$-algebra $\cB$ and a diagram 
\[
\begin{tikzcd} 
 & \cB  \ar[d, "\sim"] \\ 
 \ul{A} \ar[r, "f"'] \ar[ur] & \ul{B} 
 \end{tikzcd}
\]
We remind that the symbol $\xrightarrow{\sim}$ denotes a weak equivalence, and $\ul{A}$ (resp. $\ul{B}$) denotes the constant simplicial commutative ring on $A$ (resp. $B$). Furthermore, we can even arrange that the formation of this diagram is functorial in $B$. 
\end{prop}

\begin{proof}[Proof sketch] We will give a ``canonical'' resolution, which is functorial in $B$. Note that for any $A$-algebra $R$, we have a canonical algebra homomorphism
\begin{equation}\label{3eq: underlying}
A[R] \rightarrow R
\end{equation}
sending $[r]$ to $r$. This construction gives rise to two maps $A[A[B]] \rightarrow A[B]$:
\begin{itemize}
\item Apply \eqref{3eq: underlying} with $R = B$ to get a map $A[B] \rightarrow B$, and then take $A[-]$. For example, this map sends $[a[b]] \mapsto [ab]$. 
\item Apply \eqref{3eq: underlying} with $R = A[B]$. For example, this maps sends $[a[b]] \mapsto a[b]$. 
\end{itemize}
We also have a map $A[B] \rightarrow A[A[B]]$ sending $[b]$ to $[[b]]$, which is a section of both maps above. 

Contemplating the combinatorics of the situation, we get a simplicial object
\begin{equation}\label{eq: canonical resolution}
\ldots  A[A[B]] \rightrightarrows A[B] 
\end{equation}
where the maps defined previously are the first face and degeneracy maps. Clearly this resolution is level-wise a polynomial algebra over $A$, and one can check that it is free. Why is it a resolution of $B$? 

This comes from an ``extra degeneracy'' argument, which is part of a more general pattern. If we have a simplicial set $\ul{X}$ 
\[
\ldots X_1 \rightrightarrows X_0
\]
with an augmentation $d_0 \co X_0 \rightarrow X_{-1}$ and ``extra degeneracies'' $s_{-1} \co X_{n-1} \rightarrow X_n$ satisfying the natural extensions of the simplicial identities, then the map $\ul{X} \rightarrow X_{-1}$ (where $X_{-1}$ is regarded as a constant simplicial set) is a weak equivalence. 

Going back to \eqref{eq: canonical resolution}, there are extra degeneracies $B \rightarrow A[B]$ sending $b \mapsto [b]$ and $s_{-1} \co A[B] \mapsto A[A[B]]$ sending $a[b] \mapsto [a[b]]$, etc. This can be used to show that the map from \eqref{eq: canonical resolution} to $B$ is a weak equivalence. 

\end{proof}

\begin{remark}\label{rem: monad}
The construction of the proof is a special case of a more general one, and we give the general statement because we feel that it elucidates the situation. Let $T$ be a \emph{monad} on a category $\msf{C}$, meaning the data of an endofunctor $T \co \msf{C} \rightarrow \msf{C}$ plus natural transformations $\eta \co \Id_{\msf{C}} \rightarrow T$ and $\mu \co T \circ T \rightarrow T$ satisfying various coherence conditions detailed \cite[\S VI]{Mac71}. We remark that monads often arise in practice as compositions of adjoint functors. (The running example of interest is: $\msf{C}$ is the category of sets, and $T$ is the functor taking a set $S$ to the underlying set of $A[S]$. Here $T$ can be viewed as arising from a free-forgetful adjunction.) 

An \emph{algebra} over the monad $T$ consists of the data of $x \in \msf{C}$ plus a morphism $a \co T x \rightarrow x$ satisfying coherence conditions. (In the running example, the category of algebras over $T$ is equivalent to the category of $A$-algebras.) Given an algebra $x$ over $T$, one can form a simplicial set called the \emph{bar construction} $B(T,x)_{\bu}$, which has $B(T,x)_n = T^{\circ (n+1)} x$. By formal combinatorial analysis, the bar construction is equipped with extra degeneracies. The augmentation $B(T,x)_{\bu} \rightarrow x$ will therefore be a weak equivalence in great generality. 

\end{remark}

\subsubsection{Uniqueness} If $f \co \cA \rightarrow \cB$ is a morphism of simplicial commutative rings, we will refer to a diagram 
\[
\begin{tikzcd} 
 & \cB'  \ar[d, "\sim"] \\ 
 \cA \ar[r, "f"'] \ar[ur] & \cB 
 \end{tikzcd}
\]
where $\cB'$ is a free simplicial $\cA$-algebra, as a \emph{free resolution} of $\cB$ as a simplicial $\cA$-algebra. Let us discuss ``how unique'' these resolutions are. In the classical situation of homological algebra, projective resolutions are unique up to homotopy. Analogously, there is a sense in which free resolutions of a simplicial commutative $\cA$-algebra are ``unique up to homotopy''. The best language for formulating such statements is model category theory, which would take us too far afield to explain, but we can say what it specializes to in this particular situation. 

We write $\cA[X] := \cA \otimes_{\ul{\Z}} \ul{\Z[X]}$ for the one-variable polynomial algebra over $\cA$, the tensor products being formed levelwise (as will be the case for all tensor products formed in this paragraph; they will be the ``same'' as the derived tensor products to be discussed in \S \ref{sssec: derived global deformation functor}). Let us call $\cA[X_0,X_1] := \cA[X] \otimes_{\cA} \cA[X]$. This represents the ``homotopy coproduct''; in particular, given two $\cA$-algebra homomorphisms $f,g \co \cA[X] \rightarrow \cB$ we get an $\cA$-algebra homomorphism $f \otimes g \co \cA[X_0,X_1] \rightarrow \cB$. There is a product morphism $\mu \co \cA[X_0, X_1] \rightarrow \cA[X]$. Let $\cA[X_0,X_1,Y]$ be a simplicial resolution of $\cA[X_0,X_1] \xrightarrow{\mu} \cA[X]$. Then we say that $f$ and $g$ are \emph{homotopic} if $f \otimes g$ extends to a commutative diagram 
\[
\begin{tikzcd}
\cA[X_0,X_1] \ar[r] \ar[dr, "f \otimes g"']  & \cA[X_0,X_1,Y] \ar[d, dashed] \\
& \cB
\end{tikzcd}
\]
Such a diagram is called a \emph{homotopy} between $f$ and $g$. 

We say that two simplicial resolutions 
\[
\begin{tikzcd} 
 & \cB'  \ar[d, "\sim"] \\ 
 \cA \ar[r, "f"'] \ar[ur] & \cB 
 \end{tikzcd}\quad \begin{tikzcd} 
 & \cB''  \ar[d, "\sim"] \\ 
 \cA \ar[r, "f"'] \ar[ur] & \cB 
 \end{tikzcd} 
\]
are \emph{homotopic} if there are maps $\alpha \co \cB' \rightarrow \cB''$ and $\beta \co \cB'' \rightarrow \cB'$ making the appropriate diagrams commute, and such that $\alpha \circ \beta$ and $\beta \circ \alpha$ are homotopic to the identity. Lemma \ref{lem: lifting resolutions} may be used to show that any two simplicial resolutions are homotopic.

\subsection{A vista of model categories}\label{ssec: model cats}

Let us hint at the more general framework underlying these procedures, which is Quillen's theory of \emph{model categories}. Our discussion will be very brief; a favorite introductory reference is \cite{DS95}.

The purpose of a model category is to incorporate homotopy theory into category theory. The basic idea is that in practice, one often wants to introduce a notion of maps that are ``weak equivalences'', but not necessarily isomorphisms. For example, in the category of topological spaces there is the usual notion of homotopy equivalences, and in the category of chain complexes there is the notion of quasi-isomorphism. A model category is a framework to do homotopy-theoretic constructions in such a situation. 

More formally, a model category is a category equipped with several distinguished classes of morphisms, satisfying various properties. The most important are the \emph{weak equivalences}, which designate maps that are ``homotopy equivalences''. The weak equivalences are morphisms that we want to ``invert''. Familiar examples are:
\begin{itemize}
\item In the category of chain complexes of $R$-modules, weak equivalences should be the quasi-isomorphisms. Inverting them leads to the \emph{derived category of $R$-modules}. 
\item In the category of simplicial commutative rings, weak equivalences should be as defined in Definition \ref{defn: weak equivalence}.
\end{itemize} 

The other data in a model category are classes of morphisms called the \emph{cofibrations} and \emph{fibrations}, which are roughly the generalizations of projective and injective resolutions. There are various axioms placed on these classes of morphisms: for example, they are preserved by compositions and retracts. An object is called \emph{fibrant} if its map to the terminal object is a fibration, and an object is called \emph{cofibrant} if its map from the initial object is a cofibration. 

\begin{example}
In the standard model structure on simplicial sets, the \emph{Kan complexes} referred to in \S \ref{sssec: htpy groups} are the fibrant objects. Kan complexes are characterized by lifting conditions possessed by the singular simplicial set of a topological space; in particular, $\Sing Y$ is a Kan complex for any  topological space $Y$. This is the reason why, when defining homotopy groups of simplicial sets as maps from a sphere, it was necessary to replace the target by a Kan complex. 

On the other hand, in the standard model category structure on topological space the fibrations are Serre fibrations, and in particular all topological spaces are fibrant. This is the reason why, when defining homotopy groups of a topological space, it was \emph{not} necessary to replace the target.

\end{example}

The axioms of a model category imply that every object admits a weak equivalence from a cofibrant object, and a weak equivalence to a fibrant object. We think of such a map as a ``resolution'' by a cofibrant object or a fibrant object. As we know, derived functors on abelian categories are defined using such resolutions. A model category structure provides a more general notion of ``resolution'' in categories which are not necessarily abelian, and therefore allows to construct ``non-abelian derived functors''.

Our particular constructions with simplicial commutative rings are examples of general constructions with model categories. In particular, our notion of ``free simplicial $\cA$-algebra $\cR$'' means that $\cR$ is a cofibrant $\cA$-algebra in the standard model structure; meanwhile, all simplicial commutative rings are fibrant. 

\subsection{The cotangent complex}

We may now define the \emph{cotangent complex} of a morphism $f \co A \rightarrow B$, as invented by Quillen and developed by Illusie in \cite{Ill71, Ill72}. Let $\cB \xrightarrow{\sim} B$ be a simplicial resolution of $B$ as a free simplicial $A$-algebra. Then let $\bL_{\cB/A}$ be the simplicial $\cB$-module obtained by forming K\"{a}hler differentials level-wise: $\bL_{\cB/A}[n] := \Omega^1_{\cB_n/A}$. Finally we define the \emph{cotangent complex} to be 
\[
\bL_f := \bL_{B/A}  := \bL_{\cB/A} \otimes_{\cB} B.
\]

\begin{exercise}
Show that if $\cB$ and $\cB'$ are two free simplicial resolutions of $B$, then $\bL_{\cB/A} \otimes_{\cB'} B$ is homotopic to $\bL_{\cB'/A} \otimes_{\cB} B$. Therefore, $\bL_{B/A}$ is well-defined up to homotopy. In particular, it gives a well-defined object of the derived category of $B$-modules. 
\end{exercise}

\begin{remark}[Derived generalizations] Even better, if $\cA$ is an animated ring, and $\cB$ is an animated $\cA$-algebra, then $\bL_{\cB/\cA}$ can be constructed as an animated $\cB$-module -- see \cite[\S 3.2]{Lur04}.
\end{remark}

We expect $\bL_{B/A}$ to have reasonable finiteness properties when the morphism $f$ has good finiteness properties. However, the canonical resolutions used to build $\cB$ in the proof of Proposition \ref{prop: resolutions exist} are extremely large. Therefore, we would like to know that there exist ``smaller'' resolutions. We will examine this question next. 

\subsection{Economical resolutions} We will show the existence of resolutions with ``good'' finiteness properties, following \cite[\S 4]{Iy07}.

\begin{prop}\label{prop: economical resolution}
Let $A$ be a noetherian commutative ring and $f \co A \rightarrow B$ a finite type morphism. Then there exists a simplicial $A$-algebra resolution of $B$ by a free simplicial $A$-algebra $A[X]$ with each $X_n$ a finite set. 
\end{prop}

\begin{cor}\label{cor: affine pseudocoherence}
Let $A$ be a noetherian commutative ring and $f \co A \rightarrow B$ a finite type morphism. Then there exists a representative of $\bL_{B/A}$ by a complex of finite free $B$-modules. 
\end{cor}

\begin{remark}
Typically $\bL_{B/A}$ cannot be represented by a perfect complex; that is, it will typically have homology groups in infinitely many degrees. For example, if $A$ is a field and $\bL_{B/A}$ has finite Tor-amplitude, then in fact $\bL_{B/A}$ has Tor-amplitude at most $2$, and $B$ is LCI over $A$. This is a result of Avramov \cite{Avr99}, resolving a Conjecture of Quillen.  
\end{remark}

As motivation, we recall how to construct ``efficient'' resolutions of a finite type $A$-module $M$ by a complex of free $A$-modules. We can build a sequence of complexes of finite free $A$-modules $\{F_i\}$ whose homology approximates $M$ in degrees up to $i$. For $F_0$, we pick any surjection from a free module $A^n \surj M$. Then we pick generators for $I := \ker (F_0 \surj M)$, which induces a map $A^m \xrightarrow{\partial} A^n$ with image $I$, whose map to $M$ induces an isomorphism on $H_0$. So we may take $ F_0 = [A^m \xrightarrow{\partial} A^n]$. We then inductively build $F_{i}$ from $F_{i-1}$ by picking representatives in $F_{i-1}$ for generators of $H_{i}(F_{i-1})$, and then adding free summands in degree $i+1$ that bound these generators. 

To perform an analogous construction for simplicial commutative rings, what we want is a way to ``kill cycles'' in a simplicial commutative ring. This is accomplished by the following Lemma. 

\begin{lemma}\label{lem: adjoin}
Given a simplicial commutative ring $\cA$, and a homotopy class $[z] \in \pi_i(\cA)$, there exists a free simplicial $\cA$-algebra $\cA'$ such that: 
\begin{itemize}
\item $\cA_n'$ is a free $\cA_n$-algebra on finitely many generators for each $n$, and $\cA_n' = \cA_n$ for $n \leq i$. 
\item The map $\cA \rightarrow \cA'$ induces an isomorphism on $\pi_n$ for $n<i$, and for $n=i$ an exact sequence 
\[
0 \rightarrow \pi_0(\cA) \cdot [z] \rightarrow \pi_i(\cA) \rightarrow \pi_i(\cA') \rightarrow 0 
\]
\end{itemize}
\end{lemma}

Indeed, supposing Lemma \ref{lem: adjoin} is true, we can make a resolution as in Proposition \ref{prop: economical resolution} by picking a surjection $\ul{A[x_1, \ldots, x_d]} \surj B$. Then, by repeatedly applying Lemma \ref{lem: adjoin}, we may build a sequence of free simplicial $A$-algebras with the ``correct'' finiteness properties and homotopy groups in all finite degrees. Taking their colimit gives the desired resolution. 

Now to prove Lemma \ref{lem: adjoin}, we want to ``adjoin'' a variable whose boundary is $[z]$. However, because of the simplicial identities, the process of ``adjoining'' variables is necessarily quite complicated. For this, we will imitate what happens to the singular simplices when attaching cells to a topological space $X$ to kill a homology class in degree $i-1$. If we attach a cell in degree $i$, then this creates additional degenerate simplices in all higher dimensions. For this reason, the construction is considerably more complicated. 

\begin{proof}[Proof of Lemma \ref{lem: adjoin}] 
We define 
\[
X_n := \{x_t \mid  t \co [n] \surj [i+1] \in \Delta\}
\]
and we take $\cA_n' = \cA_n[X_n]$, with faces and degeneracies defined as follows. 
\begin{itemize}
\item We have $s_j(x_t) = x_{t \circ \sigma_j}$. 
\item Note that $X_{i+1} = \{x_{\Id}\}$. We set $d_0(x_{\Id}) = z$ and $d_j(x_{\Id}) = 0$ for $j>0$. For $n>i+1$, we define $d_j(x_t) = x_{t \circ \delta_j}$ if $t \circ \delta_j$ is surjective, and otherwise it factors through a face map $d_{j'} \co [i] \rightarrow [i+1]$ so we define it to make the following diagram commute 
\[
\begin{tikzcd}
{[n]} \ar[r, "t"] & {[i+1]} \\ 
{[n-1]} \ar[u, "\delta_j"]  \ar[r] & {[i]} \ar[u, "\delta_{j'}"] 
\end{tikzcd}
\]
\end{itemize}
One can check that this $\cA_n'$ is a simplicial commutative ring satisfying the conclusions of Lemma \ref{lem: adjoin} (see \cite{MSCR} for more details). 
\end{proof}

\begin{remark}
A more conceptual way to phrase this construction is as follows (following \cite{MSCR}). For a simplicial $A$-algebra $\cA$, a class in $\pi_i(\cA)$ is represented by a map of simplicial sets $\partial \Delta^{i+1} \rightarrow \cA$ (here we are using Remark \ref{rem: kan complex}, that $\cA$ is fibrant). This induces a map of simplicial commutative rings $A[\partial \Delta^{i+1}] \rightarrow \cA$, where $A[\partial \Delta^{i+1}]$ is the free simplicial $A$-algebra on the simplicial set $\partial \Delta^{i+1}$. We may then form 
\[
\cA' := \cA \otimes_{A[\partial \Delta^{i+1}]} A[\Delta^{i+1}]
\]
and we claim that it satisfies the conclusions of Lemma \ref{lem: adjoin}. The first bullet point is satisfied by construction. To compute the effect on homotopy groups, we note that $\pi_*(A) \xrightarrow{\sim}\pi_*(A[\Delta^{i+1}])$ because $\Delta^{i+1}$ is contractible, and the lowest positive-degree homotopy group of $A[\partial \Delta^{i+1}]$ is in degree $i$ and maps to $[z]$ in $\pi_i(\cA)$ by construction. Then conclude using the Tor spectral sequence \cite[(5.2)]{Qu70}
\[
\Tor_p^{\pi_*(A[\partial \Delta^{i+1}])}(\pi_*(\cA), \pi_*(A[\Delta^{i+1}]))_q \implies \pi_{p+q} (\cA \otimes_{A[\partial \Delta^{i+1}]} A[\Delta^{i+1}]). 
\]
\end{remark}

\subsection{Derived tensor products}\label{ssec: derived tensor}

Let $\cA \rightarrow \cB$ and $\cA \rightarrow \cB'$ be two maps of simplicial commutative rings. Let $\cA \inj \cP \xrightarrow{\sim} \cB$ be a free resolution of $\cB$ as an $\cA$-algebra and $\cA \inj \cP' \xrightarrow{\sim} \cB'$ be a free resolution of $\cB'$ as an $\cA$-algebra. The \emph{derived tensor product} ``$\cB \dotimes_{\cA} \cB'$'' is represented by the simplicial commutative ring $\cP \otimes_{\cA} \cP'$, with 
\[
(\cP \otimes_{\cA} \cP')_n = \cP_n \otimes_{\cA_n} \cP'_n. 
\]
Although this is not unique since it depends on a choice of resolutions, it is unique up to weak equivalence. In fact, it suffices to resolve only one of the terms, e.g., the simplicial commutative $\cP \otimes_{\cA} \cB'$ also represents the derived tensor product. Indeed, there is an evident map 
\[
\cP \otimes_{\cA} \cP' \rightarrow \cP \otimes_{\cA} \cB'
\]
and this is a weak equivalence of free $\cB'$-algebras. 

If $A$ is classical and $\cP$ is a free simplicial $A$-algebra, then the underlying simplicial $A$-module of $\cP$ associates under the Dold-Kan correspondence to a complex of free $A$-modules. In particular, the Dold-Kan correspondence takes free simplicial resolutions to the familiar notion of free resolutions of chain complexes. In particular, we see that if $\cA \xrightarrow{\sim}  A$, $\cB \xrightarrow{\sim}  B$, and $\cB' \xrightarrow{\sim}  B'$ are all classical, then we have
\[
\pi_i ( \cB \dotimes_{\cA} \cB') \cong \Tor^A_i(B, B').
\]
This has the following consequence. 

\begin{cor}\label{cor: tor-independent}
Let $A \rightarrow B$ and $A \rightarrow B'$ be maps of classical commutative rings, such that $\Tor^i_A(B,B') = 0$ for all $i>0$. Then $B \dotimes_A B' \rightarrow \pi_0(B \dotimes_A B') \cong B \otimes_A B'$ is a weak equivalence. 
\end{cor}

\subsection{Properties of the cotangent complex}\label{ssec: properties of cotangent complex}

Now we establish some properties of the cotangent complex which are familiar for the module of K\"{a}hler differentials, at least in the smooth case. As a reference we recommend \cite{Qu68}. 

\begin{prop}\label{prop: tensor cotangent}
If $B,B'$ are $A$-algebras with $\Tor^i_A(B,B') = 0$ for $i>0$, then 
\begin{enumerate}
\item $\bL_{B \otimes_A B'/B'} \cong \bL_{B/A} \otimes_A B'$.
\item $\bL_{B \otimes_A B'/A} \cong (\bL_{B/A} \otimes_{A} B') \oplus (B\otimes_A \bL_{B'/A})$. 
\end{enumerate}
\end{prop}

\begin{remark} It is true in general that 
\begin{enumerate}
\item $\bL_{\cB \dotimes_{\cA} \cB'/\cB'} \cong \bL_{\cB/\cA} \otimes_{\cA} \cB'$.
\item $\bL_{\cB \dotimes_{\cA} \cB'/\cA} \cong (\bL_{\cB/\cA} \otimes_{\cA} \cB') \oplus (\cB\otimes_{\cA} \bL_{\cB'/\cA})$. 
\end{enumerate}
The Tor-vanishing assumption is used to ensure that $B \dotimes_A B' \xrightarrow{\sim} B \otimes_A B'$. 
\end{remark}

\begin{proof}
Both arguments are completely formal from the standard properties of formation of K\"{a}hler differentials with respect to tensor products, using the observation that: 
\begin{itemize}
\item If $\cP$ is a free simplicial $\cA$-algebra, and $\cB$ is any simplicial $\cA$-algebra, then $\cP \otimes_{\cA} \cB$ is a free simplicial $\cB$-algebra. 
\end{itemize}
\end{proof}

\begin{prop}\label{prop: cotangent triangle}
If $A \rightarrow B \rightarrow C$ is a composition of morphisms, then there is an exact triangle in the derived category of $C$-modules: 
\[
\bL_{B/A} \otimes_B C \rightarrow \bL_{C/A} \rightarrow \bL_{C/B}
\]
\end{prop}

\begin{proof}
Choose a free resolution $A \inj \cB \xrightarrow{\sim}   B$. This gives $C$ the structure of a $\cB$-algebra, so we may then choose a free resolution $\cB \inj \cC \xrightarrow{\sim} C$. This gives a diagram 
\[
\begin{tikzcd} 
& & \cC  \ar[ddr, "\sim"] \\
& \cB \ar[d, "\sim"]  \ar[ur, hook]\\
A \ar[ur, hook]  \ar[r] & B \ar[rr] & & C
\end{tikzcd}
\]
Then we form the tensor product 
\[
\begin{tikzcd} 
& & \cC  \ar[ddr, "\sim"] \ar[d]  \\
& \cB \ar[d, "\sim"]  \ar[ur, hook] & B \otimes_{\cB} \cC \ar[dr]  \\
A \ar[ur, hook]  \ar[r] & B \ar[rr] \ar[ur] & & C
\end{tikzcd}
\]
From the exact sequence of K\"ahler differentials for smooth morphisms \cite[Tag 02K4]{stacks-project}, applied level-wise, we obtain an exact triangle
\[
\bL_{\cB/A} \otimes_{\cB} \cC \rightarrow \bL_{\cC/A} \rightarrow \bL_{\cC/\cB}.
\]
Applying $- \otimes_{\cC} C$, we get an exact triangle. Let us verify that the terms are as claimed. 
\begin{itemize}
\item $\bL_{\cB/A} \otimes_{\cB} \cC \otimes_{\cC} C \cong \bL_{\cB/A} \otimes_{\cB} C \cong \bL_{B/A} \otimes_B C$. 
\item $\bL_{\cC/A} \otimes_{\cC} C \cong \bL_{C/A}$. 
\item We have $\bL_{\cC/\cB} \otimes_{\cB} B \cong \bL_{B \otimes_{\cB} \cC / B}$. Since $B \otimes_{\cB} \cC$ is a free resolution of $C$ as a $B$-algebra, tensoring with $C$ over $\cC$ gives $\bL_{C/B}$. 
\end{itemize}
\end{proof}

\begin{cor}[Cotangent complex of localizations]\label{cor: cotangent of localizations} \hfill
Let $A$ be a commutative ring. 
\begin{enumerate}
\item If $S$ is a multiplicative system in $A$, then $\bL_{S^{-1}A/A} \cong 0$. 

\item Let $A \rightarrow B$ be a ring homomorphism. If $S$ is a multiplicative system in $A$ and $T$ is a multiplicative system in $B$ containing the image of $S$, then 
\[
\bL_{T^{-1} B / S^{-1} A } \cong  \bL_{B/A} \otimes_B T^{-1} B.
\]
\end{enumerate}
\end{cor}

\begin{proof} (i) Multiplication induces $S^{-1} A \otimes_A S^{-1}A \xrightarrow{\sim} S^{-1} A$. Putting this into Proposition \ref{prop: tensor cotangent}(ii), we get 
\[
\bL_{S^{-1} A /A} \cong \bL_{S^{-1} A \otimes_A S^{-1}A /A } \cong \bL_{S^{-1} A /A} \oplus \bL_{S^{-1} A /A}
\]
via the diagonal map. This forces $\bL_{S^{-1} A /A}\cong 0$. 

(ii) By Proposition \ref{prop: tensor cotangent}(i) We have $\bL_{S^{-1} B / S^{-1}A } \cong \bL_{B/A} \otimes_B S^{-1} B$. Apply Proposition \ref{prop: cotangent triangle} to get an exact triangle 
\[
\bL_{S^{-1} B / S^{-1} A} \otimes_{S^{-1} B} T^{-1} B \rightarrow \bL_{T^{-1} B / S^{-1} A}  \rightarrow \bL_{T^{-1} B/ S^{-1} B}.
\]
The rightmost term vanishes according to (i), from which the result follows. 

\end{proof}

\subsection{Some computations of the cotangent complex} Now we will identify the cotangent complex in some special cases. Simplicial resolutions are almost always too unwieldy to compute with by hand, so we will instead need to make clever use of the formal properties explained above.

\begin{prop}\label{prop: smooth cotangent complex} Let $A$ be a noetherian commutative ring and $f \co A \rightarrow B$ a finite type morphism. Then 
\begin{enumerate}
\item $f$ is \'{e}tale if and only if $\bL_f \cong 0$. 
\item $f$ is smooth if and only if $\bL_f \cong \Omega_f^1$ and is finite projective. 
\end{enumerate}
\end{prop}

\begin{proof}[Proof of one implication] 
We will only prove the forward directions for now. The converses will be established after we have developed the connection between $\bL_f$ and deformation theory.

Suppose $f$ is \'{e}tale. Then the multiplication map $ B \otimes_A B \rightarrow B$ is a localization, so Corollary \ref{cor: cotangent of localizations} implies that $\bL_{ B/B \otimes_A B } \cong 0$. By Proposition \ref{prop: tensor cotangent} we have $\bL_{B \otimes_A B/B} \cong \bL_{B/A} \otimes_A B$. Then applying Proposition \ref{prop: cotangent triangle} to $B \xrightarrow{\Id_B \otimes 1_B} B \otimes_A B \rightarrow B$, we find that $\bL_{B/A} \cong 0$. 

Next suppose $f$ is smooth. The local structure theorem for smooth maps says that locally on $\Spec B$, $f$ can be factored as an \'{e}tale map over an affine space: 
\[
\begin{tikzcd}
\Spec B \ar[r] \ar[dr, "f"'] &  \A^n_A \ar[d] \\
& \Spec A
\end{tikzcd}
\]
At the level rings, this means that after some localization, we have a factorization $A \rightarrow P \rightarrow B$ where $P$ is a polynomial ring over $A$ and $P \rightarrow B$ is \'{e}tale. Then the morphism of constant simplicial $A$-algebras $A \rightarrow P$ is already free, so $\bL_{P/A} \cong \Omega_{P/A}$ is finite free. From Proposition \ref{prop: cotangent triangle} we have an exact triangle $\bL_{P/A} \otimes_P B \rightarrow \bL_{B/A}  \rightarrow \bL_{B/P}$ where $\bL_{P/A} \otimes_P B \cong \Omega_{B/A}$ and $\bL_{B/P} = 0$, which completes the proof.

\end{proof}

Recall that a \emph{regular sequence} in a commutative ring $A$ is a finite sequence of elements $r_1, \ldots, r_n \in A$ such that  each $r_i$ is a non-zerodivisor modulo $A/(r_1, \ldots, r_{i-1})$. 

\begin{prop}\label{prop: LCI cotangent complex}
If $B = A/I$ and $I$ is generated by a regular sequence, then 
\[
\bL_{B/A} \xrightarrow{\sim} I/I^2[1].
\]
\end{prop}

\begin{proof}
Let us first treat the case $A = \Z[x]$ and $B = \Z$, $I = (x)$. In this case we have a section $\Z \rightarrow \Z[x]$. Apply Proposition \ref{prop: cotangent triangle} to the sequence $\Z \rightarrow \Z[x] \rightarrow \Z$ to get an exact triangle 
\[
\bL_{\Z[x]/\Z} \otimes_{\Z[x]} \Z \rightarrow \bL_{\Z/\Z} \rightarrow \bL_{\Z/\Z[x]}.
\]
This shows that $\bL_{\Z/\Z[x]} \cong \bL_{\Z[x]/\Z}[1] \otimes_{\Z[x]} \Z \cong (x)/(x^2)$. 

Now we consider the more general situation of the Proposition. By induction, it suffices to handle the case where $I = (r)$ is principal. Note that a choice of $r$ induces a map $\Z[x] \rightarrow A$ sending $x \mapsto r$, and this map induces $A/(r) \cong \Z \otimes_{\Z[x]} A$. We claim that $r$ is a non-zerodivisor if and only if $\Tor^i_{\Z[x]}(\Z, A) = 0$ for all $i>0$, or in other words (by Corollary \ref{cor: tor-independent}) if and only if $\Z \dotimes_{\Z[x]} A \xrightarrow{\sim} A/(r)$. The claim is an elementary computation in homological algebra, using the resolution $\Z[x] \xrightarrow{x} \Z[x]$ of $\Z$ over $\Z[x]$. Then by Proposition \ref{prop: tensor cotangent}, we have 
\[
\bL_{B/A} \cong \bL_{\Z/\Z[x]} \otimes_{\Z[x]} A \cong (x)/(x^2) [1] \otimes_{\Z[x]} A \cong (r)/(r^2)[1].
\]
\end{proof}

\subsection{Globalization}
The theory we have discussed can be globalized to schemes, and then to stacks. 
For example, let $Y \rightarrow X$ be a morphism of schemes. Then we define $\bL_{Y/X}$ by finding a free replacement of $\cO_Y$ as a simplicial $\cO_X$-algebra, and then forming $\Omega$-levelwise, etc. Traditionally, it was important to develop the theory this way instead of trying to ``glue'' the cotangent complexes from affine open subspaces, because of inadequate categorical technology: on the one hand the cotangent complex is not unique at the level of chain complexes, and on the other hand the derived category is not suitable for gluing. 

However, in modern language one can indeed glue the cotangent complex affine-locally in the desired way, viewing the cotangent complex as an animated module and using that the categories of animated quasicoherent sheaves satisfy Zariski descent.

\subsection{Andr\'{e}-Quillen homology}

Let $A \rightarrow B$ be a map of commutative rings. We have defined the cotangent complex $\bL_{B/A}$ as a simplicial $B$-module up to homotopy. 

\begin{defn}
Let $M$ be a $B$-module. We define the \emph{Andr\'{e}-Quillen homology} groups 
\[
D_i(B/A; M) := H_i(\bL_{B/A} \otimes_B M)
\]
and the \emph{Andr\'{e}-Quillen cohomology} groups 
\[
D^i(B/A; M) := H^i(\Hom_B(\bL_{B/A}, M)).
\]
\end{defn}

As a consequence of the properties of the cotangent complex, we have various properties of Andr\'{e}-Quillen (co)homology: 
\begin{itemize}
\item Let $B$ be an $A$-algebra and $M' \rightarrow M \rightarrow M''$ be a short exact sequence of $B$-modules. Then there is a long exact sequence 
\[
D^i(B/A; M') \rightarrow D^i(B/A; M) \rightarrow D^i(B/A; M'') \rightarrow D^{i+1}(B/A; M') \rightarrow \ldots 
\]
\item Let $B,B'$ be $A$-algebras such that $\Tor^A_i(B,B') = 0$ for $i>0$. Then for any $B \otimes_A B'$-module $M$, we have
\[
D^i(B \otimes_A B'/B'; M) \cong D^i(B/A; M)
\]
and
\[
D^i(B \otimes_A B'/A; M) \cong D^i(B/A; M) \oplus D^i(B'/A; M). 
\]
\item If $A \rightarrow B \rightarrow C$ is a sequence of ring homomorphisms and $M$ is a $C$-module, then we get a long exact sequence
\[
D^i(B/A; M) \rightarrow D^i(C/A; M) \rightarrow D^i(C/B; M) \rightarrow D^{i+1}(B/A;M) \rightarrow \ldots 
\]
\end{itemize}

We will next apply Andr\'{e}-Quillen cohomology to study deformation theory.	

\subsubsection{Deformation theory setup} 

Let $f \co X \rightarrow S$ be a scheme and $S \inj S'$ a square-zero thickening with ideal sheaf $\cI$. We consider deformations of $X$ to $X' \rightarrow S'$, as in the diagram 
\[
\begin{tikzcd}
X \ar[d]  \ar[r, dashed] & X' \ar[d, dashed] \\ 
S \ar[r, hook] & S'
\end{tikzcd}
\]
Then this is a matter of constructing an extension 
\[
\begin{tikzcd}
0 \ar[r] &  \cJ  \ar[r, dashed] &  \cO_{X'} \ar[r, dashed] & \cO_X  \ar[r] & 0 \\
0 \ar[r]  & \cI \ar[r] \ar[u]   & \cO_{S'} \ar[r] \ar[u, dashed] & \cO_S \ar[r] \ar[u] & 0 
\end{tikzcd}
\]
where $\cJ$ may be regarded as a sheaf on $\cO_X$ because of the square-zero property, and as such it is isomorphic to $f^* \cI \cong \cI \otimes_{\cO_S} \cO_X$. 

\subsubsection{Square-zero algebra extensions}
Motivated by this, we consider the following problem. Let $f \co A \rightarrow B$ a homomorphism of (classical) commutative rings and let $M$ be a $B$-module. An \emph{extension of $B$ by $M$ as $A$-algebras} is an exact sequence 
\[
0 \rightarrow M \rightarrow E \rightarrow B \rightarrow 0
\]
which presents $E$ as a square-zero $A$-algebra extension of $B$ by $M$. We would like to understand the structure of all such extensions. It turns out that this is controlled by Andr\'e-Quillen homology.

\begin{prop}\label{prop: extensions and AQ cohomology}
Let $M$ be a $B$-module. Then $D^1(B/A; M)$ is the set of isomorphism classes of extensions of $B$ by $M$ as $A$-algebras. 
\end{prop}

\begin{example}
The zero element $0 \in D^1(B/A; M)$ corresponds to the trivial square-zero extension $B \oplus M$, with multiplication 
\[
(b_1, m_1)(b_2, m_2) = (b_1b_2, b_1 m_2 + b_2 m_1). 
\]
\end{example}

\begin{proof}
Let us give maps in both directions. First suppose we have an extension $M \rightarrow E \rightarrow B$. Pick a free $A$-algebra resolution $A \rightarrow \cB \xrightarrow{\sim} B$. In particular, since $\cB_0$ is polynomial over $A$, we have some lifting 
\[
\begin{tikzcd}
 & \cB_0 \ar[r] \ar[d, dashed, "\phi_0"] & B \ar[d, equals] \\
M \ar[r] & E \ar[r] & B
\end{tikzcd}
\]
Now, we have two maps $d_0, d_1 \co \cB_1  \rightrightarrows \cB_0$ and we know that the image of $(d_0-d_1)$ is an ideal in $\cB_0$, with $\cB_0 / \mrm{Im}(d_0-d_1) \cong \pi_0(\cB) \cong B$. So the difference 
\[
\phi_0 \circ (d_0 - d_1) \co \cB_1 \rightarrow M
\]
is a $B$-linear derivation from $\cB_1$ to $M$, hence corresponds to an element of $\Hom_B(\Omega_{\cB_1}, M)$. One checks that it is a cocycle, hence induces a class in $D^1(B/A; M)$. The different choices of $\phi_0$ modify this cocycle by a coboundary, so we get a well-defined map from isomorphism classes of extensions to $D^1(B/A; M)$.

This also suggests how to construct the inverse map. Given a class in $D^1(B/A; M)$ represented by a derivation $\delta \co \cB_1 \rightarrow M$, define $E$ to be the pushout of $B$-modules. 
\[
\begin{tikzcd}
\cB_1 \ar[r, "d_0-d_1"] \ar[d, "\delta"] & \cB_0 \ar[d, dashed] \\
M \ar[r] & E
\end{tikzcd}
\]
We equip $E$ with the multiplication induced by its structure as a quotient of $\cB_0 \oplus M$. For this to be well-defined one checks that the image of $B_1$ under $(d_1-d_0, \delta)$ is a square-zero ideal, which is straightforward.

\end{proof}

\begin{example}
Consider $\Z$-algebra extensions of $\F_p$ by $\F_p$. Since we are considering arbitrary extensions of commutative rings, these are controlled by $D^1(\F_p/\Z; \F_p)$. By Proposition \ref{prop: LCI cotangent complex}, the cotangent complex of $\Z \rightarrow \F_p$ is $\bL_{\F_p/\Z} \cong \F_p[1]$. So 
\[
D^1(\F_p/\Z; \F_p) \cong \Hom_{\F_p}(\F_p[1], \F_p[1]) \cong \F_p.
\]
The 0 class corresponds to the extension $\F_p[\epsilon]/\epsilon^2$. The non-zero classes all have underlying ring $\Z/p^2\Z$, with the maps to $\F_p$ being the natural projection composed with an automorphism of $\F_p$. 
\end{example}

\begin{example}
Let $A \rightarrow B$ a smooth map of commutative rings and let $M$ be any $B$-module. Then by Proposition \ref{prop: smooth cotangent complex} we have $\bL_{B/A} \cong \Omega_{B/A}^1$ is finite projective, so $D^1(B/A; M) =0$. This says that there is a unique square-zero $A$-algebra extension of $B$ by $M$, which is necessarily the split extension. 
\end{example}

\begin{proof}[Completion of the proof of Proposition \ref{prop: smooth cotangent complex}]
We want to show that if $f \co A \rightarrow B$ is a finite type morphism of Noetherian rings, and $\bL_{B/A}$ is a finite projective $B$-module concentrated in degree $0$, then $f$ is smooth. Thanks to the finiteness hypotheses, it suffices to show that $f$ is formally smooth, which means that for any square-zero extension $S \inj \wt{S}$, the diagram 
\[
\begin{tikzcd}
S \ar[d, hook] \ar[r] & \Spec B \ar[d] \\
\wt{S} \ar[r] \ar[ur, dashed] & \Spec A
\end{tikzcd}
\]
has a lift. First suppose that the map $S \rightarrow \Spec B$ is an isomorphism. Then $\wt{S}$ is the spectrum of a square-zero $A$-algebra extension of $B$, and we want to show that it splits. But the assumption on $\bL_{B/A}$ implies that $D^1(B/A; M) = 0$ for any $B$-module $M$. So the class of $[\wt{S}] \in D^1(B/A; M)$ vanishes, which means that there is a lift as desired. 

Now we can treat the general case. Suppose $S = \Spec R$, $\wt{S} = \Spec \wt{R}$, and let $I = \ker(\wt{R} \rightarrow R)$. Choose a surjection $P \surj R$ from a polynomial $A$-algebra $P$, say with kernel $J$. By choosing lifts of generators, we find a lift $P \rightarrow \wt{R}$, which carries $J$ to $I$, hence factors over $P/J^2$. Then $P/J^2$ is a square-zero $A$-algebra extension of $B$, so by the case handled in the previous paragraph there is a splitting $B \rightarrow P/J^2$, whose composition to $\wt{R}$ gives the desired lift. 
\[
\begin{tikzcd}
P \ar[r, "\sim"] & B  =  P/J\ar[r] \ar[d, bend left, dashed] & R = \wt{R}/I \\
P \ar[r] \ar[u, equals] & P/J^2  \ar[u, twoheadrightarrow] \ar[r] &  \wt{R} \ar[u, twoheadrightarrow]
\end{tikzcd}
\]

\end{proof}

\subsection{Deformation theory of schemes}\label{ssec: deformation theory 1}

Let $\wt{A} \rightarrow A$ be an extension by a square-zero ideal $I$. We consider the problem of finding flat deformations 
\[
\begin{tikzcd}
J \ar[r] & \wt{B} \ar[r] & B \\
I \ar[u] \ar[r] & \wt{A} \ar[u] \ar[r] &  A \ar[u]
\end{tikzcd}
\]
or geometrically, deformations of $\Spec B \rightarrow \Spec A$ to a flat family over $\Spec \wt{A}$. 
\[
\begin{tikzcd}
\Spec B \ar[d] \ar[r, hook]  & \Spec \wt{B} \ar[d] \\
\Spec A \ar[r, hook] & \Spec \wt{A}
\end{tikzcd}
\]

Since $I$ is square-zero, we can regard $I$ as an $A$-module. Similarly, we can regard $J$ as a $B$-module. 

\begin{lemma}
In the situation above, $\wt{B}$ is flat over $\wt{A}$ if and only if $I \otimes_A B \xrightarrow{\sim} J$. 
\end{lemma}

\begin{proof} 
If $\wt{B}$ is flat over $\wt{A}$, then applying $\wt{B} \otimes_{\wt{A}} (-)$ to the short exact sequence $I \rightarrow \wt{A} \rightarrow A$, and using that the $\wt{A}$-action on $I$ factors through $A$, shows that $I \otimes_A B \xrightarrow{\sim} J$.

For the other direction, see \cite[\href{https://stacks.math.columbia.edu/tag/00MD}{Tag 00MD}]{stacks-project}.

\end{proof}

Henceforth we may assume that $J = I \otimes_A B$. We know that $\wt{A}$-algebra extensions of $B$ by $J$ are classified by $D^1(B/\wt{A};  J)$.  The composition $\wt{A} \rightarrow A \rightarrow B$ gives a long exact sequence 
\[
\begin{tikzcd}[row sep = tiny]
\ar[r] & D^1(B/A; J) \ar[r] &  D^1(B/\wt{A};  J) \ar[r] &  D^1(A/\wt{A}; J) \\
\ar[r] & D^2(B/A; J)
\end{tikzcd}
\]
The equivalence class of the extension $\wt{A}$ may be viewed as an element $[\wt{A}] \in  D^1(A/\wt{A}; J) $, which comes from $ D^1(B/\wt{A};  J)$ if and only if its image in $D^2(B/A;J)$ vanishes. By the long exact sequence, its pre-image in $D^1(B/\wt{A};  J) $ forms a torsor (possibly empty) for $D^1(B/A; J) $. Hence we find: 

\begin{prop}\label{prop: deformation-obstruction} Let the situation be as above. 
\begin{itemize}
\item Attached to the extension $\wt{A} \rightarrow A$ is an \emph{obstruction class} $\mrm{obs}(\wt{A}) \in D^2(B/A; J)$, which vanishes if and only if an extension $\wt{B}$ exists. 
\item If $\mrm{obs}(\wt{A}) \in D^2(B/A; J)$ vanishes, then the isomorphism classes of extensions form a non-empty torsor for $ D^1(B/A; J) $. 
\item The automorphism group of any extension is $D^0(B/A;J) \cong \Hom_B(\Omega_{B/A}, J)$. (This is elementary; the theory of the cotangent complex is not relevant.)
\end{itemize}
\end{prop}

\subsubsection{Application: Witt vectors}\label{sssec: Witt vectors} We will give an application of the preceding theory to the construction of \emph{Witt vectors} for perfect $\F_p$-algebras.

\begin{lemma}
If $B$ is a perfect $\F_p$-algebra, then $\bL_{B/\F_p} = 0$.
\end{lemma}

\begin{proof} The Frobenius endomorphism on $B$ is an isomorphism by definition, hence induces an automorphism of $\bL_{B/\F_p}$. To compute the action of Frobenius on $\bL_{B/\F_p}$, pick a free resolution of $\cB \xrightarrow{\sim} B$ over $\F_p$, and note that $\Frob$ on $B$ lifts to $\Frob$ on $\cB$ (the $p$-power map level-wise on $\cB$):
\[
\begin{tikzcd}
\cB \ar[r, "\Frob"] \ar[d, "\sim"] & \cB \ar[d, "\sim"] \\
B \ar[r, "\Frob"] & B 
\end{tikzcd}
\]
Hence the action of $\Frob$ on $\bL_{B/\F_p}$ is tensored from its action on $\bL_{\cB/\F_p}$, which is zero because the Frobenius endomorphism induces $0$ on K\"{a}hler differentials. So on the one hand we have that $\Frob$ induces an automorphism of $\bL_{B/\F_p}$, but on the other hand it is the zero map. Therefore, $\bL_{B/\F_p} =0$. 
\end{proof}

\begin{lemma}There is a unique flat deformation of $B$ to $\Z/p^n$. 
\end{lemma}

\begin{proof}The base case $n=1$ is tautological. We may suppose by induction that we have a unique deformation $B_n$ over $A_n := \Z/p^n$, and we want to extend it over $\wt{A} = \Z/p^{n+1}$. Let $I_n \subset \wt{A}$ be the kernel of $\wt{A} \rightarrow A_n$, a module over $A_n$ which in fact is pulled back from $A_1$. Let $J_n = I_n \otimes_{A_n} B_n$, which is then pulled back from $B$. Since the action of $B_n$ on $J_n$ factors over $B_1 = B$, and $B_n$ is flat over $A_n$, we have $D^i(B_n/A_n; J_n) \cong D^i(B/A_1;  J_n) \cong 0$ for each $i$. In particular, the obstruction class in Proposition \ref{prop: deformation-obstruction} lies in $D^2(B_n/A_n; J_n) = 0$, so an extension exists; and then since $D^1(B_n/A_n; J_n) = 0$, the extension is unique.
\end{proof}

Therefore, there is a unique (up to unique isomorphism) diagram 
\[
\begin{tikzcd}
B & B_2 \ar[l, twoheadrightarrow]  & B_3 \ar[l, twoheadrightarrow]  & \ldots  \ar[l, twoheadrightarrow]  \\ 
\F_p  \ar[u] & \Z/p^2 \Z  \ar[l, twoheadrightarrow] \ar[u]  & \Z/p^3 \Z  \ar[l, twoheadrightarrow] \ar[u]   & \ldots \ar[l, twoheadrightarrow] \ar[u] 
\end{tikzcd}
\]
with all squares co-Cartesian and all vertical arrows flat. The inverse limit $\varprojlim_n B_n$ is called the \emph{ring of Witt vectors of $B$}. This construction of the Witt vectors was pointed out in \cite[Remark 5.14]{Sc12}.

\subsection{Deformation theory of maps}

Our next situation is motivated by the deformation theory of maps: given a square-zero extension $X \inj \wt{X}$, and a map $X \rightarrow Y$, we want to understand extensions of the map 
\[
\begin{tikzcd}
X \ar[r, hook] \ar[d] & \wt{X} \ar[dl, dashed] \\
Y
\end{tikzcd}
\]
The local problem is: given an $A$-algebra $B$ and a square-zero extension 
\[
J \rightarrow \wt{B} \rightarrow B ,
\]
equip $\wt{B}$ with a compatible $A$-algebra structure. The equivalence class of $\wt{B}$ as a commutative $\Z$-algebra extension of $B$ by $J$ can be viewed as an element $[\wt{B}] \in D^1(B/\Z; J)$. We want to lift this a class in $D^1(B/A; J)$. We have the long exact sequence 
\[
\ldots \rightarrow D^0(A/\Z; J) \rightarrow D^1(B/A; J) \rightarrow D^1(B/\Z; J) \rightarrow D^1(A/\Z; J)  \rightarrow \ldots
\]
The map $D^1(B/\Z; J) \rightarrow D^1(A/\Z; J) $ can be interpreted as sending $\wt{B}$ to the fibered product $A \times_B \wt{B}$, which is a commutative algebra extension of $A$ by $J$. We can then interpret the long exact sequence as follows.

\begin{prop}\label{prop: map-obstruction} Let the situation be as above. 
\begin{itemize}
\item Attached to the extension $\wt{B} \rightarrow B$ is an \emph{obstruction class} $\mrm{obs}(\wt{B}) \in D^1(A/\Z; J)$, which vanishes if and only if there exists an $A$-algebra structure on $\wt{B}$ compatible with the given one on $B$. 
\item If $\mrm{obs}(\wt{B}) \in D^1(A/\Z; J)$ vanishes, then the isomorphism classes of $A$-algebra structures form a non-empty torsor for $ D^0(A/\Z; J)  = \Hom_A(\Omega_{A/\Z}, J)$. (This is elementary; the theory of the cotangent complex is not relevant.)
\end{itemize}
\end{prop}

\subsubsection{Application: Fontaine's map}\label{sssec: Fontaine's map}
Let $R$ be a $p$-adically complete and $p$-torsion-free ring. Let $R^{\flat} := \varprojlim_{\Frob_p} (R/p)$. This is a perfect ring, characterized uniquely up to unique isomorphism as ``the'' final perfect ring mapping to $R/p$. 

It is a very important fact in $p$-adic Hodge Theory that there is a \emph{unique} map $W(R^{\flat}) \rightarrow R$ that reduces mod $p$ to the canonical map $R^{\flat} \rightarrow R/p$, called \emph{Fontaine's map}, although we will not be able to explain its significance here. Here $W(R^{\flat})$ are the Witt vectors, which as discussed in \S \ref{sssec: Witt vectors} is the inverse limit of the unique (up to unique isomorphism) family of lifts
\[
\begin{tikzcd}
R^{\flat}  & W_2(R^{\flat}) \ar[l, twoheadrightarrow] &  W_3(R^{\flat})  \ar[l, twoheadrightarrow] & \ldots  \ar[l, twoheadrightarrow]  \\ 
\F_p \ar[u]  & \Z/p^2 \Z  \ar[l, twoheadrightarrow] \ar[u] & \Z/p^3 \Z  \ar[l, twoheadrightarrow] \ar[u]   & \ldots \ar[l, twoheadrightarrow] \ar[u] 
\end{tikzcd}
\]
Let us build up a family of maps step-by-step, starting with $W_2(R^{\flat})  \rightarrow R/p^2$.
\[
\begin{tikzcd}
R/p & R/p^2 \ar[l, twoheadrightarrow] & R/p^3 \ar[l, twoheadrightarrow] & \ldots \ar[l, twoheadrightarrow]  \\ 
R^{\flat} \ar[u]  & W_2(R^{\flat}) \ar[u, dashed]  \ar[l, twoheadrightarrow] & W_3(R^{\flat}) \ar[u, dashed] \ar[l, twoheadrightarrow] & \ldots  \ar[l, twoheadrightarrow]  \ar[u] \\
\F_p \ar[u] & \Z/p^2 \Z \ar[l, twoheadrightarrow] \ar[u]  & \Z/p^3 \Z  \ar[l, twoheadrightarrow] \ar[u] &  \ldots  \ar[l, twoheadrightarrow]  \ar[u] 
\end{tikzcd}
\]
(The $p$-torsionfree property of $R$ is used to see that $R/p^n$ is flat over $\Z/p^n \Z$.) By similar arguments as above, the obstruction to such an extension is a class in 
\[
D^1(W_2(R^{\flat}) /  (\Z/p^2 \Z); pW_2(R^{\flat}) ),
\]
which is isomorphic to $D^1(R^{\flat} / \F_p; R^{\flat})$, which vanishes because $R^{\flat}$ is perfect. Then, the set of such an extensions is a torsor for $D^0(R^{\flat} / \F_p; R^{\flat}) = 0$, so the extension is unique. By induction, we have constructed compatible maps $W_n(R^\flat) \rightarrow R/p^n$ for each $n \geq 1$. Finally, take the inverse limit of these maps to construct Fontaine's map (here is where we use that $R$ is $p$-adically complete).

\subsection{Global deformations}

So far we have discussed deformations of \emph{affine} schemes. Now we begin the study of deformations of possibly non-affine schemes.  

\begin{example}
Let $X$ be a \emph{smooth} (but not necessarily affine) variety over a field $k$. A \emph{first-order deformation} of $X$ is a flat family $\wt{X}$ over $k[\epsilon]/(\epsilon^2)$, whose fiber over $\Spec k$ is $X$. Then we claim that the set of first-order deformations of $X$ is in bijection with $H^1(X, TX)$. 

To see this, write $X = \bigcup U_i$ as a union of affines $U_i = \Spec A_i$. A first-order deformation of $X$ is then uniquely specified by a compatible collection of first-order deformations of the $U_i$. Since the $U_i$ are \emph{smooth}, they have unique first-order deformations $\wt{U}_i$. The gluing data is an automorphism of $\wt{U}_i \cap \wt{U}_j := \Spec \wt{A}_{ij}$ lifting the given automorphism on $U_i \cap U_j := \Spec A_{ij}$, which is of the form 
\[
a + b \epsilon \mapsto a + (\delta a + b) \epsilon
\]
where $\delta$ is a $k$-linear derivation from $A_{ij}$ to itself, i.e., a section of the tangent sheaf $TX$ over $U_{ij}$. In order to glue to a first-order deformation of $X$, these sections must satisfy the cocycle condition on triple intersections. Finally, modifying each $\wt{U}_i$ by an automorphism (as a first-order deformation) corresponds to adding a coboundary. Therefore, we see that the data of a first-order deformation is exactly that of a class in the first Cech cohomology $H^1(X, TX)$. 

As an application, consider the \emph{moduli space of (smooth, projective) genus $g$ curves} $\cM_g$. A $k$-point to this moduli space is a smooth, projective genus $g$ curve $C$. The tangent space to $\cM_g$ at the $k$-point corresponding to $C$ has underlying set the set of extensions of this $k$-point to a $k[\epsilon]/\epsilon^2$-point, which is to say the set of first-order deformations of $C/k$. According to what we have said, this is $H^1(C, TC)$. Furthermore, the obstructions are valued in $H^2(C, TC)$, which vanishes because $C$ is a curve -- this tells us that $\cM_g$ is \emph{smooth} at the $k$-point corresponding to $C$. 
\end{example}

\subsection{Globalizing the cotangent complex}\label{ssec: globalization of cotangent complex}

By globalizing Proposition \ref{prop: deformation-obstruction}, one proves:

\begin{prop}\label{prop: deformation-obstruction global} Let $f \co X \rightarrow S$ be a map of schemes. Let $S \inj \wt{S}$ be a square-zero thickening with ideal sheaf $\cI$, and $\cJ = f^* \cI$. 
\begin{itemize}
\item Attached to the extension $S \inj \wt{S}$ is an \emph{obstruction class} $\mrm{obs}(\wt{S}) \in \Ext^2(\bL_{X/S}, \cJ)$, which vanishes if and only if a flat extension $\wt{X} \rightarrow \wt{S}$ exists. 
\item If $\mrm{obs}(\wt{S}) \in \Ext^2(\bL_{X/S}, \cJ)$ vanishes, then the isomorphism classes of flat extensions form a non-empty torsor for $\Ext^1(\bL_{X/S}, \cJ) $. 
\item The automorphism group of any flat extension is $\Ext^0(\bL_{X/S}, \cJ)$. 
\end{itemize}
\end{prop}

\begin{example}[LCI curves]
Recall that a finite type morphism of Noetherian rings $A \rightarrow B$ is said to be \emph{local complete intersection (LCI)} if it can be Zariski-locally written as a composition $A \rightarrow A[x_1, \ldots, x_n]  \surj B$ where the second map is a regular quotient. More generally, a morphism of schemes is said to be LCI if it is Zariski-locally LCI. For example, a curve with at worst nodal singularities over a field is LCI. By Proposition \ref{prop: LCI cotangent complex}, if $Y \rightarrow X$ is LCI then $\bL_{Y/X}$ is concentrated in degrees $-1,0$. 

Suppose $C$ is an LCI curve over a field $k$, which is is smooth outside a finite number of closed points of $C$. Then we claim that $\Ext^2(\bL_{C/k}, \cJ) = 0$, which in particular shows that all infinitesimal deformations of $C/k$ are unobstructed. By the local-global spectral sequence, it suffices to show the vanishing of 
\begin{itemize}
\item $H^2(C, \ul{\Ext}^0(\bL_{C/k}, \cJ))$,
\item $H^1(C, \ul{\Ext}^1(\bL_{C/k}, \cJ))$, and
\item $H^0(C, \ul{\Ext}^2(\bL_{C/k}, \cJ))$. 
\end{itemize}
The first vanishes because $C$ is a curve, so it has cohomological dimension $1$. The third vanishes because $\ul{\Ext}^2(\bL_{C/k}, \cJ)$ vanishes, thanks to $\bL_{C/k}$ being supported in degrees $[-1, 0]$. Finally, since $C$ is smooth away from a finite collection of closed points, $\ul{\Ext}^1(\bL_{C/k}, \cJ)$ is a torsion sheaf by the localizing property of the cotangent complex plus its calculation in the smooth case (\S \ref{ssec: properties of cotangent complex}), so its higher cohomology vanishes. 

\end{example}

For completeness, we state the globalization of Proposition \ref{prop: map-obstruction}: 

\begin{prop}\label{prop: map-obstruction global} Let $X \inj \wt{X}$ be a square-zero extension with ideal sheaf $\cJ$ and $f \co X \rightarrow Y$. Let $\bL_Y$ be the absolute cotangent complex of $Y$ (i.e., the cotangent complex of $Y \rightarrow \Spec \Z$). 
\begin{itemize}
\item Attached to $\wt{X}$ is an \emph{obstruction class} $\mrm{obs}(\wt{X}) \in \Ext^1(f^* \bL_{Y}, \cJ)$, which vanishes if and only if $f$ extends to a map $\wt{X} \rightarrow Y$. 
\item If $\mrm{obs}(\wt{X}) \in \Ext^1(f^* \bL_{Y}, \cJ)$ vanishes, then the isomorphism classes of extensions of $f$ form a non-empty torsor for $\Ext^0(f^* \bL_{Y}, \cJ)$. 
\end{itemize}
\end{prop}

\subsection{Geometric interpretations of Andr\'{e}-Quillen homology}\label{ssec: geometric interpretation} In this section we reflect on the relation between derived algebraic geometry and the cotangent complex; a more substantial exposition of this topic may be found in \cite{V10}. As was discussed in \S \ref{sec: derived moduli}, one of the characteristic features of derived algebraic geometry in practice is its ability to construct spaces whose cotangent complexes can be ``identified'' in a meaningful way, for example as some natural cohomology theory. The reason this happens is because derived algebraic geometry supplies a ``geometric'' interpretation of the cotangent complex. 

It may be helpful to consider the analogy between schemes and reduced schemes, explained in \S \ref{ssec: visualizing derived schemes}. Even if we are only interested in reduced schemes, allowing nonreduced schemes into our vocabulary allows us to gives a geometric interpretation of tangent spaces as maps from dual numbers. In the context of moduli spaces, this allows to express tangent spaces as first-order deformations, and then often in practice as a cohomology group. 

We could contemplate defining a ``reduced moduli space'' as a functor on reduced commutative rings. However, it would be a priori unclear how to describe its tangent space explicitly. If the functor can be extended to non-reduced rings, then its tangent space could be calculated as in the previous paragraph, and if for example the moduli space turned out to be smooth then one would know a posteriori the tangent space of the underlying reduced scheme. The reader might enjoy contemplating this thought experiment in the example of $\cM_g$, the moduli space of smooth projective genus $g$ curves. 

Similarly, the tangent \emph{complex} of a classical moduli space is a priori hard to describe explicitly, because the ``higher'' tangent spaces are not prescribed by some direct process from the moduli problem. However, in the derived world the full tangent complex does admit such a direct description, which we now explain. 

Let $A \rightarrow B$ be a morphism of commutative rings, and $\eta \co B \rightarrow A$ an augmentation. It may be helpful to focus on the special case $A = \Spec k$ for a field $k$, in which case $\eta$ corresponds to a $k$-point of $B$. 

Let $M$ be a $B$-module. It can be turned into a simplicial $B$-module $M[0]$, which is characterized as the unique simplicial $B$-module with $\pi_0(M[0]) = M$ and $\pi_i(M[0]) = 0 $ for $i>0$. For $n \geq 0$, there is a simplicial $B$-module $M[n]$ with $\pi_n(M[n]) = M$ and $\pi_i(M[n]) = 0 $ for $i \neq n$. (For example, $M[n]$ can be constructed using the Dold-Kan correspondence.) 

Then we form the level-wise simplicial square-zero extension $A \oplus M[n]$. Its graded homotopy ring is the square-zero extension of $A$ (in degree $0$) by $M$ (in degree $n$). 

\begin{example}
Let $A = k$ be a field and $M=k$. Then $k \oplus k[0] = k[\epsilon]/\epsilon^2$ is the ring of \emph{dual numbers} over $k$. More generally, $k \oplus k[n]$ is called the ring of \emph{derived order $n$ dual numbers}. 
\end{example}

Then the Andr\'{e}-Quillen cohomology group $D^n(B/A; M)$ can be viewed as the group of homomorphisms from $B$ to $A \oplus M[n]$ lifting $\eta$, in the homotopy category of simplicial $A$-algebras. The reason for this is relatively formal, once the relevant definitions are in place; for now we just give a sketch. To compute these homomorphisms, by definition one needs to resolve $B$ by a free simplicial $A$-algebra $\cB$. (This is analogous to how maps in the derived category of $A$-modules are computed by first resolving the source by a projective resolution.) Then, lifts 
\[
\begin{tikzcd}
& A \oplus M[n] \ar[d] \\
\cB \ar[ur, dashed] \ar[r] & A
\end{tikzcd}
\]
are computed by derivations from $\cB$ into $M[n]$, which can be converted into a description as homomorphisms from K\"{a}hler differentials of $\cB$ into $M[n]$, which is basically the definition of Andr\'{e}-Quillen cohomology. 

Geometrically, we can think of ``$\Spec (A \oplus M[n])$'' as being some ``higher derived infinitesimal thickening'' of $A$ by $M$. A derived moduli problem can be evaluated on such a space, giving a geometric interpretation of Andr\'{e}-Quillen homology groups; this makes it \emph{easier} to understand the cotangent complex of a derived moduli space than that of a classical moduli space, in the same sense that it is ``easier'' to understand the tangent space of a classical moduli space than that of a ``reduced moduli space''.

\bibliographystyle{amsplain}
\providecommand{\bysame}{\leavevmode\hbox to3em{\hrulefill}\thinspace}
\providecommand{\MR}{\relax\ifhmode\unskip\space\fi MR }
\providecommand{\MRhref}[2]{%
  \href{http://www.ams.org/mathscinet-getitem?mr=#1}{#2}
}
\providecommand{\href}[2]{#2}

\end{document}